\documentclass{article}
\usepackage{chngcntr}
\counterwithin{equation}{section}

\usepackage{amsmath, amsthm, amssymb, amsfonts}
\usepackage{thmtools}
\usepackage{graphicx}
\usepackage{setspace}
\usepackage{geometry}
\usepackage{float}
\usepackage{hyperref}
\usepackage[utf8]{inputenc}
\usepackage[english]{babel}
\usepackage{framed}
\usepackage[dvipsnames]{xcolor}
\usepackage{tcolorbox}
\usepackage{color}
\usepackage{tikz-cd}
\usepackage[title]{appendix}
\usepackage{biblatex}
\newtheorem{theorem}{Theorem}[section] 

\newtheorem{definition}[theorem]{Definition} 

\newtheorem{lemma}[theorem] {Lemma}
\newtheorem{corollary}[theorem]{Corollary}
\newtheorem{proposition}[theorem]{Proposition}

\theoremstyle{remark}
\newtheorem{remark}[theorem]{\bf Remark}
\newtheorem{example}[theorem]{\bf Example}

\newcommand{\C}{\mathbb{C}}
\newcommand{\K}{\mathbb{K}}
\newcommand{\R}{\mathbb{R}}
\newcommand{\Z}{\mathbb{Z}}
\newcommand{\bp}{\bar{\partial}}
\begin{document}
\setcounter{section}{-1}

\author{Jian Han}

\thispagestyle{empty}
\addvspace{5mm}  


\begin{center}
\begin{doublespace}
{\textbf{\large A Construction of Formal Frobenius Manifold from Deformation of Complex Structure}}
\end{doublespace}

\vspace{10mm}
{Master Thesis by: Jian Han} 

\vspace{30mm}

{ Advisor: Prof. Noemie Combe}\\[12pt]
{ Examiner: Prof. Matthias Schwarz}\\[12pt]
{ For the Degree of}\\[12pt]
{Masters of Mathematical Physics} \vfill
{Universität Leipzig }\\
{Leipzig, Germany}
\vfill

\begin{onehalfspace}
{ February, 2024}

\end{onehalfspace}

\end{center}
\newpage
\thispagestyle{empty}

\tableofcontents

\newpage 
\section{Acknowledgements}

 I would like to thank my supervisor Prof. Noemie Combe, for her support and encouragement throughout the duration of my master thesis.  Her expertise and patience are instrumental not only in the completion of this thesis but also in my academic growth. 

\medskip
 
 \noindent I am also grateful to Prof. Matthias Schwarz for the generous assistance, insightful feedback and amazing courses he gave. 

\newpage
\section{Introduction}
\subsection{Motivation}
Quantum theory has revealed a new world in mathematics since twentieth century.  Frobenius manifolds were initially discovered in the context of physics in \cite{Dubrovin} for describing two dimensional topological quantum field theory. It appears in various fields in mathematics and physics, see for instance \cite{Manin,BK,Katz,Lurie2,CoMa, Saito}. This makes Frobenius manifold a potential language to describe diverse phenomenons in mathematics and physics such as  mirror symmetry.    

\smallskip

In the following we list different sources of Frobenius manifolds:

\begin{enumerate}
    \item Saito's theory of unfolding of isolated singularities, which is related to Landau-Ginzburg models. See \cite{Saito}
    
    \item Quantum cohomology.  A quantum cohomology ring is  a deformation of a cohomology ring where the deformed cup product (or quantum cup product) is defined by Gromov-Witten invariants. This corresponds to the $A$ twisted sigma model, see \cite{Manin}.

    \item  Formal moduli spaces of solutions to the Maurer–Cartan equations modulo gauge equivalence, which corresponds to the $B$ twisted sigma model, see \cite{BK}.
    
    \item (Flat) manifolds of probability distributions of exponential type. This comes from information geometry and has applications in machine learning, see \cite{CoMa,GSI}. 
    
\end{enumerate}


In this thesis we study (graded) Frobenius manifolds coming from (extended) deformations of complex structures on compact Calabi--Yau manifolds. This relies on investigations of the formal moduli spaces of solutions to the Maurer–Cartan equations modulo gauge equivalence. 

\smallskip




\subsection{Contents of this Thesis}
In this thesis we will introduce some background on complex structures and deformation theory. Then,  following  \cite{BK}, we expose the proof of the construction of a formal Frobenius manifold coming from (extended) deformations of complex structures, but starting from the differential geometry context. This corresponds to the source of dGBV algebras, which can be interpretated as a formal deformation problem in context of derived algbraic geometry.

Most of results in this thesis are not original. 

\smallskip

$\bullet$ In chapter 2 we introduce $G$ structures on a manifold and the integrability condition. In particular, we give the definition of an almost complex structure.

\medskip

$\bullet$ In chapter 3, we provide the definition of a deformation of a complex structure in terms of deformations of $G$ structures such that the integrability condition remains satisfied during the deformation. This relies on \cite{D.Huybrechts}.

\medskip

The central idea in deformation theory is that:

\smallskip

{\it ``If $X$ is a moduli space over a field $k$ of characteristic zero, then a formal neighborhood of any point
$x \in X$ is controlled by a differential graded Lie algebra
".}  

\smallskip

This was developed by Deligne, Drinfeld, and Feigi (see \cite{Lurie1}, for more).

\medskip

$\bullet$ Chapter 4 follows the philosophy outlined above. So we introduce differential graded Lie algebra  and explain how this construction is used in the deformation problems. This is applies to the example of deformations of complex structures on Calabi--Yau manifolds. 

We discuss the construction of the deformation functor and of its smoothness and suggest an equivalent definition of smoothness of the deformation functor.  This allows a practical  approach to the problem starting from differential geometry. 

We prove that an isomorphism between two deformations of a complex structure is the same as a gauge equivalence of their representatives in the corresponding differential graded Lie algebra. 

Finally, we apply the above construction to the case of Calabi-Yau manifolds. The unobstructness of a Calabi-Yau manifold can be considered as a special case of the extended deformation problem in \cite{BK}.

\smallskip

Chapter 5 presents a graded version of the deformation functor. 

For the ungraded case, a classical result is that an abelian DGLA is representable. We show that the same holds for the graded version.

This construction is used to define the extended moduli space of complex structures on a Calabi-Yau manifold. This is done following \cite{BK,Manin}.  

\smallskip

Chapter 6 mainly recalls notions on (formal) Frobenius manifolds. We discuss the proof that the extended moduli space is in fact a formal Frobenius manifold. We will  use methods in \cite{Manin}. 

\newpage

\section{$G$ structure on Manifold}
We introduce the notion of a $G$ structure  and its integrability in this section.  The main reference the author used in this section is \cite{S.Sternberg} and \cite{Kobayashi1}. 

\subsection{Definition of $G$ structure on Manifold}

We firstly recall the definition of principle bundle. 

Let $M$ be an $n$ dimensional smooth manifold. Given a Lie group $G$, a principle $G$ bundle  is a triple $(P,M,\pi)$ where is $(P,M,\pi)$ a fiber bundle and $G$ is a right action on $P$  such that  $G$ preserves the fibers and acts on each fiber freely and transitively. Hence any fiber of $P$ is diffeomorphic to the Lie group $G$. For simplicity sometimes we denote a principle $G$ bundle $(P,M,\pi)$ by $\pi: P \longrightarrow M$, or  just by $P$ if there is no confusion. 

\medskip 

Let $\pi_1:P_1\longrightarrow M_1$ and $\pi_2:P_2\longrightarrow M_2$ be two principle $G$ bundles. A morphism from $P_1$ to $P_2$ is a pair of maps $(\bar{f}:P_1\longrightarrow P_2,f:M_1 \longrightarrow M_2)$ such that the following  diagram commutes: 
$$\begin{tikzcd}
	{P_1} & {P_2} \\
	{M_1} & {M_2}
	\arrow["{\bar{f}}", from=1-1, to=1-2]
	\arrow["{\pi_1}"', from=1-1, to=2-1]
	\arrow["f"', from=2-1, to=2-2]
	\arrow["{\pi_2}", from=1-2, to=2-2]
\end{tikzcd}
$$
and $\Bar{f}(p \cdot g)=\Bar{f}(p)\cdot g$ for any $p\in P_1$ and $g\in G$. 

\smallskip

A subbundle of principle bundle $(P,M,\pi)$ is a principle $G'$ bundle $(P',M,\pi)$ where $G'\subset G$ is a Lie subgroup of $G$ and and there is an embedding  

$$\begin{tikzcd}
	P' & {P} \\
	M & M
	\arrow["\pi'"', from=1-1, to=2-1]
	\arrow["i", hook, from=1-1, to=1-2]
	\arrow["{\pi}", from=1-2, to=2-2]
	\arrow["id", from=2-1, to=2-2]
\end{tikzcd}$$
such that for any $p\in P'$, $g\in G' \subseteq G$, we have $i(pg)=i(p)g$.

\medskip

    For any point $x\in M$, on a smooth manifold $M$ denote $Fr_x$  the set of all frames of $T_xM$.  We have a principle $GL_n(\R)$ bundle $$\pi:\coprod_{x\in M}Fr_x \longrightarrow M,$$ where an element $g\in GL_n(\R)$ acts on any frame naturally by matrix multiplication. We call this principle bundle the \textbf{frame bundle} of $M$, denoted by $Fr(M)$. A section of a frame bundle on an open set $U\subset M$ is called a \textbf{frame field} on $U$.

\medskip

\begin{definition}
A\textbf{ $\mathbf{G}$ structure }on a smooth manifold $M$ is a principle $G$ bundle $\pi:P \longrightarrow M$ such that $(P,M,\pi)$ is a subbundle of the frame bundle $Fr(M)$.      
\end{definition}

We will call a frame $e_x$ of $T_x M$ which lies in $\pi^{-1}(x)$ an \textbf{admissible frame} of the $G$ structure.

\smallskip

Let   $\pi_1:P_1 \longrightarrow M_1, \pi_2:P_2\longrightarrow M_2$ be two $G$ structures.  If there is a diffeomorphism between the base manifolds $f: M_1 \longrightarrow M_2$,  $f$ induces an isomorphism between frame bundles $\bar{f_*}: Fr(M_1)\longrightarrow Fr(M_2)$ by mapping a frame $(e_1,...,e_n)_x$ at $x\in M_1$ to $(f_*e_1,...,f_* e_n)_{f(x)}$ at $f (x) \in M_2$. Without causing confusion we will write $\bar{f_*}$ also as $f_*$. If $f_*$ restricts to an isomorphsim between $P_1$ and $P_2$ we  call $f$ an isomorphism between the two $G$ structures. If $P_1,P_2$ are the same principle bundle, we will call such a $f$ an automorphism of $P_1$. 

\begin{definition}
   Given two $G$ structures  $\pi_1:P_1 \longrightarrow M_1, \pi_2:P_2\longrightarrow M_2$, and let $x\in M_1$, $y \in M_2$, we say $P_1$  is \textbf{locally equivalent } to  $P_2$ at $x$ if there exist a neighborhood $U$ of $x$ on $M_1$, a neighborhood of $V$ of $y$ on $M_2$ and an isomorphism $f: U \longrightarrow V$ such that $f$ induces an isomorphism between two $G$ structures restrict on $U,V$ respectively, i.e. ${P_1}|_U \cong P_2|_V$ .
\end{definition}

\smallskip

A Riemannian metric $g$ on a manifold $M$ is a positive definite inner product on tangent space at each point $p\in M$ and the inner product varies smoothly w.r.t $p$. 

\begin{proposition}
    Given a smooth manifold $M$ with dimension $n$, there is one to one  correspondence between Riemannian metrics on $M$ and  $O(n)$ structures on $M$, where $O(n)$ is the orthogonal  group in dimension $n$.
    \end{proposition}
\begin{proof}

    Given a Riemannian metric $g$ on $M$, at any $p\in M$, one can find the set the all orthogonal frames at $x$, denoted by $$S_x:=\{(e_1,...,e_n)\in Fr_x(M)|\, g(e_i,e_j)=\delta_{ij}\}\subset Fr_x(M),$$ where $\delta_{ij}$ is Kronecker delta. We define $P=\bigsqcup\limits_{x\in M} S_x $ with subset topology induced from $Fr(M)$ and projection $$\pi:P\longrightarrow M,\quad S_x \mapsto x.$$ Since $g$ varies smoothly, $P$ is a   submanifold of $Fr(M)$ and $\pi:P\longrightarrow M$ is a subbundle of $Fr(M)$ with local trivialization induced from the trivialization of $Fr(M)$. For any $x\in M$, $O(n)$ acts on $S_x$ freely and transitively. So $P$ is a principle $O(n)$ bundle, which is a subbundle of $Fr(M)$. Hence it is a $O(n)$ structure on $M$.

    Conversely, take an $O(n)$ structure $P$ on $M$. For any $x\in M$ and arbitrary admissible frame $(e_1,...,e_n)\in \pi^{-1}(x)\subset Fr(M),$ we define a positive definite inner product $g_x$ on it as $g_x(e_i,e_j)=\delta_{ij}$. This definition of inner product does not depend on the choice of admissible frame since for any other admissible frame $$(e_1',...,e_n')\in \pi^{-1}(x),$$ there exist a marix $[A_{ij}]\in O(n)$ such that $(e_1',...,e_n')=(e_1,...,e_n)A$ and $$g_x(e_i',e_j')=g_x(\sum\limits_k e_k A_{ik}, \sum\limits_l e_l A_{jl})=\delta_{ij}.$$ Hence $g_x$ is well-defined and a smooth choice of $g$ would be a Riemannian metric on $M$.
    \end{proof}

It is also easy to check that an isomorphism between two $O(n)$ structures is an isometry between the base manifolds. Hence knowing an $O(n)$ structure on a manifold is equivalent to knowing a Riemannian structure on a manifold.

\begin{definition}
    
    \label{ex:2.3}
    Given a smooth real manifold $M$ with dimension $2n$, an almost complex structure on $M$ is an endomorphism $J$ of tangent bundle $TM$ such that $J^2=-1$.  

    \end{definition}

   \begin{proposition}
  There is a one to one correspondence between almost complex structures on $2n$ dimensional manifold $M$ and $GL_n(\C)$ structures on $M$,  where $GL_n(\C)$ is identified with a subgroup of $GL_{2n}(\R)$ by the embedding $A+iB \mapsto \begin{bmatrix}
A & B \\
-B & A 
\end{bmatrix}$, where $A,B \in GL_n(\R)$.
      
   \end{proposition}
 \begin{proof}
     
    Given an almost complex structure $J$ on $M$,  locally on an open set $U\subset M$, we determine a frame field on $U$ by  firstly finding a vector field $v_1$ on $U$, then $J(v_1)$ is a vector field on $U$ which is independent with $v_1$  on $U$. We choose a vector field $v_2$ which is independent from the former two vector fields on $U$, then $J(v_3)$ is a vector field on $U$ which is independent from the first three vector fields on $U$, we can proceed this till we get  $v_n$ and  $J(v_n)$. In this way we find a desirable basis. Denote $J(v_i)$ by $w_i$, we define a frame field  $\{ v_1,...,v_n, w_1,...,w_n \}$ such that $$J(v_i)=w_i, \quad J(w_i)=-v_i.$$ We define a $GL_n(\C)$ structure on $M$ by restricting the fibres of $Fr(M)$ at every point $p$, precisely, the frames of $T_pM$ to the set of frames $\{ (v_1,...,v_n,w_1,...,w_n)g, g\in GL_n(\C)   \}$.

    Conversely, if we have a $GL_n(\C)$ structure on $M$, locally on  open set $U\subset M$, we choose an admissible frame field $(v_1,...,v_n,w_1,...,w_n)$ on $U$. We define an almost complex structure on $U$ by defining $J(v_i)=w_i, J(w_i)=-v_i$. The resulting almost complex structure is independent of the choice of the frame since if we choose another admissible frame field $(v'_1,...,v'_n,w'_1,...,w'_n)$ on $U$,  we have $$(v'_1,...,v'_n,w'_1,...,w'_n)=(v_1,...,v_n,w_1,...,w_n)\begin{bmatrix}
A & B \\
-B & A 
\end{bmatrix} $$ for some $A,B\in GL_n(\R)$. For simplicity we write $(v,w)$ instead of $(v_1,...,v_n,w_1,...,w_n)$, similar for $(v',w')$. Then we have  $$J(v',w')=J(v,w)\begin{bmatrix}
A & B \\
-B & A 
\end{bmatrix}=(w,-v) \begin{bmatrix}
A & B \\
-B & A 
\end{bmatrix}= (Aw+Bv,Bw-Av)=(w',-v').$$ 
 \end{proof}   

\medskip

For any $G \subset GL_n(\mathbb{R})$, we can construct a natural $G$ structure on an $n$-dimensional vector space $V$, by seeing $V$ as a smooth manifold and identify the tangent space at any $v\in V$ with $V$, choosing frames comes from the standard basis of $V$ acted by $G$.  Precisely, 
\begin{definition}
    Let $\mathcal{V}$ be an $n$-dimensional $\K$-vector space ($\K=\R$ or $\C$) and $(e_1,...,e_n)$ be its basis. A flat $G$ structure on $\mathcal{V}$ is the subbundle $(P,\mathcal{V},\pi)$ of its frame bundle $Fr(\mathcal{V})$ where the fibre of $P$ at any point $x\in \mathcal{V}$ is $G(\{e_1,...,e_n\})$, that is, the fibre at any point is the orbit of a canonical frame at this point under the action of $G$. We denote a flat structure on $\mathcal{V}$ by $G(\mathcal{V})$.
\end{definition}

\begin{definition}
    A $G$ structure $P$ on $M$ is called \textbf{integrable}  if it is locally equivalent to flat $G$ structure. That is,  given vector space $\mathcal{V}$ with same dimension as $M$, for any $x\in M$, there exists a neighborhood $U$ of x and an open set $V$ in $\mathcal{V}$ such that the $G$ structure $P|_U$ is isomorphic to $G(\mathcal{V})|_{V}$ (we will denote it by $G(V)$) as a $G$ structure. 
\end{definition}

The following proposition gives a criterion when is a $G$ structure integrable. 

\begin{proposition} {\label{t2.5}}
    Let $\pi:P\longrightarrow M$ be a $G$ structure on $M$, $P$ is integrable if and only if there exists an atlas $\{(U_\alpha,x_\alpha)\}_{\alpha\in A}$ of $M$, where $x_\alpha=(x_1,...,x_n)$ is the coordinate function on $U_\alpha$, such that the canonical frame field $(\frac{\partial}{\partial x_1},...,\frac{\partial}{\partial x_n})$  on  $U_\alpha$ always lies in $P$.
\end{proposition}

\begin{proof}
 For any $x\in M$, there is an open set $U_\alpha \subset M$ such that $x\in U_\alpha$. If $P$ is integrable, $P|_{U_\alpha}$ is isomorphic to a flat $G$ structure $G(V)$, which is a restriction of flat structure of $G(\mathcal{V})$ on an open set $V\subset \mathcal{V}$, where $\mathcal{V}$ is an $n$-dimensional $\R$ vector space. This means that there exists a homeomorphism $f: {V} \longrightarrow U_\alpha, f(x_1,...,x_n)=x$ such that for any $g\in G$, the frame  $$\{\frac{\partial}{\partial x_1},...,\frac{\partial}{\partial x_n}  \} \cdot g \in G(V),$$ and $$f_* \{\frac{\partial}{\partial x_1},...,\frac{\partial}{\partial x_n}  \} \cdot g \subset P|_U.$$ Let $g=$id,  we see $f^{-1}$ as the coordinate function on $U_\alpha$, then $\{\frac{\partial}{\partial x_1},...,\frac{\partial}{\partial x_n}  \}$  is the canonical frame field on $U_\alpha$, which lies in $P$. 
\smallskip

 Conversely, if we have an altas such that on a chart $(U_\alpha, x_\alpha)$ of $M$ with coordinate  $x_\alpha=(x_1,...,x_n)$ and canonical frame field $\{\frac{\partial}{\partial x_1},...,\frac{\partial}{\partial x_n}  \} \in P$,  we have a homeomorphism $f:U_\alpha \longrightarrow V\subset \mathbb{R}^n \cong \mathcal{V}$. Hence we have a canonical isomorphism $\Bar{f}:Fr(U_\alpha)\longrightarrow Fr(V)$. We restrict $\Bar{f}$ on $P|_{U_\alpha}$, we have an isomorphism $\Bar{f}|_{U_\alpha}:P|_{U_\alpha} \longrightarrow G(V)\subset Fr(V)$. It remain to prove that $G(V)$ is a restriction of flat $G$ structure on $V$. This can be seen as a fact that we can identiy $\{\frac{\partial}{\partial x_1},...,\frac{\partial}{\partial x_n}  \} \cdot g \in P$ with  $\{\frac{\partial}{\partial x_1},...,\frac{\partial}{\partial x_n}  \} \cdot g \in Fr(V)$. Since for any $g\in G$, the former always lies in $P$, the later always lies in $G(V)$. Therefore by definition, $G(V)$ is restriction of a flat $G$ structure.
\end{proof}

\medskip

\subsection{Complex Manifold and Integrability of Almost Complex Structures }
\begin{definition}
   We call a Hausdorff connected space a \textbf{complex manifold} if it is locally homeomorphic to an open set  $V\subset\mathbb{C}^n$, that is, it admits an atlas of coordinate charts $\{(U_i,h_i)\}$,  where $U_i \subset M$,  $h_i : U_i \rightarrow V_i \subset \mathbb {C}^n$ is a homeomorphism and   for any pair $U_i$, $U_j$ with homeomorphism $h_i$, $h_j$ respectively, we have that either $U_i \cap U_j = \emptyset$, or $h_i \circ {h_j}^{-1}$ is holomorphic on $U_i \cap U_j$. We call such an atlas a complex structure.
\end{definition}

\begin{proposition}\label{P:4.2}
Given an $n$ dimensional complex manifold $M$, its local charts determine an almost complex structure on its underlying $2n$ dimensional smooth real manifold $M$, i.e. an almost complex manifold $(M,J)$. 
\end{proposition}

\begin{proof}
We interpret $M$ as a $2n$ dimensional real manifold $M$, locally, in the holomorphic chart $U\subset M$ with local complex coordinates, a point $(z_1,...,z_n)\in \mathbb{C}^n$ is identified with  \\ $p=(x_1,...,x_n,y_1,... ,y_n) \in \mathbb{R}^{2n}$ by $z_i=x_i+\imath y_i$. The tangent space at this point is generated by $\frac{\partial}{\partial x_1},...\frac{\partial}{\partial x_n},\frac{\partial}{\partial y_1},...,\frac{\partial}{\partial y_n} $. We define an almost complex structure $J$ on the real tangent bundle $T_\R M$ locally at any point $p=(x_1,...,x_n,y_1,...,y_n)\in U\subset M$. 

\begin{equation}\label{E:1}J: T_p U \longrightarrow T_p U, \quad J(\frac{\partial}{\partial x_i}) = \frac{\partial}{\partial y_i}, \quad J(\frac{\partial}{\partial y_i}) =-\frac{\partial}{\partial x_i},\end{equation}

Obviously  $J^2 = -\mathrm{id}$ hence $J$ is indeed an almost complex structure. 

The almost complex structure we just defined does not depend on the choice of chart. If we choose another chart $(z'_1,...,z'_n)$ for $U$,  then we have $(z_1,...,z_n)=h (z'_1,...,z'_n)$ for some invertible holomophic map $h$. Define $x'_i,y'_i$ by $z'_i=x'_i+iy'_i$, and write $h$ as $h=(h_1(x'_1,...,y'_n),...,h_{2n}(x'_1,...,y'_n))$, we have that 
\begin{align*} 
\frac{\partial}{\partial x'_i} &=\frac{\partial h_j}{\partial x'_i} \frac{\partial}{\partial x_j}+\frac{\partial h_{j+n}}{\partial x'_i} \frac{\partial}{\partial y_j}    \\
\frac{\partial}{\partial y'_i} &=\frac{\partial h_j}{\partial y'_i} \frac{\partial}{\partial x_j}+\frac{\partial h_{j+n}}{\partial y'_i} \frac{\partial}{\partial y_j}   
\end{align*}
So $J$ maps on $\frac{\partial}{\partial x'_i}, \frac{\partial}{\partial y'_i}$ as 

\begin{align*}
\frac{\partial}{\partial x'_i} &\mapsto \sum_j \frac{\partial h_j}{\partial x'_i} \frac{\partial}{\partial y_j}-\sum_j \frac{\partial h_{j+n}}{\partial x'_i} \frac{\partial}{\partial x_j}    \\
\frac{\partial}{\partial y'_i} &\mapsto \sum_j \frac{\partial h_j}{\partial y'_i} \frac{\partial}{\partial y_j}-\sum_j \frac{\partial h_{j+n}}{\partial y'_i} \frac{\partial}{\partial x_j}  .
\end{align*}

Since $h$ is holomorphic, we have Cauchy-Riemann equation 

\begin{align*}
    \frac{\partial h_j}{\partial x'_i} &= \frac{\partial h_{j+n}}{\partial y'_i} \\
    \frac{\partial h_j}{\partial y'_i} &= -\frac{\partial h_{j+n}}{\partial x'_i},
\end{align*}
substitute this into the last set of equations we have again $J:\frac{\partial}{\partial x'_i} \mapsto \frac{\partial}{\partial y'_i}, \frac{\partial}{\partial y'_i} \mapsto -\frac{\partial}{\partial x'_i} $.

\end{proof}

\medskip

 For a complex manifold $M$ and  any $p \in M$, define $T_{\C,p}(M)=T_{\mathbb{R},p}(M) \otimes \mathbb C$ to be the complexified tangent space to $M$ at $p$. It can be realized as the space of $\mathbb{C}$-linear derivations in the ring of complex-valued $C^\infty$ functions on $M$ around $p$. We  write 
 $$\begin{aligned}
T_{\mathbb{C}, p}(M) & =\mathbb{C}\{\frac{\partial}{\partial x_i}, \frac{\partial}{\partial y_i}\} \\
& \cong\mathbb{C}\left\{\frac{\partial}{\partial z_i}, \frac{\partial}{\partial \bar{z}_i}\right\}
\end{aligned}$$
where
$$
\frac{\partial}{\partial z_i}=\frac{1}{2}\left(\frac{\partial}{\partial x_i}-i \frac{\partial}{\partial y_i}\right), \quad \frac{\partial}{\partial \bar{z}_i}=\frac{1}{2}\left(\frac{\partial}{\partial x_i}+i \frac{\partial}{\partial y_i}\right)
$$
Extend its almost complex structure $J$ on $T_\mathbb{R} M$ defined in Proposition \ref{P:4.2} linearly to $T_\mathbb{C} M$. We see that
$$
\begin{aligned}
& J \frac{\partial}{\partial z_i}=i \frac{\partial}{\partial z_i}, \\
& J \frac{\partial}{\partial \bar{z}_i}=-i \frac{\partial}{\partial \bar{z}_i}
\end{aligned}
$$ 
So we have a direct sum decomposition of complexified tangent space of $M$,
$$
T_{\mathbb{C}} U=T^{1,0} U \oplus T^{0,1} U
$$
such that $T^{1,0} U$ and $T^{0,1} U$ are the eigenspace of (the $\mathbb C$-linear extension of) $J$ with eigenvalues $i$ and $-i$ respectively.
Tangent vectors in $T^{1,0}M$ (resp. $T^{0,1} M$) are called holomorphic (resp. anti-holomorphic) tangent vectors. The corresponding vector bundle is called an holomorphic (resp. anti-holomorphic) tangent bundle. (The settings of tangent space on complex manifold here can be found in detail in \cite{P.Griffiths})

\smallskip

We have similar setting on almost complex manifold. Let $(M,J)$ be an almost complex manifold, we define $T_{\mathbb{C},p}(M)=T_{\mathbb{R},p}(M) \otimes \mathbb C$ the complexified tangent space to $M$ at $p$. Then we have

\begin{theorem} \label{t:2.10}
    
    On an almost complex manifold $(M,J)$ there exists a direct sum decomposition
$$
T_{\mathbb{C}} M=T^{1,0} M \oplus T^{0,1} M
$$ 
of complex vector bundles on $M$, such that  $T^{1,0} M$ is the eigenspace of the $\mathbb{C}$-linear extension of $J$ with eigenvalue $i$ and   $T^{0,1} M$ is the eigenspace with eigenvalue $-i$.
\end{theorem}
 If $(M,J)$ is induced from $M$ as a complex manifold as in Proposition \ref{P:4.2}, then $T^{1,0} M$ is naturally isomorphic (as a complex vector bundle) to the holomorphic tangent bundle of $M$. We call $T^{1,0} M$ the holomorphic vector bundle of $(M,J)$ and $T^{0,1}M$ the anti-holomorphic tangent bundle of $(M,J)$. A section of (anti-)holomorphic tangent bundle a (anti-)holomorphic vector field.

\begin{proof}
  For any $x\in M$, $v\in  T_{\C,x}M$, $v=\frac{1}{2}(v-iJv)+\frac{1}{2}(v+iJv)$, one can check that this is a direct sum decomposition into  eigenspaces with eigenvalues $i,-i$, respectively.  The (anti-)holomorphic tangent space of complex manifold $M$ at any point $p\in M$ is by definition the eigenspace with eigenvalue $i$ (resp. $-i)$ of $J$, which coincide the definition of $T^{1,0}$ (resp. $T^{0,1}$).
\end{proof}
\begin{remark} \label{r:2.19}
   Conversely, an almost complex structure is also determined by the decomposition in  direct sum of mutually-conjugated vector subbundles of complexified tangent bundle with eigenvalues $i$ and $-i$. For any $v \in \Gamma(T_{\mathbb {C}}M)$, that is, $v$ is a vector field on an open set $U\subset M$. If it can be decomposed as $$v=av_1+bv_2,$$ where $v_1 \in T^{1,0}M, v_2 \in T^{0,1}M$, and $a,b$ are complex valued smooth functions on $M$, then we have 
   $$J v= iav_1-ibv_2.$$ Therefore, once we know the (anti-)holomorphic part of the tangent bundle, we know the almost complex structure. If $M$ is a complex manifold, we sometimes write it is as $(M,J)$ where $J$ is the almost complex structure induced from the complex manifold $M$.
\end{remark}

{Given $M$ a complex manifold with local coordinates $(z_i)$, consider it as an almost complex  manifold with local coordinates $(x_i,y_i)$ such that $z_i=x_i+iy_i$. We can identify its real tangent space  with the holomorphic part of the complexified tangent space of the complex manifold $M$ by, in locally coordinates $(z_i=x_i+iy_i)$, identifying $\frac{\partial}{\partial x_i}$ with $\frac{\partial}{\partial z}$, $\frac{\partial}{\partial y_i}$ with $i\frac{\partial}{\partial z}$. With this identification on $M$, a holomorphic
tangent field generates a one parameter group of diffeomorphisms as its flow and any one parameter group of diffeomorphisms has a generator as a holomorphic vector field. 

}
\medskip

\textbf{Notation}: From now , by $TM$, we mean the complexified tangent  bundle of the manifold $M$.
\medskip

\begin{proposition}
    
    Let $(M,J)$ be an almost complex manifold. It is integrable if and only if it is induced from a complex structure on $M$ in the way in Proposition \ref{P:4.2}. 
\end{proposition}
\begin{proof}
    Let $(M,J)$ be a $2n$ dimensional almost complex manifold, if it is  integrable, by Proposition \ref{t2.5}, locally on an open set $U\subset M$, we choose a chart $(x_1,...,x_n,y_1,...,y_n)$ such that $\{ \frac{\partial}{\partial x_1},...,\frac{\partial}{\partial x_n},\frac{\partial}{\partial y_1},...,\frac{\partial}{\partial y_n}  \}$ forms an admissible frame field on $U$. By Definition \ref{ex:2.3}, we have \begin{equation}{\label{eq:7}}      
    J\left(\frac{\partial}{\partial x_i}\right)=\frac{\partial}{\partial y_i}, \quad J\left(\frac{\partial}{\partial y_i}\right)=-\frac{\partial}{\partial x_i}.
    \end{equation}
For simplicity we write $(x,y)$ for $(x_1,...,x_n,y_1,...,y_n)$ and $( \frac{\partial}{\partial x}), (\frac{\partial}{\partial y})$ for $( \frac{\partial}{\partial x_1},...,\frac{\partial}{\partial x_n}),(\frac{\partial}{\partial y_1},...,\frac{\partial}{\partial y_n})  $ respectively,  i.e. we have defined a homeomorphism $h:U\longrightarrow V\subset \C^n$ and $h(x,y)=(x_1+iy_1,...,x_n+iy_n)$. To see that $M$ is a complex manifold, it is  sufficient to verify that if there is another chart $(x'_1,...,x'_n,y'_1,...,y'_n)$ with canonical admissible frame field on $U'$ such that $U' \cap U \neq \emptyset$, then the transition map on $U' \cap U $ is holomorphic, that is, on $U' \cap U $,  Cauchy-Riemann equations

\begin{equation} \label{eq:8}
    \begin{split}
        \frac{\partial x'_i}{\partial x_j}&=\frac{\partial y'_i}{\partial y_j}, \\
        \frac{\partial x'_i}{\partial y_j}&=-\frac{\partial y'_i}{\partial x_j}
    \end{split}
\end{equation}
are satisfied.

Again for simplicity we write $x',y', \frac{\partial}{\partial x'}, \frac{\partial}{\partial y'}$ as abbreviation of $(x'_1,...,x'_n),  (y'_1,...,y'_n)$,   $( \frac{\partial}{\partial x'_1},...,\frac{\partial}{\partial x'_n}), \\ (\frac{\partial}{\partial y'_1},...,\frac{\partial}{\partial y'_n})$ respectively. Since $( \frac{\partial}{\partial x'}, \frac{\partial}{\partial y'} )$ is another admissible frame field on $U\cap U'$, we have $( \frac{\partial}{\partial x}, \frac{\partial}{\partial y} )= ( \frac{\partial}{\partial x'}, \frac{\partial}{\partial y'} )\begin{bmatrix}
A & B \\
-B & A 
\end{bmatrix}$ for some $A,B \in GL_n(\R)$. Using this we can compute terms in Equation \ref{eq:8}  
\begin{align*}
\frac{\partial x'_i}{\partial x_j} &= \frac{\partial x'_i}{\partial x'_k}A^k_j - \frac{\partial x'_i}{\partial y'_k}B^k_j =A^i_j\\
\frac{\partial x'_i}{\partial y_j} &= \frac{\partial x'_i}{\partial x'_k}B^k_j + \frac{\partial x'_i}{\partial y'_k}A^k_j =B^i_j\\
\frac{\partial y'_i}{\partial y_j} &= \frac{\partial y'_i}{\partial x'_k}B^k_j + \frac{\partial y'_i}{\partial y'_k}A^k_j=A^i_j \\
\frac{\partial y'_i}{\partial x_j} &= \frac{\partial y'_i}{\partial x'_k}A^k_j - \frac{\partial y'_i}{\partial y'_k}B^k_j =-B^i_j
\end{align*}
    which implies that Equation \ref{eq:8} is satisfied.
\end{proof}

\begin{definition}
   The Nijenhuis tensor is a  $(2,1)$-tensor field on a smooth (real) manifold such that for any vector fields $X$  and $Y$ on $M$ we have: 
    $$
    N_J(X, Y)=[X, Y]+J([J X, Y]+[X, J Y])-[J X, J Y].
    $$
\end{definition}

\begin{theorem}\label{T: 2.14.} 
Let $(M,J)$ be an almost complex manifold. The Nijenhuis tensor vanishes  if and only if $[T^{0,1} M, T^{0,1} M] \subset T^{0,1} M$.     
\end{theorem} \label{t:4.7}
 \begin{proof}(See \cite{D.Huybrechts}, chapter 2.)
     If  $[T^{0,1} M, T^{0,1} M] \subset T^{1,0} M$, for any $v,w \in \Gamma(TM)$, they can be  decomposed as $v=v_1+v_2$, $w=w_1+w_2$, where $v_1, w_1 \in \Gamma (T^{1,0}M)$, $v_2, w_2 \in \Gamma (T^{0,1}M)$. Compute $N_J(v,w)$ componentwisely we get $N_J=0$.
     
     To prove the other direction, suppose we have $v,w \in \Gamma(T^{0,1}M)$. If $N_J(v,w)=0$, expanding $N_J$ by definition shows that $J[v,w]=-i[v,w]$, which means $[v,w]\subset T^{0,1}M$.
 \end{proof}

\begin{theorem} {(Newlander-Nirenberg)}\label{C:NN}
    Let $(M,J)$ be an almost complex manifold. 
    $(M,J)$ is integrable if and only if its Nijenhuis tensor vanishes.
\end{theorem}

\begin{lemma}
    
     Any two dimensional almost complex manifold is integrable. 
\end{lemma}
   
\begin{proof}
Let $S$ be a two dimensional manifold and $v \in \Gamma (TS)$ be a vector field on $S$. Since the Lie bracket is anti-commutative, we have
$$
N_J(v,v)=[v, v]+J([J v, v]+[X, J v])-[J v, J v]=0 $$ and $$
N_J(v,Jv)=[v, Jv]+J([J v, Jv]+[v, J^2 v])-[J v, J^2 v]=0 .
$$
Because our manifold is two dimensional and $v$, $Jv$ are linearly independent, they linearly span the space of vector fields on the manifold. So $N_J$ vanishes on the two dimensional manifold. By Lemma \ref{C:NN}, we conclude that any two dimensional almost complex manifold is integrable.
\end{proof}

We now define some notations of vector bundles on  (almost) complex manifold.

Given an almost complex manifold $M$, we define vector bundle of complex $k$ forms on $M$ as $T^{k} M^*= \wedge^k {T M}^*$, vector bundle of $(p,q)$ form on $M$ as $$ T^{p, q} M^*=\wedge^p\left(T^{1,0} M\right)^* \otimes \wedge^q\left(T^{0,1} M\right)^*.$$ Hence there is a canonical isomorphism $$ T^{k} M^*=\wedge^k({T^{1,0}M \oplus T^{0,1}M})^* =\bigoplus_{p+q=k} T^{p, q} M^*.$$

We write the space of smooth sections of vector bundle $ T^{k}M^*$, $ T^{p, q} M^*$  as $\Omega^k (M)$, $\Omega^{p, q} (M)$ resp. and called an element in them a (complex)  $k$ form,  $(p,q)$ form respectively.

The exterior derivative $$d :\Omega^k M=\bigoplus\limits_{p+q=k} \Omega^{p, q} (M) \rightarrow \Omega^{k+1} M=\bigoplus \limits_{p+q=k+1} \Omega^{p, q} (M) $$can be decomposed. We denote the projection from $\bigoplus\limits_{i+j=k}\Omega^{i, j} (M)$ to $\Omega^{p, q} (M)$ by $\pi^{p,q}$. 

We  define
$$
\partial := \pi^{p+1,q} \circ d
$$
$$
\bar{\partial} := \pi^{p,q+1} \circ d .
$$
We call a complex form $\omega$ holomorphic if $\Bar{\partial}\omega=0$.

\begin{corollary} \label{c:2.25}
If the  almost complex manifold $(M,J)$ is integrable, then $\partial^2={\Bar{\partial}}^2=\partial \Bar{\partial}+\Bar{\partial}\partial=0$
\end{corollary}
\begin{proof}
    See \cite{D.Huybrechts}, chapter 2.
\end{proof}

Therefore, if $(M,J)$ is integrable (or equivalently $M$ is a complex manifold), we can define \textbf{Dolbeault cohomology} groups of $M$ as $$
H^{p,q}(M)= \frac{ker (\Bar{\partial}: \Omega^{p,q}(M)\longrightarrow \Omega^{p,q+1}(M))}{Im (\Bar{\partial} :\Omega^{p,q-1}(M)\longrightarrow \Omega^{p,q}(M))}.
$$
Later we will use differential form with value in holomorphic tangent bundle $\Omega^{*,*}(T^{1,0}M)$ and Dolbeault cohomology with values in holomorphic tangent bundle $H^{*,*}(T^{1,0}M)$, which is defined as 
$$
\Omega^{p,q}(T^{1,0}M)=\Omega^{p,q}(M)\otimes T^{1,0}M,
$$

$$
H^{p,q}(T^{1,0}M)=H^{p,q}(M)\otimes T^{1,0}M.
$$

\smallskip

\begin{definition}
    A complex vector bundle on a complex manifold $M$, $\pi:E \rightarrow M$ is called a \textbf{holomorphic vector bundle }if $E$ is a complex manifold and the projection $\pi$ is a holomorphic map. 
\end{definition}

We will be interested in a special class of complex manifolds which is called Kähler manifold. 

\begin{definition}
    A (integrable almost) complex manifold $(M,J)$ is called a \textbf{Kähler manifold} if it admits a compatible symplectic structure, that is, there exists a symplectic form $\theta\in \Omega^2(M)$ such that $\theta(\cdot, J\cdot)$ defines a Riemannian metric on $M$. 
\end{definition}

We need a technical lemma for  later use, 
\begin{lemma} {($\partial \Bar{\partial}$ lemma)}  \label{l:pp}
    If a complex manifold $M$ is Kähler, then on complex $\Omega^{*,*}(M)$, an element $\alpha \in \Omega^{*,*}(M)$ is $\partial$-exact $\iff$ $\alpha$ is $\Bar{\partial}$-exact $\iff$ $\alpha$ is $\partial \Bar{\partial}$-exact.
\end{lemma}
\begin{proof}
    See \cite{D.Huybrechts}.
\end{proof}

   \newpage 

\section{Deformation of Complex Structure}
In this thesis we consider a complex structure as an integrable almost complex structure. The deformation of a complex structure is a family of almost complex structures such that integrability is always satisfied in this family.  We refer to  \cite{D.Huybrechts}.

\subsection{One Parameter Deformation of Complex Structure}

In last chapter we see that when  $(M,J)$ is an integrable almost complex manifold, it is naturally a complex manifold. Conversely, when $M$ is a complex manifold, it induces canonically an almost complex structure on $M$ as a smooth manifold. Hence we will sometimes call an integrable almost complex manifold $(M,J)$ a complex manifold.

\begin{definition}\label{D:def}
Let $(M,J)$ be a complex manifold.
A  deformation of a complex structure $J$  is a smooth family of integrable almost complex structures $J_t$ over $(-\epsilon, \epsilon) \subset \R$ on $M$, where $t\in (-\epsilon,\epsilon)$ and $\epsilon$ is small enough. When $t=0$,  we have $J_0=J$.
\end{definition}

As stated in Remark \ref{r:2.19}, an almost complex structure is determined by its decomposition of complexified tangent bundle 

\begin{equation}\label{E:dec}
    T M=T^{1,0} M \oplus T^{0,1} M.
\end{equation}

Hence, a smooth family of almost complex structures $J_t$ with $J_0=J$ on $M$ is determined by the family of decomposition $$TM=T_t^{1,0}M \oplus T_t^{0,1}M$$ 
over $(-\epsilon,\epsilon)\subset \R$.
\smallskip

Once  $T_t^{1,0}M $ and $T_t^{0,1}M$ are given, we have 
$$J_t|_{T_t^{1,0}M}=i \cdot \text { id, }$$ and $$J_t|_{T_t^{0,1}M}=-i  \cdot \text { id }.$$

Denote the space of sections of complexified tangent bundle by  $\Gamma(TM)$. Any vector field $v \in \Gamma(TM)$ can be decomposed as follows $$v=av_1+bv_2,$$ where $v_1 \in T_t^{1,0}M$ and $v_2 \in T_t^{0,1}M$ and $a,b$ are complex valued smooth functions on $M$. 

If one has the decomposition as in Equation \ref{E:dec}, then $$J_t (v)= iav_1-ibv_2.$$

For sufficiently small $t$,  we define the canonical projection $\pi^{0,1}|_{T^{0,1}_t M}$ from $TM$ to $T^{0,1}M$ restricted to $T_t^{0,1}M$, and write $\pi^{1,0}$ for the projection from $TM$ to $T_t^{1,0}M$.

\begin{lemma}
    The projection $\pi^{0,1}|_{T^{0,1}_t}$ is an isomorphism from $T_t^{0,1}M$ to $T^{0,1}M$.
\end{lemma}
\begin{proof}
 For any point $p$ on $M$, $J_t$ is always a diagonalizable matrix with two eigenvalues $i $, $-i$. By   choosing a suitable basis, we can write $J$ as \[\begin{bmatrix}
i &  &  &  &  &\\
 & ...&  &  &  &\\
 &  & i &  &  &\\
 &  &  & -i &  &\\
 &  &  &  & ...&\\
 &  &  &  & & -i 
\end{bmatrix}. \] It is a diagonal matrix where the first $n$ diagonal elements are $i$; the other diagonal elements are $-i$. Since $J_t$ is diagonalizable, we write $J_t$ as $M_t J {M_t}^{-1}$, where $\{M_t \}$ is a smooth family of invertible matrices with $M_0=\text{id}$. Since ${T_p}^{0,1}M$ is the eigenspace of $J$ with eigenvalue $-i$, vectors in ${T_p}^{0,1}M$ are those column vectors whose first $n$ rows are 0 and $\pi^{0,1}|_{T^{0,1}_t M}=\begin{bmatrix}
    0 &0 \\
    0& {id}_{n\times n}
\end{bmatrix}$.\\

For arbitrary $p\in M$, denote anti-holomorphic tangent space of $J_t$ at point $p$ by $T_{t,p}^{0,1}M$. We omit $p$ in  $T_{t,p}^{0,1}M$ if there is no confusion. Take an arbitrary vector $v$ in anti-holomorphic tangent space of $J_t$ at point $p \in M$, i.e.  $v\in T_{t,p}^{0,1}M$,  $v$ is an eigenvector of $J_t$ with eigenvalue $-i$. Since $J_t=M_t J {M_t}^{-1}$, $M_t^{-1}v$ is an eigenvector of $J$ with eigenvalue $-i$,  {$M_t^{-1}v=\begin{bmatrix}     
\textbf{0} \\
\textbf{w}
\end{bmatrix}$.} Therefore we have $v=M_t\begin{bmatrix}     
\textbf{0} \\
\textbf{w}
\end{bmatrix}$, where $\textbf{w}$ is $(z_1,...,z_n)^\intercal$.

We discuss now the injectivity of $\pi^{0,1}|_{T^{0,1}_t M}$. For any $v\in T_{t,p}^{0,1}M $,  if $\pi^{0,1}|_{T^{0,1}_t}v=0$, then the lower $n$ entries of $v$ are zero, that is, $v=M_t \begin{bmatrix}     
\textbf{0} \\
\textbf{w}
\end{bmatrix}=
\begin{bmatrix}     
\textbf{z} \\
\textbf{0}
\end{bmatrix}
$. The lower right $(n\times n)$-block of $M_t$ (denoted by $M_t'$) is invertible  as $t=0$. (In fact one has $M_0'=id_{n\times n}$). Hence $M_t'$ remains invertible when $t$ is small. The equation $v=M_t \begin{bmatrix}     
\textbf{0} \\
\textbf{w}
\end{bmatrix}=
\begin{bmatrix}     
\textbf{z} \\
\textbf{0}
\end{bmatrix}
$ has a solution if $M_t' \textbf{w}=0$, which is only possible if $\textbf{w}=0 $ since $M_t'$ is invertible and $v=M_t \begin{bmatrix}     
\textbf{0} \\
\textbf{w}
\end{bmatrix}=0.$ 

To show the surjectivity of $\pi^{0,1}|_{T^{0,1}_t M}$, for any $\begin{bmatrix}
    \textbf{0} \\
\textbf{y}
\end{bmatrix}\in T_p^{0,1}M$,  we need to solve the equation $\pi^{0,1}|_{T^{0,1}_t M} \textbf{x}=\begin{bmatrix}
    \textbf{0} \\
\textbf{y}
\end{bmatrix}$, where $\textbf{x} \in T_{t}^{0,1}M $. It can be written as $M_t\begin{bmatrix}     
\textbf{0} \\
\textbf{a}
\end{bmatrix} $ as above. Thus we want to  solve the equation $$\begin{bmatrix}
    0 &0 \\
    0& {id}_{n\times n}
\end{bmatrix} M_t\begin{bmatrix}     
\textbf{0} \\
\textbf{a}
\end{bmatrix}=\begin{bmatrix}
    \textbf{0} \\
\textbf{y}
\end{bmatrix}.$$ 
Again, when $t$ is small, $M_t'$ is invertible. So this equation is solvable, i.e. $\pi^{0,1}|_{T^{0,1}_t M}$ is surjective.
\end{proof} 

\begin{remark}\label{compact}
    To make the choice of  sufficiently small $t$ consistent on $M$, we will assume from now on that \textbf{all manifolds are compact}.
\end{remark}

\begin{proposition}\label{p:3.4}
    For small $t$, a deformation of complex structure $J_t$ is encoded in the map $\phi_t$, which is defined as 
    \begin{equation} \label{d:encode}
        \phi_t:=\pi ^{1,0} \circ ({\pi^{0,1}|_{T^{0,1}_t M}})^{-1} : T^{0,1}M \longrightarrow T^{1,0}M.
    \end{equation}

\end{proposition}
\begin{proof}
At any point $p$ on the manifold $M$, for any $v_t \in T_{t,p}^{0,1}M$ 
    $$
    \begin{aligned}
v_t & =\pi^{1,0} \cdot v_t+\pi^{0,1}|_{T^{0,1}_t M} \cdot v_t \\
& =\pi^{1,0} \circ (\pi^{0,1}|_{T^{0,1}_t M})^{-1} \circ \pi^{0,1}|_{T^{0,1}_t M} \circ v_t+\pi^{0,1}|_{T^{0,1}_t M} \circ v_t \\
& =(\phi_t+\mathrm{id})(\pi^{0,1}|_{T^{0,1}_t M} \circ v_t) \\
&  \subset (\mathrm{id}+\phi_t) (T^{0,1}_{p}M).
\end{aligned}
$$
Hence  $T_{t,p}^{0,1}M\subset (\text{id}+\phi_t) (T_p^{0,1}M)$. Since $\pi^{0,1}|_{T^{0,1}_t M}$ is an isomorphism, for any  $v \in T_p^{0,1}M$, there exists a vector $v_t \in T_{t,p}^{0,1}M$ such that $v=\pi^{0,1}|_{T^{0,1}_t M} \cdot v_t$. Hence $$(\text{id}+\phi_t) \cdot v = (\mathrm{id}+\phi_t) \cdot \pi ^{0,1} \cdot v_t = (\pi^{0,1}|_{T^{0,1}_t M}+\pi^{1,0})v_t=v_t   \in T_{t,p}^{0,1}M. $$ Therefore we get $(\text{id}+\phi_t) (T_p^{0,1}M) \subset T_{t,p}^{0,1}M$. Then we have $T_{t,p}^{0,1}M = (\text{id}+\phi_t) (T_p^{0,1}M)$.

From Remark 2.13, the direct sum decomposition (Eq.\eqref {E:dec}) determines the almost complex structure. So given $\phi_t$ we know $T_{t}^{0,1}M$ and thus also the almost complex structure $J_t$.

Conversely, given $J_t$, which gives the decomposition $ T_{t}M=T_{t}^{1,0}M \oplus T_{t}^{0,1}M $, we have the  projection maps $\pi^{0,1}|_{T^{0,1}_t M}$ and $\pi^{1,0}$. So we get $\phi_t$ by definition. 

Therefore, knowing $J_t$ is equivalent to knowing $\phi_t$. In other words, the deformation of $J$ is encoded in $\phi_t$.  
\end{proof} 
\smallskip

\begin{definition}{\label{R:3.5}}
We call $\phi_t$ defined by Eq. \ref{d:encode} the \textbf{encoding map} of $J_t$. 
\end{definition}

    From the proof of last proposition, we know that given a deformation $J_t$, for any $v \in T^{0,1} M$, there is only one $w\in T^{1,0}M$ such that $v+w \in T_t^{0,1}M$, hence we can also define the encoding map $\phi_t$ as  $\phi_t$ sends $v$ to the only corresponding $w$ such that $v+w \in T_t^{0,1}M$.

\medskip

We want  to deform the complex structures on $M$ (i.e. the integrable almost complex structures on $M$). Now we use the integrability property of $J_t$ and express it via $\phi_t$ in Proposition.\ref{p:3.4}. To the end of this, we define the following  bilinear map
\begin{equation} \label{d:bracket}
\begin{split}
        [\cdot,\cdot]: \Omega^{0, 1}\left(T^{1,0}M\right) \times \Omega^{0, 1}(T^{1,0}M) &\longrightarrow \Omega^{0, 2}\left(T^{1,0}M\right) \\
 (\alpha,\beta)&\mapsto [\alpha,\beta]
 \end{split}
 \end{equation}

where in  local coordinates, 
\begin{equation} 
     \alpha=\sum_{i,j} a_i^j d \bar{z}^i \otimes \frac{\partial}{\partial z^j}  \text { and } \beta=\sum_{k,l} b_k^l d \bar{z}^k \otimes \frac{\partial}{\partial z^l}, 
\end{equation}

we have:
\begin{equation}\label{d:bracket}
    \begin{aligned}
{[\alpha, \beta] } & = \sum_{i,j,k,l}d \bar{z}^i \wedge d \bar{z}^k \otimes [a_i^j \frac{\partial}{\partial z^j}, b_k^l \frac{\partial}{\partial z^l}] \\
& = \sum_{i,j,k,l}d \bar{z}^i \wedge d \bar{z}^k \otimes\left(a_i^j \frac{\partial b_k^l}{\partial z^j} \frac{\partial}{\partial z^l}-b_k^l \frac{\partial a_i^j}{\partial z^l} \frac{\partial}{\partial z^j}\right) \in \Omega^{0, 2}\left(T^{1,0}M\right).
\end{aligned}
\end{equation}

We also introduce the partial differential 
\begin{equation}
    \bar{\partial}:\Omega^{0, 1}(T^{1,0}M) \longrightarrow \Omega^{0, 2}\left(T^{1,0}M\right).
    \end{equation}

In local coordinates, this is given by:
\begin{equation}\label{d:partial}
\bar{\partial} \alpha=\sum_{i,j,k} \frac{\partial a_i^j}{\partial \bar{z}^k} d \bar{z}^k \wedge d \bar{z}^i \otimes \frac{\partial}{\partial z_j}.
\end{equation}

\smallskip

\begin{lemma}
  Both operations $[\cdot,\cdot]$ and $\bp$  defined by Eq. \ref{d:bracket} and Eq.\ref{d:partial} resp. are independent of the choice of (holomorphic) local coordinates.
\end{lemma}
\begin{proof}
      For simplicity, we use the Einstein summation convention here. Let $$\alpha=\alpha_i^j d\bar{z}^i \otimes \frac{\partial}{\partial z^j}, \, \, \beta=\beta_k^l d\bar{z}^k \otimes \frac{\partial}{\partial z^l}$$ in coordinates $(z_1,...,z_n)$. By definition,  $[\alpha, \beta]_z=d\bar{z}^i \wedge d\bar{z}^k\otimes (a_i^j \frac{\partial b_k^l}{\partial z^j}\frac{\partial}{\partial z^l}-b_k^l\frac{\partial a_i^j}{\partial z^l}\frac{\partial}{\partial z^j})$.       If we work in another coordinates $(y_1,...,y_n)$, by chain rule we have  $$\alpha = a_i^j(\frac{\partial \bar{z}^i}{\partial y^p}dy^p+\frac{\partial \bar{z}^i}{\partial \bar{y}^p}d\bar{y}^p)\otimes (\frac{\partial y^q}{\partial z^j}\frac{\partial}{\partial y^q}+\frac{\partial \bar{y}^q}{\partial z^j}\frac{\partial}{\partial \bar{y}^q})$$, $$\beta= b_k^l(\frac{\partial \bar{z}^k}{\partial y^r}dy^r+\frac{\partial \bar{z}^k}{\partial \bar{y}^r}d\bar{y}^r)\otimes (\frac{\partial y^s}{\partial z^l}\frac{\partial}{\partial y^s}+\frac{\partial \bar{y}^s}{\partial z^l}\frac{\partial}{\partial \bar{y}^s}).$$ Since the transition functions are holomorphic,  we have $\frac{\partial \bar{z}^i}{\partial y^p}= \frac{\partial \bar{z}^k}{\partial y^r}=\frac{\partial \bar{y}^q}{\partial z^j}=\frac{\partial \bar{y}^s}{\partial z^l}=0$, and  $$[\alpha,\beta]_x=\left(\frac{\partial \bar{z}^i}{\partial \bar{y}^p}d\bar{y}^p \wedge \frac{\partial \bar{z}^k}{\partial \bar{y}^r}d\bar{y}^r \right)\otimes [a_i^j\frac{\partial y^q}{\partial z^j} \frac{\partial}{\partial y^q}, b_k^l \frac{\partial y^s}{\partial z^l} \frac{\partial}{\partial y^s}],$$ $$[\alpha,\beta]_y=\left(d\bar{y}^p \wedge d\bar{y}^r\right)\otimes [a_i^j \frac{\partial\bar{z}^i}{\partial \bar{y}^p}\frac{\partial y^q}{\partial z^j} \frac{\partial}{\partial y^q}, b_k^l \frac{\partial\bar{z}^k}{\partial \bar{y}^r}\frac{\partial y^s}{\partial z^l} \frac{\partial}{\partial y^s}]. $$ Again since transition fucntions are holomorphic,  $\partial (\frac{\partial\bar{z}^i}{\partial \bar{y}^p})/\partial y^s = \partial (\frac{\partial\bar{z}^k}{\partial \bar{y}^r})/\partial y^q=0$. Hence we can freely take the terms $\frac{\partial\bar{z}^i}{\partial \bar{y}^p},\frac{\partial\bar{z}^k}{\partial \bar{y}^r} $ out of Lie bracket, then it is  easy to see that $[\alpha, \beta]_z=[\alpha, \beta]_y$. 

The proof that $\bp$ is independent of the choice of coordinates is similar, using the holomorphicity of transition functions, we omit it here.

\end{proof}

By direct computation we can show that  the following is satisfied. 
\begin{lemma}\label{l:3.7}
    For any $\alpha,\beta \in \Omega^{0,1}(T^{1,0}M)$ , we have $$\bp[\alpha,\beta]=[\bp \alpha, \beta]-[\alpha,\bp \beta]$$. 
\end{lemma}

\begin{theorem}
    $J_t$ is integrable if and only if $\phi_t$ satisfies the Maurer--Cartan (MC) equation $$\Bar{\partial}\phi_t+ \frac{1}{2}[\phi_t,\phi_t]=0.$$
\end{theorem}
\begin{proof} (We follow closely the proof of \cite{D.Huybrechts}, Chapter 6, Lemma. 6.1.2.)
The proof is done in local coordinates. For simplicity we denote $\phi_t$ by $\phi$. 

By Theorem \ref{T: 2.14.}, $J_t$ is integrable if and only if $[T^{0,1}_tM,T^{0,1}_tM]\subset T^{0,1}_tM$. We use this result as the criterion of integrablity.  

Consider two vectors in the basis of $T^{0,1}M$:  $\frac{\partial}{\partial \Bar{z_i} }$, $\frac{\partial}{\partial \Bar{z_j} }$. From Prop.\ref{p:3.4} we know that $(\text{id}+\phi)\frac{\partial}{\partial \Bar{z_k} } \in T^{0,1}_t$.\\

If $[T^{0,1}_tM,T^{0,1}_tM]\subset T^{0,1}_tM$, we have 

$$
[\frac{\partial}{\partial \Bar{z_i} }+ \phi \frac{\partial}{\partial \Bar{z_i} },\frac{\partial}{\partial \Bar{z_j} }+\phi \frac{\partial}{\partial \Bar{z_j} }] \in T^{0,1}_tM.
$$
Expressing $\phi$ as $\sum\limits_{k,l} \phi_{kl}d\Bar{z_k}\otimes \frac{\partial}{\partial z_l }$, we have 
\begin{align*}
    [\frac{\partial}{\partial \Bar{z_i} }+ \phi \frac{\partial}{\partial \Bar{z_i} },\frac{\partial}{\partial \Bar{z_j} }+\phi \frac{\partial}{\partial \Bar{z_j} }]  &=[\frac{\partial}{\partial \Bar{z_i} }+ \sum_l \phi_{il} \frac{\partial}{\partial {z_l} },\frac{\partial}{\partial \Bar{z_j} }+\sum_l \phi_{jl} \frac{\partial}{\partial {z_l} }] \\
    & =  [\frac{\partial}{\partial \Bar{z_i} },\sum_l \phi_{jl} \frac{\partial}{\partial {z_l} }] +[\sum_l \phi_{il} \frac{\partial}{\partial {z_l} }, \frac{\partial}{\partial \Bar{z_j} }]+[\sum_l \phi_{il} \frac{\partial}{\partial {z_l} },\sum_l \phi_{jl} \frac{\partial}{\partial {z_l} }] \\
    &= \sum_l \frac{\partial \phi_{jl}}{\partial \Bar{z_i}} \frac{\partial}{\partial z_l}-\sum_l \frac{\partial \phi_{il}}{\partial \Bar{z_j}} \frac{\partial}{\partial z_l}+[\sum_l \phi_{il} \frac{\partial}{\partial {z_l} },\sum_l \phi_{jl} \frac{\partial}{\partial {z_l} }] \\
    &=\bar{\partial}\phi (\frac{\partial}{\partial \Bar{z_i} },\frac{\partial}{\partial \Bar{z_j} }) + \frac{1}{2}[\phi,\phi](\frac{\partial}{\partial \Bar{z_i} },\frac{\partial}{\partial \Bar{z_j} })
\end{align*}

Hence we have $$(\Bar{\partial}\phi +\frac{1}{2}[\phi,\phi])(\frac{\partial}{\partial \Bar{z_i} },\frac{\partial}{\partial \Bar{z_j} })\in T^{0,1}_tM \cap T^{1,0}M.$$

But when $t$ is small $T^{0,1}_tM \cap T^{1,0}M= \emptyset$. So  $\Bar{\partial}\phi +\frac{1}{2}[\phi,\phi]=0$.  
\smallskip
Conversely, if $\Bar{\partial}\phi +\frac{1}{2}[\phi,\phi]=0$, we have $$[\frac{\partial}{\partial \Bar{z_i} }+ \phi \frac{\partial}{\partial \Bar{z_i} },\frac{\partial}{\partial \Bar{z_j} }+\phi \frac{\partial}{\partial \Bar{z_j} }]=(\Bar{\partial}\phi +\frac{1}{2}[\phi,\phi])(\frac{\partial}{\partial \Bar{z_i} },\frac{\partial}{\partial \Bar{z_j} })\in T^{0,1}_tM $$ for any basis vectors $\frac{\partial}{\partial \Bar{z_i} }, \frac{\partial}{\partial \Bar{z_j} }$ in $T^{0,1}M$. Since  $T^{0,1}_tM=(\text{id}+\phi)T^{0,1}M$, we conclude that $[T^{0,1}_tM,T^{0,1}_tM]\subset T^{0,1}_tM$.

\end{proof}

 \medskip

\subsection{First Order Deformation and its Obstruction}    
Previously  we have seen that a   deformation of a complex structure $(M,J)$ can be seen as a family of solutions to the Maurer-Cartan equation. To solve MC equation formally, i.e., to find formal power series which satisfies this differential equation, it is natural to look for the first order solution firstly. Consider the formal power series expansion of $\phi_t$ (we write $\phi(t)$ instead of $\phi_t$),
$$
\phi (t)= \phi_0+\phi_{1}t+\phi_{2} t^2+\ldots .
$$
Since $\phi(0)=0$,  $\phi_0=0$. Substitute the series into  MC equation then we get 
 
\begin{equation} \label{eq:9}
\begin{split}  
& 0= \bar{\partial} \phi_{1} \\
& 0=\bar{\partial} \phi_{2}+\frac{1}{2}[\phi_{1}, \phi_{1}] \\
& \vdots\\
& 0= \bar{\partial} \phi_{k}+\frac{1}{2}\sum_{i+j=k}[\phi_{i}, \phi_{j}]\\
& \vdots \\
\end{split}
\end{equation}
Consider the first order term of the equations in Eq. \ref{eq:9}.
\begin{definition}
     If there exists an element $\phi_1 \in \Omega ^{0,1}({T^{1,0}M})$  satisfying the first equation in Eq.\eqref{eq:9}$,$ i.e. $\Bar{\partial}(\phi_1)=0$, we call $\phi_1$ a \textbf{first order deformation} of the complex structure. 
\end{definition}

Similarly, we call $\sum\limits^{k}\limits_{i=1}t^k \phi_k$  a $k$  \textbf{order deformation} if it satisfies Maurer-Cartan equation up order $k$, that is, it satisfies the first $k$ equations in Eq.\ref{eq:9}.

\begin{remark}
  For any deformation of complex structure $\phi_t$, the first order term of its power series expansion is a first order deformation. The converse may not be true. Given a first order deformation, i.e., given $\phi_1$ such that $\bp \phi_1=0$, there may not exist a solution of MC equation $\phi_t$ such that $\phi_1$ is the first order term of $\phi_t$. For example, from Equation \ref{eq:9} we see that a first order deformation $\phi_1$ cannot be the first order term of a solution of the MC equation if $[\phi_1,\phi_1] $ is not exact.  
\end{remark}

The set of first order deformation is a subset of the set of closed elements in $\Omega^{0,1}(T^{1,0}M)$, actually, after we identify isomorphic deformations, we will see that the set of first order deformation is in bijection with $H^{0,1}(T^{1,0}M)$.

\medskip

\subsection{Isomorphism between Deformations}
Two complex structures $J,J'$ on $M$ are isomorphic if there is a diffeomorphism $f: M\longrightarrow M$ such that $J'$ is the push-forword of $J$ by $f$. Similarly we define a notion of isomorphism between deformations.

\begin{definition}
    We say two  deformations $J_t$, $J'_t$ of complex structure on $M$ are isomorphic if there exists a one-parameter group of  diffeomorphisms $I_t$ which is  also parametrized by $t\in(-\epsilon,\epsilon)$  $$I_t: 
    (-\epsilon,\epsilon) \times M\to M,$$ such that  $J'_t\circ dI_t=dI_t\circ J_t$ for any $t$. In other words, $J'_t$ is the push forward of $J_t$ by $ I_t$.
\end{definition}

\begin{lemma} \label{l:6.9}
    Let  $J'_t,J_t$ be deformations of complex structures $(M,J)$. The corresponding deformation encoding maps are denoted by $\phi'(t), \phi(t)$, inducing first order deformations $\phi'_1, \phi_1$, respectively. Then $ \phi'_1- \phi_1  \text{ is the trivial element in }  H^{0,1}(T^{1,0}M)$ if and only if $J'_t$ and $J_t$ are isomorphic deformations.
\end{lemma}
\begin{proof}
We firstly prove the case when $J_t$ does not change with $t$, that is, $J_t\equiv J $. By definition of encoding map, $\phi(t)=\pi ^{1,0} \circ ({\pi^{0,1}|_{T^{0,1}_t}})^{-1}\equiv 0$. Therefore $ \phi'_1-\phi_1 $  is exact in $ \Omega^{0,1}(T^{1,0}M)$  is equivalent to say that $ \phi'_1 $  is exact in $ \Omega^{0,1}(T^{1,0}M)$. It is enough to check this locally on an small open set $U$. We choose holomorphic coordinates $(z_1,...,z_n)$ on $U$ such that $\frac{\partial}{\partial z_i} \in T^{1,0}U$. For simplicity we omit $U$ in e.g. $T^{1,0}U, T_t^{0,1}U$.

If $J'_t$ is isomorphic to $J_t$, there exists a family of diffeomorphisms $I_t$ such that $J_t$ is the push-forward of $J_0(t)$ by $I_t$. So for any $p\in U$ and $I_t(p)\in U$,  the holomorphic part of $J_t$ at $I_t(p)$, denoted by $T_{t,I_t(p)}^{1,0}$, is just $dI_t ( T_p^{1,0})$, similar for the anti-holomorphic part. Now we determine the first order term of the encoding map of $J'_t$, which is the first order deformation we want to check.

We write $I_t$ in this coordinate chart as $$I_t(z_1, \ldots, z_n) = (y_1(t, z_1, \ldots, z_n,\bar{z_1},\ldots, \Bar{z_n}), \ldots, y_n(t, z_1, \ldots, z_n,\Bar{z_1},\ldots, \Bar{z_n})),$$ and $y_k$ in a power series:

$$
y_k = z_k + t f_k(z_1, \ldots, z_n,\bar{z_1},\ldots, \Bar{z_n}) + o(t^2).
$$

    For any $\frac{\partial}{\partial \Bar{z_i}} \in T^{0,1}_p$, $$d{I_t}  (\frac{\partial}{\partial  \Bar{z_i}})=\sum\limits_j\frac{\partial y_j}{\partial  \Bar{z_i}}\frac{\partial }{\partial  {z_j}}+\sum\limits_j\frac{\partial \bar{y_j}}{\partial  \Bar{z_i}}\frac{\partial }{\partial  \Bar{z_j}} = \frac{\partial }{\partial  \Bar{z_i}} + \sum\limits_j t\frac{\partial f_j}{\partial  \Bar{z_i}}\frac{\partial }{\partial  {z_j}} +\sum\limits_j t\frac{\partial \bar{f_j}}{\partial  \Bar{z_i}}\frac{\partial }{\partial  \Bar{z_j}} +o(t^2)\in T^{0,1}_{t,I_tp}$$  and $\{dI_t(\frac{\partial}{\partial \bar{z_i}}) \}$ is the basis of $T^{0,1}_{t,I_t(p)}$. To find $\phi(t)=\pi^{1,0}\circ ({\pi^{0,1}|_{T^{0,1}_t}})^{-1}$, we first determine $({\pi^{0,1}|_{T^{0,1}_t}})^{-1}  (\frac{\partial}{\partial \bar{z_k}})$, that is, we want to find $\lambda^k_l \in \C$, $l=1,...,n$ such that $$\pi^{0,1}|_{T^{0,1}_t}\circ \left(\sum\limits_l \lambda^k_l dI_t\left(\frac{\partial}{\partial \bar{z_l}}\right)\right)=\frac{\partial}{\partial \bar{z_k}}.$$
    
    In matrix form, we define $n \times n$ matrix $F$ as $[F]^p_{q}=\frac{\partial \bar{f_p}}{\partial \bar{z_q}}$ and $\lambda^k=(\lambda^k_1,...,\lambda^k_n)^t$.  After writing  $\frac{\partial}{\partial \bar{z_k}}$ in matrix form as $e_k=(0,0,...0,1,0,...0)^t$, what we want to solve now is 
    \begin{equation}\label{eqn:star}
   (\text{I}+tF)\lambda^k = e_k  \tag{*}
  \end{equation}
    
    Since we want to determine first order deformation, we omit higher order terms. Then we have $(1+tF)^{-1}=1-tF$ and  $\lambda^k=e_k-tF e_k$. Therefore $$({\pi^{0,1}|_{T^{0,1}_t}})^{-1} \left(\frac{\partial}{\partial \bar{z_k}}\right)=\sum \limits_l \lambda_l dI_t\left(\frac{\partial}{\partial \bar{z_l}}\right).$$  After omitting higher order terms we have $$\phi'(t)=\pi^{1,0}\circ ({\pi^{0,1}|_{T^{0,1}_t}})^{-1}=\sum\limits_k dz_k \otimes ({\pi^{1,0} \circ \pi^{0,1}|_{T^{0,1}_t}})^{-1}\left(\frac{\partial}{\partial \bar{z_k}}\right)  +o(t^2)=t\sum\limits_{j,k} dz_k \otimes \frac{\partial f_j}{\partial \bar{z_k}} \frac{\partial}{\partial {z_j}}+o(t^2).$$ Then
    \begin{equation} \label{eq:6.9}
        {\phi'_1=\sum_{j,k} dz_k \otimes \frac{\partial f_j}{\partial \bar{z_k}} \frac{\partial}{\partial {z_j}}= \sum_j \Bar{\partial}\left(f_j \frac{\partial}{\partial z_j}\right)}
    \end{equation}
    hence it is exact.

Conversely, given $\phi'_1$. The first order term of $\phi'(t)$, if $\phi'_1$ is exact in $ \Omega^{0,1}(T^{1,0}M)$, we can write $\phi_1$ as  $\phi_1=\sum\limits_j \bar{\partial}(  f_j \frac{\partial}{\partial z_j})$. For vector  field $\sum\limits_j f_j \frac{\partial}{\partial z_j}$, we define $I_t$ as the flow of this vector field, which is a one parameter group of diffeomorphisms on $M$. To prove $J'_t$ is isomorphic to $J$, it is  sufficient to prove 
\begin{equation} \label{}
dI_t\left(\frac{\partial}{\partial \bar{z_i}}\right)=\frac{\partial}{\partial \Bar{z_i}}+\phi'(t)\left(\frac{\partial}{\partial \Bar{z_i}}\right).
\end{equation}
Since $\frac{\partial}{\partial \Bar{z_i}}+\phi'(t)\left(\frac{\partial}{\partial \Bar{z_i}}\right)$ spans the anti-holomorphic tangent space of $J'_t$, hence if the equation is satisfied, the change of anti-holomorphic part of tangent space is induced by the push forward by $I_t$. 

When $t=0$, the equation is clearly satisfied. To prove the equation for arbitrary $t$, it is  sufficient to show that the time derivatives on both sides of the equation are equal. Hence we make derivative with respect to $t$ at $t=0$ for both side, then   $$\text{LHS}=-[\sum_j f_j \frac{\partial}{\partial z_j}, \frac{\partial}{\partial \Bar{z_i}}],$$ and $$\text{RHS}=\phi'_1(\frac{\partial}{\partial \Bar{z_i}})=  \Bar{\partial}(\sum \limits_j f_j\frac{\partial}{\partial {z_j}})(\frac{\partial}{\partial \Bar{z_i}})=\sum \limits_j \frac{\partial f_j}{\partial \Bar{z_i}} \frac{\partial}{\partial \Bar{z_i}}.$$ We see that LHS=RHS. \\

To prove the general case (in which $J_t$ is arbitrary) notice that we can still use the steps in which we proved when $J_t=J$ for any $t$ but replace $\frac{\partial}{\partial \bar{z_j}}$ by $\frac{\partial}{\partial \Bar{z}_{j,t}}=\frac{\partial}{\partial \bar{z_j}}+\phi(t)\frac{\partial}{\partial \bar{z_j}}$ which is the anti-holomorphic part of $J_t$,  also $\phi'(t)$ by $\tilde{\phi'}(t)$, where $\tilde{\phi'}(t)$ is the encoding map of $J'_t$ if we replace the underlying complex structure $J$ by $J_t$. Hence we have 
$$\frac{\partial}{\partial \bar{z_j}}+\phi(t)\frac{\partial}{\partial \bar{z_j}} + \tilde{\phi'}(t)\left(\frac{\partial}{\partial \bar{z_j}}+\phi(t)\frac{\partial}{\partial \bar{z_j}}\right)= \frac{\partial}{\partial \bar{z_j}}+\phi'(t)\frac{\partial}{\partial \bar{z_j}}.$$ When considering only orders lower than 2, we have $\tilde{\phi'}(t)= \phi'(t)-\phi(t)$ and we can replace   $\frac{\partial}{\partial z_{j,t}}$ with $\frac{\partial}{\partial \bar{z_j}}$. Therefore, we can just replace $\phi'(t)$ by $\phi'(t)-\phi(t)$ hence we have the same result. That is,  $J'_t$ is isomorphic to $J_t$ if and only if $\phi'_1-\phi_1$ is exact
\end{proof}

\begin{corollary}\label{c:6.11}
There is a canonical bijective map between the set of equivalence classes of first order deformations on $M $ and the set $H^{0,1}(T^{1,0}M)$ .
\end{corollary}
\begin{proof}
By definition, first order deformations are closed elements in $\Omega^{0,1}(T^{1,0}M)$. By Lemma \ref{l:6.9} we see that the equivalence classes of first order deformations are the cohomological equivalence classes. Hence the canonical projection map from a first order deformation to its equivalence class of first order deformations send it to an element in $H^{0,1}(T^{1,0}M)$.
\end{proof} 
We  call an equivalence class of first order deformation, or equivalence class in $H^{0,1}(T^{1,0}M)$ Kodaira-Spencer class.

 Hence we say that $H^{0,1}(T^{1,0}M)$ is the space of first order deformations. If every first order deformation of a complex manifold $(M,J)$ is integrable, we say that this complex manifold is \textbf{unobstructed}.

\medskip

\subsection{Formal Deformation, Obstruction of First Order Deformation}
From last section we see that $H^{0,1}(T^{1,0}M)$ is the space of first order deformations, but not every first order deformation is the first order term of a solution of the MC equation. To solve the  MC equation,  we now want to determine which first order deformations can be extended to a solution of  the MC equation. To avoid facing convergence problem, we will consider formal solutions of the MC equation and define formal deformation of complex structure as  formal solution of the MC equation.

\begin{definition}\label{d:3.13}
   We call the formal power series  $\phi [[t]]= \phi_0+\phi_{1}t+\phi_{2} t^2+\ldots \in \Omega^{0,1}(T^{1,0}M)[[t]]$
    with coefficients  $\phi_i\in \Omega^{0,1}({T^{1,0}M})$ and formal variable $t$  a \textbf{formal deformation} of a complex structure  if $\phi [[t]]$ satisfies the Maurer-Cartan equation $$\Bar{\partial}\phi [[t]]+ \frac{1}{2}[\phi [[t]],\phi [[t]]]=0,$$ or equivalently, satisfies the set of Equations \ref{eq:9}.  
\end{definition} 

Additionally, if a first order deformation $\phi$ is the first order term of a formal deformation, we say $\phi$ is integrable. If every element in $H^{0,1}(T^{1,0}M)$ is integrable, we say that $M$ is \textbf{unobstructed}.

\begin{theorem}{\label{T:3.13}} 
    If $H^{0,2}(T^{1,0}M)$ is trivial, then $M$ is unobstructed, that is, every Kodaira-Spencer class $[\phi]$ in $H^{0,1}(T^{1,0}M)$ can be integrated. Hence we say that obstruction of the integrability of first order deformation lies in $H^{0,2}(T^{1,0}M)$.
\end{theorem}
\begin{proof}
    For arbitrary closed element $\phi_1\in \Omega^{0,1}(T^{1,0}M) $, we want to extend  it to a formal power series $\phi[[t]]$ that satisfies Maurer-Cartan equation. To the end of this, we extend solutions by induction. Suppose we have a $k$ order solution $$S_k := \sum \limits^{k}_{i=1} t^i \phi_i,$$ where $\phi_i \in \Omega^{0,1}(T^{1,0}M) $. We want to find $\phi_{k+1}$ such that $S_{k+1}=S_k+t^{k+1}\phi_{k+1}$ satisfies MC equation up to order $k+1$. Hence, we need to solve 
    $$
    \Bar{\partial}(S_k+ t^{k+1} \phi_{k+1}) +\frac{1}{2}[S_k+ t^{k+1} \phi_{k+1},S_k+ t^{k+1} \phi_{k+1}]=0.
    $$
The vanishing of coefficients with $t^{k+1}$ leads to the equation
$$
\frac{1}{2}\sum \limits_{i+j=k+1 \atop  0<i<j<k+1}[\phi_i,\phi_j]=\Bar{\partial}(\phi_{k+1}).
$$
To determine if this equation is solvable, notice that the left hand side of this equation a closed element in $\Omega^{0,1}{(T^{1,0}M)}$ by Lemma \ref{l:3.7}. If $H^{0,2}(T^{1,0}M)=0$ (which means that every closed element in $\Omega^{0,1}{(T^{1,0}M)}$ is exact) then there exists a $\phi_{k+1}$ such that the equation above is satisfied.
\end{proof}

\medskip

Hence any complex manifold $M$ with trivial $H^{0,2}(T^{1,0}M)$ is unobstructed. However $H^{0,2}(T^{1,0}M)$ is only a sufficient but not necessary condition for a manifold to be unobstructed. As we will see Calabi-Yau manifolds with non-vanishing $H^{0,2}(T^{1,0}M)$ are unobstructed.

\begin{example}{(Complex Structures on 2-Torus)}

    The complex torus $\mathbf{T}_0=\mathbb{C}/(\mathbb{Z}+\mathbb{Z}i)$ with the canonical complex structure inherited from complex plane is a complex manifold. If $\mathbf{T}_0$ has another complex structure, it would induce another complex structure on its universal covering space $\R^2$. By the uniformization theorem, $\R^2$ with a complex structure is biholomorphic to complex plane. Hence $\mathbf{T}_0$ with any other complex structure is again a quotient of complex plane by a rank 2 lattice which is not equivalent to the original lattice $(1,i)$.  
    
    Therefore to define a deformation of complex structures on $\mathbf{T}_0$, we first define a family of tori over $t\in (-\epsilon, \epsilon)\subset \R$,\, ${\mathbf{T}_t}=\mathbb{C}/(\mathbb{Z}+\mathbb{Z}\tau(t))$ by the quotient of lattice, where $\tau: (-\epsilon,\epsilon)\longrightarrow \mathbb{H}$ with $\tau(0)=i$, here $\mathbb{H}$ is the upper half complex plane.

    For any $t$ there is a homeomorphism $f(t)$ that maps $\mathbf{T}_t$ onto $\mathbf{T}_0$. The push-forward of the canonical complex structure on $\mathbf{T}_t$ to $\mathbf{T}_0$ gives a family of complex structures $J_t$ on $\mathbf{T}_0$, which is a deformation of complex structure. Indeed, the  almost complex structures induced in this way are integrable since $J_t$ and Lie bracket commute with the push-forward by $f(t)$, which implies that the Nijenhuis tensor vanishes. In particular, $(\mathbf{T}_0,J_t)$ is biholomorphic to $\mathbf{T}_t$ as complex manifolds, see figure \ref{fig:tori}.

    We now do this explicitly. Firstly consider the tori as (integrable) almost  complex manifolds with real coordinates inherited from complex plane, i.e. we identify the point with real coordinates  $(x,y)$   with the point $(x+yi)$ in complex coordinates. In the canonical real coordinates  chart, the canonical complex structure of $\mathbf{T}_0$ is $J_0= \begin{bmatrix}
        0&-1 \\
        1&0
    \end{bmatrix}$. We write  the complex structure $J_t$ on ${\mathbf{T}_t}$   as $\bar{J_t}=\begin{bmatrix}
        0&-1 \\
        1&0
    \end{bmatrix}$ since the complex structure of $\mathbf{T}_t$ is also inherited from complex plane. For simplicity we denote $$R_t:= Re(\tau(t)), \,\, I_t=Im(\tau(t)).$$ Define 
    \begin{equation*}
        \begin{split}
            f(t):{\mathbf{T}_t} &\longrightarrow {\mathbf{T}_0} \\
            [x,y]&\longmapsto [x-y\frac{R_t}{I_t},\frac{y}{I_t}].
        \end{split}
    \end{equation*}
      It is easy to verify that $f(t)$ is a well-defined  homeomorphism from $\mathbf{T}_t$ to $\mathbf{T}_0$. Now we define a  family of complex structures over $\mathbf{T}_0$ by setting $J_t=Df \cdot \bar{J_t} \cdot (Df)^{-1}$, in this way we have defined a deformation of complex structure.

\smallskip
    
Expressed in canonical coordinates we have $Df= \begin{bmatrix}
    1&-\frac{R_t}{I_t}\\
    0&\frac{1}{I_t}
\end{bmatrix}$, substitute in $J_t$ we have $$J_t=\begin{bmatrix}
    -R_t&-I_t-\frac{R^2_t}{I_t} \\
    1&\frac{R}{I_t}. 
\end{bmatrix}.$$

\smallskip

To identify isomorphic deformations, we see that a necessary condition for two deformations $J_t,J'_{t}$ induced by $\tau, \tau'$ are isomorphic is that for any $t$ there exists an diffeomorphism $I:\mathbf{T}_0\longrightarrow \mathbf{T}_0$ such that $J_t\circ dI=dI\circ J'_{t}$, which implies that as complex manifolds, $(\mathbf{T}_0,J_t)$ is isomorphic (or biholomorphic) to $(\mathbf{T}_0,J'_t)$. Since $(\mathbf{T}_0,J_t), (\mathbf{T}_0,J'_t)$ are biholomorphic to $\mathbf{T}_t,\mathbf{T}'_t$ respectively, $J_t$ is isomorphic to $J'_{t}$ only if $\mathbf{T}_t \cong \mathbf{T}'_t$ as complex manifolds. The uniformization theorem for elliptic curves says that two complex tori are isomorphic if and only they are generated by the equivalent lattice. In our case, when $\epsilon$ is small, $(1,\tau(t))$ and $(1,\tau'(t))$ represent the same lattice if and only if $\tau(t)=\tau'(t)$. Hence two deformations are isomorphic if only if $\tau(t)=\tau'(t)$ for any $t$, that is, we don't need to worry about isomorphic deformations.

Therefore a  deformation of complex structure on $\mathbf{T}_0$ is represented by a short path in the neighborhood of $i$ starting from $i$.  Two deformations are isomorphic if the corresponding two short paths coincide around  $i$. Hence locally a small neighborhood around $i$ represents the set of complex structures which are close to the standard one. The moduli space of complex structures on $\mathbf{T}_0$ locally at the standard complex structure has complex dimension 1.

Since the complex dimension of $\mathbf{T}_0$ is one, $H^{0,2}(T^{1,0} \mathbf{T}_0)$ is trivial. By Theorem \ref{T:3.13},  $\mathbf{T}_0$ is unobstructed and the space first order deformation $H^{0,1}(T^{1,0}\mathbf{T}_0)$ has complex dimension one. Every first order deformation is now integrable and first order deformation can be seen as the tangent space of solution space of Maurer-Cartan equation, which has dimension 1.

 \begin{figure}[h]
    \centering
    \includegraphics[width=8cm]{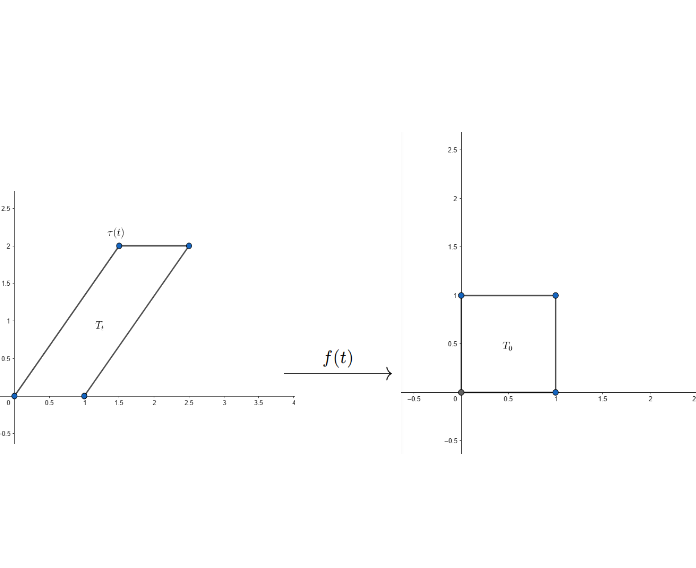}
    \caption{A family of complex structures on torus}
    \label{fig:tori}
\end{figure}

\end{example}

\begin{example}{(Complex Structures on $S^2$)}
Let $\mathbb{P}^1$ be one dimensional complex projective space. It is known that $H^{0,1}(T\mathbb{P}^1)$ is trivial,  so we cannot deform the complex structure of $\mathbb{P}^1 $ (we call this complex structure rigid).

    Indeed, by the uniformization theorem, a one-dimensional complex manifold is holomorphically isomorphic to one of the following, an open unit disk, complex plane, or Riemann sphere. In this case, this implies that $S^2$ has only one complex structure inherited from $\mathbb{P}^1$. 
\end{example}

\newpage
\section{Deformation theory via Differential Graded Lie algebra}

In the last chapter, we have considered the problem of deformations of complex structures. A formal deformation of a complex structure is a formal solution $(\phi_t)$ to a MC equation.

We have tackled this problem by considering first order solutions  to the MC equation $\phi_1$ (Sec. 3.2). To get a formal solution one needs to extend first order solutions to an arbitrary order solution (Sec. 3.4). The MC equation is encoded in 
$$\Omega^{0,1}(T^{1,0}M)\xrightarrow{\bar{\partial}} \Omega^{0,2}(T^{1,0}M).$$

\smallskip

In this chapter, we study this problem using another formalism.  
Notice that the complex $(\Omega^{0,*}(T^{1,0}M),\bar{\partial})$ has not only a differential complex structure, but also a Lie bracket. A differential complex with a compatible Lie bracket is called a {\it differential graded Lie algebra}. The idea 

\smallskip

``{\it If $X$ is a moduli space over a field $k$ of characteristic zero, then a formal neighborhood of any point
$x\in X$ is controlled by a differential graded Lie algebra}'' (see for instance \cite{Lurie1}). 

\smallskip 
 was developed in unpublished work of Deligne, Drinfeld, and Feigin.

In this chapter, we introduce  deformation theory using the differential graded Lie algebra formalism and reformulate  deformations of complex structures in this setting. Finally, we follow the strategy in \cite{BK} and give the proof of unobstructness of a Calabi-Yau manifold. The main references used are \cite{BK} and  \cite{Manetti}.

\medskip 

\subsection{Differential Graded Lie Algebra and Maurer-Cartan Equation}

\begin{definition} 
A differential graded Lie algebra $L$ (DGLA in short) consists of:  
\begin{enumerate}
    \item A differential graded vector space $\bigoplus\limits_{i\in \mathbb{Z}}L^i$, i.e. a complex with differential $d:L^n \xrightarrow{}L^{n+1}$. An element in $L^i$ is said to be a homogeneous element of degree $i$.
    \item A Lie algebra structure which is compatible with (1), i.e. for any $i,j$ there exists a bilinear map such that 
    $$[\cdot, \cdot]: L^{i} \otimes L^{j}\rightarrow L^{i+j}$$
    and for any $a$,$b \in L$,
  $$
  \begin{aligned}
      &[a,b]=-(-1)^{\Tilde{a}\Tilde{b}}[b,a]  &  \text{(graded-anti-symmetry)} \\
      &[a,[b,c]]=[[a,b],c]+(-1)^{\Tilde{a}\Tilde{b}}[b,[a,c]] & \text{(graded Jacobi identity)}\\
      &d([a,b]) = [da,b] + (-1)^{\widetilde{a}}[a,db] & \text{(compatible with differential)},
  \end{aligned}
  $$  
where $\Tilde{x}$ is the degree of $x$ for any $x\in L$.     
    
\end{enumerate}
 A morphism between two DGLAs, say $L$ and $L'$ is a morphism of complexes which preserve the Lie algebra structure. Explicitly,  $f:={(f^i)}_{i\in \mathbb{Z}}:L \longrightarrow L'$ is a morphism of DGLA if the diagram 
$$\begin{tikzcd}
	{} & {...} & {L^{i-1}} & {L^i} & {L^{i+1}} & {...} \\
	{} & {...} & {L'^{i-1}} & {L'^{i}} & {L'^{i+1}} & {...}
	\arrow["{d_{i-1}}", from=1-3, to=1-4]
	\arrow[from=1-2, to=1-3]
	\arrow["{f^{i-1}}", from=1-3, to=2-3]
	\arrow["{f^{i}}", from=1-4, to=2-4]
	\arrow["{d_i}", from=1-4, to=1-5]
	\arrow["{d'_{i-1}}", from=2-3, to=2-4]
	\arrow["{d_i}", from=2-4, to=2-5]
	\arrow[from=2-5, to=2-6]
	\arrow[from=1-5, to=1-6]
	\arrow["{f^{i+1}}", from=1-5, to=2-5]
	\arrow[from=2-2, to=2-3]
\end{tikzcd}$$
commutes and $f^i [a,b]_L=[f^i (a),f^i (b)]_{L'}$ for any $i$.
\end{definition}
\begin{remark}
    
[The algebra of deformation of complex structure on a complex manifold $M$]
    We generalize the $\bar{\partial}$ and the bracket on $\Omega^{0,1}(T^{1,0}M)$ to $\Omega^{0,q}(T^{1,0}M)$ for any $q\geq 1$ in local coordinates as 
    $$
    \begin{aligned}
        &  \bar{\partial}: \Omega^{0,q}(T^{1,0}M) \longrightarrow \Omega^{0,q+1}(T^{1,0}M)\\
        &       \alpha= \sum \limits_{Ij}\alpha^{Ij}d\bar{z}^I \otimes \frac{\partial} {\partial z_j}  \longmapsto  \sum \limits_{Ijk}\frac{\partial \alpha^{Ij}}{\partial \bar{z}^k} d \bar{z}^k \wedge d \bar{z}^I \otimes \frac{\partial}{\partial z_j}  ,
        \\
       &   [\cdot,\cdot] : \Omega^{0,p}(T^{1,0}M) \times \Omega^{0,q}(T^{1,0}M) \longrightarrow \Omega^{0,p+q}(T^{1,0}M)
       \\
       &   (\sum \limits_{Ir}\alpha^{Ir}\bar{z}^I \otimes \frac{\partial} {\partial z_r},\sum \limits_{Js}\beta^{Js}d\bar{z}^J \otimes \frac{\partial} {\partial z_s} ) \longmapsto \sum \limits_{IrJs}d\bar{z}^I \wedge d\bar{z}^J \otimes [\alpha^{Ir}\frac{\partial} {\partial z_r},\beta^{Js}\frac{\partial} {\partial z_s} ],
    \end{aligned}
    $$
Here we used  uppercase letter to denote multi-index e.g. $dz^I=dz_{i_1}\wedge ...\wedge dz_{i_q}$. 
One can verify that this defines a DGLA with grading $\bigoplus\limits_i \Omega^{0,i}(T^{1,0}M)$. We call this algebra \textbf{Kodaira-Spencer algebra} of $M$, denoted by $KS_M$.

\end{remark}

\begin{lemma} \label{e:7.3}

Given a differential graded Lie algebra $L$ and a graded commutative algebra $\mathfrak{m}$, we can define a new DGLA $L \otimes \mathfrak{m}$ by setting

$$(L \otimes \mathfrak{m})^k = \bigoplus_{k=i-j} (L^i \otimes \mathfrak{m}^{j})$$
for each $k$.  The differential  $d$ on  $L \otimes \mathfrak{m}$ is 

$$d(x \otimes a) = dx \otimes a $$

where $x$ belongs to $L$ and $a$ belongs to $\mathfrak{m}$. The  Lie bracket is defined as

$$[x \otimes a, y \otimes b] = (-1)^{deg(a)deg(y)}[x,y] \otimes ab.$$

In particular, when $\mathfrak{m}$ is concentrated in degree $0$, that is, $\mathfrak{m}^i=0$ if $i\neq 0$,  
$$(L \otimes \mathfrak{m})^k = L^k\otimes \mathfrak{m} $$ and 

$$
[x \otimes a, y \otimes b] = [x,y] \otimes ab
$$
\end{lemma}

\begin{definition}
    Let $L$ be a DGLA. The \textbf{Maurar-Cartan equation} of L is $$da +  \frac{1}{2} [a,a]=0,$$ where $a \in L^1$. The solution space of the Maurer--Cartan equation is called Maurer--Cartan space and denoted by $MC(L) \subset L^1$.
\end{definition}

\medskip

\subsection{Gauge Action on Differential Graded Lie Algebra}
As the case of deformation of complex structure, in which we have a notion of isomorphism 
between deformations,  we want to  restate this isomorphism in DGLA formalism. In this section we  define gauge action on $MC(L)$, which is the Maurer-Cartan space of a DGLA $L$.

To introduce the gauge action on $MC(L)$, we need a notion of nilpotent (graded) Lie algebra. 
Recall that an algebra $A$ is called nilpotent if there exist $k\in \mathbb{N}$, such that $$A^k=\{a_1 a_2 ... a_k,\, a_i\in A \text{ for } i=1,...,k\}=\{0\}.$$

\begin{definition}
Let $L$ be a Lie algebra. Define $L^{(1)}:=L$, $L^{(2)}=[L,L]=\{ [l_1,l_2]|\, l_1,l_2 \in L   \}$, and inductively for any higher $n$, $L^{(n)}=[L,L^{(n-1)}]$. We call $L$ \textbf{nilpotent} if there exists $n>0$ such that $L^{(n)}=\{0\}$.
\end{definition}
It is easy to see that:
\begin{lemma}
    Let $L$ be a Lie algebra and $A$ be a nilpotent algebra, then the algebra $L \otimes A$  with Lie bracket $[x \otimes a, y \otimes b] = [x,y] \otimes ab$ is nilpotent.
\end{lemma}

For a nilpotent Lie algebra $L$ we can define the exponential of an element $a\in L$ as $$exp(a)=\text{id}+a+\frac{a^2}{2!}+... \quad.$$

Since $L$ is nilpotent, this  sum  has only finite terms. By the Baker–Campbell–Hausdorff formula (BCH) we can define a group structure on $exp(L)=\{ exp(a)|\, a \in L   \}$. 

\begin{proposition}
    
    Let $L$ be an associative nilpotent Lie algebra. Then $exp(L)$ forms a group with elements $\{ exp(l), \,  l\in L \}$ and group multiplication $exp(a)exp(b)=exp(a \bullet b)$, where $$a \bullet b=a+b+\frac{1}{2}[a, b]+\frac{1}{12}[a,[a, b]]+\frac{1}{12}[b,[b, a]]+\frac{1}{24}[a,[b,[b, a]]]+\ldots$$ is given by the Taylor expansion of $log(exp(a)exp(b))$, i.e. Baker–Campbell–Hausdorff (BCH) formula. We call it the BCH product. 
\end{proposition}
\begin{proof}
    The unit of $exp(L)$ is $exp(0)$. 
    The inverse of $exp(l)$ is $exp(-l)$.   Associativity $(exp(a)exp(b))\\exp(c)=exp(a)(exp(b)exp(c))$ can also be checked by the expansion of Taylor series.
\end{proof}

We call an ideal $I$ of a ring (or an algebra over $\K$) $R$ \textbf{nilpotent} if $I^k:=\{a_1a_2,...,a_k , a_i \in I\}={0}$ for some $k$.  If $L$ is a unitary associative Lie algebra (not necessarily nilpotent) and $I$ is a nilpotent ideal of $L$, then the BCH formula is well-defined on $I$. The theorem above also works for $exp(I)$.

Given a nilpotent DGLA $L$,  its $0$ degree term $L^0$ is also nilpotent. For any element $a\in L^0$, consider the adjoint operator of $a$, $ad_a:L\xrightarrow[]{}L$, $ad_a b=[a,b]$. Since $L^0$ is nilpotent, the exponential of $ad_a$ ,
$$
exp{(ad_a)}b=\sum \limits_{i\geq 0} \frac{ad_a^i}{i!}b=\sum \limits_{i\geq 0} \frac{\overbrace{[a,...[a,[a}^{i \, \, \text{times} \,\, a},b]]]}{i!}
$$
is well-defined and we have

\begin{lemma}({\cite{Manetti}}\label{L.4.8}, Ch. 6)
The exponential $exp(ad_a)$ is an automorphism of $L$ and the set $\{ [x,x]=0, x\in L \} \subset L$ is stable under the action $exp(ad_a)$, i.e., $exp(ad_a)$ preserves $\{ [x,x]=0, x\in L \}$.
\end{lemma}
\begin{proof}

By Jacobi identity $[a,[b,c]]=-[b,[c,a]]-[c,[a,b]]$, we have $ad_a([b,c])=[ad_a b, c]+[b,ad_a c]$.  Hence $ad_a ^i [b,c]= \sum \limits_{m\le i} {i \choose m} [ad_a^m b, ad_a^{i-m}c] $. Then $$exp(ad_a)[b,c]=\sum \limits_{i\geq 0}\frac{1}{i!} ad_a^i[b,c]=\sum \limits_{i\geq 0}\frac{1}{i!}\sum \limits_{m \leq i } {i \choose m}[ad_a^m b,\\ ad_a^{i-m} c]=[exp(ad_a)b,exp(ad_a)c].$$ Also $exp(ad_a)$ preserves the degrees of homogeneous elements in $L$ since $a\in L^0$. So $ad_a$ is an endomorphism of $L$. Moreover it is an automorphism of $L$ since it has inverse $ad_{-a}$. 

The stability of the set  $\{ [x,x]=0, x\in L \} \subset L$  under $exp(ad_a)$ follows from $$exp(ad_a)[b,c]=[exp(ad_a)b,exp(ad_a)c].$$

\end{proof}

Now we see that the set $\{x\in L^1, [x,x]=0  \}$ is stable under action $exp(L^0):=\{ exp(ad_a), a\in L^0   \}$. We want to find the actions which preserve the Maurer-Cartan space, hence we now modify the action $exp(L^0)$.\\

Consider a new DGLA $(L',d',[\cdot , \cdot ]')$ such that $(L,d,[\cdot , \cdot])$ is an sub-DGLA of $L'$, i.e., $L'$ is a sub-complex and a sub-Lie algebra of $L$, by setting 
$$
(L')^i=L^i   \quad  \text{for}\quad  i\neq 1 ,\quad  (L')^1=L^1 \oplus \mathbb{K}d, d \, \,\text{ is  a formal symbol of degree 1 },
$$
with Lie bracket $[\cdot, \cdot]$ 
$$
[a+vd ,b+wd]'=[a,b]+vbd+(-1)^{deg(a)}wda.
$$
and differential 
$$
d'(a+vd)=da.
$$
In particular we have $[a,d]'=da$ if $a\in L^0$.  
\medskip

Define a  morphism $\phi: L^1 \longrightarrow (L')^1, a \mapsto a+d$. We see that the Maurer-Cartan equation $da+\frac{1}{2}[a,a]=0$ is equivalent to $[a+d,a+d]'=[\phi (a), \phi (a)]'=0$. With these notations, we are prepared to define gauge action.

\begin{definition}
For any  $a\in L^0,$ the gauge action of $a$, denoted by $exp(ad_a)$, acting on $L^1$ as 
$$
exp(ad_a) *b =\phi ^{-1} (exp(ad_a) (\phi (b))  \quad \quad \text{for any} \, \,  b\in L^1.
$$
Explicitly, 
\begin{equation} \label{g:def}
\begin{split}    
exp(ad_a)* b & =\sum_{n=0}^{\infty} \frac{(L_a)^n}{n !}(b)+\sum_{n=1}^{\infty} \frac{(L_a)^n}{n !}(\mathrm{d}) \\
& =\sum_{n=0}^{\infty} \frac{(L_a)^n}{n !}(b)+\sum_{n=1}^{\infty} \frac{(L_a)^{n-1}}{n !}(\mathrm{d} a) \\
& =b+\sum_{n=0}^{\infty} \frac{(L_a)^n}{(n+1) !}([a, b]-\mathrm{d} a) .
\end{split}
\end{equation}

Denote the set of gauge actions as $G_0=\{exp(ad_a), \, a\in L^0  \}$
\end{definition}
We say that $x,y\in L^1$ are gauge equivalent if there exists $a\in L^0$ such that $y=exp(ad_a)*x$.

\begin{corollary}
    Maurer-Cartan space $MC(L)\subset L^1$ of $L$ is stable under $G^0$.
\end{corollary}
\begin{proof}
    By Lemma \ref{L.4.8} we have $\{[x,x]'=0, x\in (L')^1\}$ is stable under $G^0$. Hence $$[x,x]'=0 \iff dx+\frac{1}{2}[x,x]=0$$ for any  $x\in L^1$.
\end{proof}

\medskip
\subsection{Artin Local Algebra}

We have seen that a formal deformation of complex structure is  a formal solution of  Maurer-Cartan equation. To state this in DGLA formalism, we will use Artin local ring as a generalization of polynomial ring.  We firstly review  the definition of Artin local ring (see \cite{M.Atiyah} for an  introduction). 

\smallskip 

An algebra over $\K$ (recall that $\K = \R$ or $\C$) is called local if it has only one maximal ideal. It is called Artin if it satisfies the descending chain condition, i.e. every descending chain of ideals is  stable. An Artin local algebra is an algebra which is both local and Artin. For simplicity we will also call an Artin local algebra just  Artin algebra without confusion. For an Artin local algebra $A$, we denote its unique maximal ideal as $\mathfrak{m}_A$. 

From the descending chain condition we see that $A$ is a finite dimensional vector space. Hence it has the form $ \frac{\K [t_1,...,t_n]}{I} $ for some ideal  $I \subset \K [t_1,...,t_n] $, and one can prove that $\mathfrak{m}_A$ nilpotent, i.e. $\mathfrak{m}_A^n=0$ for some $n$.

\begin{example}
 The ring  $\frac{\K[t_1,...t_n]}{R}$, which is a polynomial ring in finite variables quotient  by some ideal $R$    is an Artin local algebra over $\K$ with maximal ideal $(t_1,...,t_n)\subset \frac{\K[t_1,...t_n]}{R}$.
\end{example}
\medskip
\subsection{Deformation Functor and Formal Moduli Space}

The following notation is used: 

\smallskip

-- $\mathbf{Art}_{\mathbb{K}}$ the category of Artin local algebra with residue field $\mathbb{K}$;

-- \textbf{Set} the category of sets,

-- \textbf{Group}  the category of groups.

\medskip

Let us define deformation theory in DGLA formalism.

\begin{definition} (\cite{Manetti}, Ch. 6)

\vspace{3pt}

    Given a differential graded Lie algebra $L$, define three functors:
\begin{itemize}

    \item Exponential functor,
             $exp_L:\mathbf{Art}_{\mathbb{K}} \longrightarrow \textbf{Group}, \, A \mapsto exp(L^0 \otimes \mathfrak{m}_A) $             

    \item  Maurer-Cartan functor,
           $MC_L:\mathbf{Art}_{\mathbb{K}} \longrightarrow \textbf{Set}, \, A \mapsto \left\{x \in L^1 \otimes \mathfrak{m}_A \mid d x+\frac{1}{2}[x, x]=0\right\}  $

    \item Deformation functor,
          $Def_L:\mathbf{Art}_{\mathbb{K}} \longrightarrow \textbf{Set}, \, A \mapsto \frac{\left\{x \in L^1 \otimes \mathfrak{m}_A \mid d x+\frac{1}{2}[x, x]=0\right\}}{exp(L^0 \otimes \mathfrak{m}_A )}$

\end{itemize}
where ${exp(L^0 \otimes \mathfrak{m}_A )}$ is the gauge action of the DGLA $L\otimes \mathfrak{m}_A$. It is well-defined since $ \mathfrak{m}_A$ is nilpotent.

\begin{remark}
    If $A=\frac{\K\{t\}}{(t^{n+1})}$ then $Def_L(A)$ is the set of equivalence classes of solutions of MC equation up to order $n$.
\end{remark}

\end{definition}

To define formal moduli space, we now introduce the concept of complex space (or analytic variety). 

\begin{definition}
    Denote the algebra of complex convergent power series with variables $t_1,...,t_n$, for some $n>0$, by $\C\{t_1,...t_n\} $. An \textbf{analytic algebra} is a local algebra over $\C$ which is isomorphic to $ \frac{\C \{t_1,...t_n\}}{I}  $  for some   ideal  $I\subset (t_1,...,t_n)$.
    
    A morphism of analytic algebras is a morphism of local algebras, which is an algebra morphism sending the maximal ideal into maximal ideal.  Denote the category of analytic algebras by \textbf{AnaAlg}.
\end{definition}

Recall that a ringed space $(X,O_X)$ is a topological space $X$, equipped with a sheaf of rings $O_X$ on $X$ (called structure sheaf). A morphism between two ringed space $(X,O_X), (Y,O_Y)$ is a continuous map $f:X\longrightarrow Y$ equipped with a morphism of sheaves from $O_Y$ to $f_* {O_X}$.

\begin{definition}
Let $U$ be an open set of $\C^n$, $f_1,...,f_k$ be holomorphic functions on $U$, $O_U$ be the sheaf of  holomorphic functions on $U$. The ringed space $(V,O_V)$ is given by 
$$
V=\{x\in U | \, f_1(x)=...=f_m(x)=0, f_1,...f_k \in O_U\}, \,\, O_V=\frac{{{O_U}|_V}}{(f_1,...f_k)}. 
$$
We call $(V,O_V)$ a local model. 

A ringed space $(X,O_X)$ is called a \textbf{complex space} if it is locally isomorphic to a local model.   A morphism between two complex spaces is a morphism of ringed spaces.

\end{definition}

\begin{definition}
    Let $(X,O_X)$ be a complex space. Take a point $x\in X$. We introduce an equivalence relation $\sim$ on the restrictions of $(X,O_X)$ to  neighborhoods of $x$. Let $U,V$ be two arbitrary neighborhoods of $x$, we say $(U,O_U) \sim (V,O_V) $ if there exists an open set $W\subset U\cap V$ such that ${O_U}|_W \cong {O_V}|_W$. We call such an equivalence class a \textbf{germ of analytic varieties} at $x$ and the set of equivalence classes the stalk of $(X,O_X)$ at $x$. We denote the stalk of $(X,O_X)$ at $x$ by $O_{X,x}$ and the germ of analytic varieties of $(X,O_X)$ at $x$ by $(X,x)$.
\end{definition}

\begin{remark}
We define a morphism between two germs of analytic varieties as the germ of morphisms of ringed spaces. Precisely, let $(U,O_U)$, $(U',O_{U'})$ be the representatives of two germs of analytic varieties respectively and $f$, $g$ are two morphisms of ringed space between $(U,O_U)$ and  $(U',O_{U'})$. We say $f \sim g$ if $f=g$ on some open set $W\subset U\cap U'$ and call an equivalence class of the morphisms of ringed space a morphism of germs of analytic varieties. We denote the category of germs of analytic varieties by  \textbf{Germs of AnaVar}.

\end{remark}

From definition we see that the stalk $O_{X,x}$ of a complex space $(X,O_X)$ at any point $x\in X$ is isomorphic to $\frac{\C\{t_1,...t_k\} }{I}$ for some ideal $I$. Hence given a complex space, we can construct an analytic algebra by stalk at any point. Conversely, given an analytical algebra $R=\frac{\C\{t_1,...,t_d\}}{I}$, since $\C\{t_1,...t_d\}$ is Noetherian, the ideal $I$ is generated by finite many holomorphic functions $f_1,...f_m$. Hence $R$ is the stalk of some complex space $(X,O_X)$ which is defined by the zero locus of  $f_1,...,f_m$ at $0$. So we see that every analytic algebra $R$ generates a germ of analytic varieties. Actually we have 

\begin{theorem}(\cite{Fischer}, chapter 0) \label{t:4.18}
 We have an equivalence of categories  
 $$
 \textbf{AnaAlg}^{op} \longleftrightarrow  \textbf{Germs of AnaVar}
 $$
\end{theorem}

\begin{remark} \label{r:4.19}
    
Therefore often we denote a germ of analytic varieties $(X,x)$ by its corresponding analytic algebra  $O_{X,x}$. For example, we denote the germ of analytic varieties  $(\C^n,0)$ by $\C\{t_1,...,t_n\}$.  

\end{remark}

\begin{definition}
    The \textbf{Zariski tangent space} of a germ of  varieties  $(X,x)$ is 

   $$ T_x X= \left\{ \text{morphism of analytical algebras} : O_{X,x} \longrightarrow \frac{\C\{t\}}{(t^2)}  \right\}$$

\end{definition}

In particular we have

\begin{corollary}
The Zariski tangent space $T_0$ of the germ of analytic varieties generated by $\C\{t_1,...,t_d\}$ is isomorphic to the vector space $V=\text{Span}_\C (t_1,...,t_d)$
\end{corollary}

\begin{proof}
    Since $v$ is local algebra morphism, it maps maximum ideals into maximum ideals, i.e. $v(t_i)=v_i x$ for some $v_i\in \C$, we can determine a local algebra morphisms $v$ from $\C\{t_1,...t_d\}$ to $\frac{\C\{x\}}{(x^2)}$ by specifying $v(t_i)$.  We define a linear map $\phi: T_0 \longrightarrow V, v\mapsto \sum_i v_i t_i$. It is clear that $\phi$ is an isomorphism of $\C$ vector spaces.

\end{proof}

\begin{definition}
  If there exist an Artin algebra $R$ such that  $Def_L$ is natural isomorphic to the functor $Hom(R,\cdot)$ (or  $Def_L \cong h_R$),  we say  $Def_L$ is representable by $R$ and call the germ of analytic varieties  $R$ the \textbf{moduli space} of $L$. 
\end{definition}

It turns out that the requirement of the input of deformation functor to be an Artin algebra is too restrictive, the reason is that a formal solution of Maurer-Cartan equation is not required to be convergent. To include formal solutions of Maurer-Cartan equation we need  the algebra of formal power series, which is not an Artin algebra, but a pro-Artin algebra, i.e., a completion (or projective limit) of Artin algebra. Recall that a local complete Noetherian algebra  $\widehat{R}$ with maximal ideal $\mathfrak{m}_R$ is the completion of a local Artin algebra $R$ with maximal ideal $\mathfrak{m}_R$, i.e.  $\widehat{R} \cong \underset{i}{\varprojlim} \frac{R}{\mathfrak{m}^i} $. Specifically any $\K$ pro-Artin algebra is the quotient of an algebra of formal power series. In particular, if $R=\K[t]$ with maximal ideal $(t)$, its completion $\widehat{R}$ would be $\K [[t]]$, the formal power series in $t$.  \\
\begin{definition}\label{d: formal deformation}
 We denote the set  $\left\{ x\in L\otimes \mathfrak{m}_{\K [[t]]}|dx+\frac{1}{2}[x,x]=0 \right\} / G^0$ as $Def_L(\K[[t]])$ and call an element in it a \textbf{formal deformation} . 
\end{definition}

\begin{definition}
    Let  $L$ be a DGLA, its deformation functor $Def_L$ is called \textbf{pro-representable} if $Def_L \cong h_{\widehat{R}}$ for some pro-Artin algebra $\widehat R$. We call $\widehat{R}$ the \textbf{formal moduli space} of $L$.
\end{definition}

\medskip

\subsection{Smoothness of Deformation Functor and Obstruction}
Let $A_k$ stand for $\frac{\K \{t\}}{(t^{k+1})}$ for any $k>0$. 
Given a differential graded Lie algebra $L$ we call an element in $Def_L \left(A_k \right)$ a \textbf{$k$ order deformation}. In particular, we call the first order deformation $Def_L\left( A_1  \right)$ the \textbf{tangent space} of $Def_L$ and write it  as $TDef_L$.  From definition we see that a $k$ order deformation is an equivalence class of  solutions of Maurer-Cartan equation up to order $k$.

To find a formal solution of Maurer-Cartan equation, we  want to iteratively extend a $k$ order solutions to $k+1$ order solution as in Corollary \ref{c:6.11}. 

\begin{definition}[Smoothness]

\hspace{3pt}

Let $L$ be a DGLA.  The deformation functor $Def_L$ is \textbf{smooth} (or  the DGLA $L$ is unobstructed) if every $k$ order deformation can be extended to a $k+1$ order deformation. Precisely, $Def_L$ is smooth if  for any element $\phi_k \in MC_L(A_k)$ there exists an element $\phi_{k+1}\in MC_L(A_{k+1})$ such that $${\phi_{k+1}}=\phi_k \in MC_L(A_k).$$ 
   \end{definition}

\medskip

In particular, if $Def_L$ is smooth, then every first order deformation can be extended to a formal solution of Maurer-Cartan equation.

\begin{theorem}
    If $H^2(L)=0$, then $Def_L$ is smooth.
\end{theorem}
\begin{proof}
    The proof is essentially the same as in Theorem \ref{T:3.13}. Suppose we have a $k$ order solution $S_k = \sum^{k}\limits_{i=1} t^i \phi_i \in L^1\otimes \mathfrak{m}_{A_k}$. To extend it to a solution of order $k+1$, we want to find $\phi_{k+1}$ such that $S_{k+1}=S_k+t^{k+1}\phi_{k+1}$ satisfies MC equation up to order $k+1$. Hence, we need to solve 
    $$
    \Bar{\partial}(S_k+ t^{k+1} \phi_{k+1}) +\frac{1}{2}[S_k+ t^{k+1} \phi_{k+1},S_k+ t^{k+1} \phi_{k+1}]=0.
    $$
The vanishing of coefficients with $t^{k+1}$ leads to the equation

\begin{equation}\label{eq:ex}
\frac{1}{2}\sum_{i+j=k+1 \atop  0<i<j<k+1}[\phi_i,\phi_j]=\Bar{\partial}(\phi_{k+1}).
\end{equation}
 
Since $H^{0,2}(L)=0$, every closed element in $L^2$ is exact. By  Jacobi identity, the left side of Eq.4.2 is closed, hence there exist $\phi_{k+1}$ solving the equation, i.e., we can extend any $k$ order solution to $k+1$ order. Hence $Def_L$ is smooth.  
\end{proof}
Therefore we say that first order deformations in $L$ live in $H^1(L)$ and obstructions live in $H^2(L)$.

\begin{lemma} \label{l:7.23}
    If $L$ is pro-representable by $\K[[t_1,...,t_h]]$, then $Def_L$ is smooth.
\end{lemma}
\begin{proof}
    To prove this, we first show that we can extend any element in  $h_{K[[t_1,...,t_h]]}(  A_{k-1})$ to an element in $h_{\K[[t_1,...,t_h]]}(A_{k})$. Then since $h_{K[[t_1,...,t_h]]}\cong Def_L$ as functors, we can extend any element in $Def( A_{k-1})$ to $Def( A_{k})$, which by definition gives the smoothness of $Def_L$. 
    
    Given  $$\alpha\in h_{K[[t_1,...,t_h]]}(A_{k-1})=Hom(\K[[t_1,...,t_h]],A_{k-1}),$$  where $\alpha$ is determined by $$\alpha_i:=\alpha(t_i)\in \mathfrak{m}_{A_{k-1}},$$ where  $ i=1,...,h$. Since the canonical projection $\pi:A_{k} \longrightarrow A_{k-1}$ is surjective, we can find  $\psi_i\in {A_{k}}$ such that $\pi(\psi_i)=\alpha_i$. Then the morphism $\psi \in Hom(\K[[t_1,...,t_h]],A_{k})$ determined by $\psi(t_i)=\psi_i$ is an extension of $\alpha \in Hom(\K[[t_1,...,t_h]],A_{k-1})$.

\end{proof}

\begin{proposition}\label{p:4.25}
       The tangent space of $Def_L$, $TDef_L$ has a natural $\K$ vector space structure. It is isomorphic to  $H^1(L)$  as $\K$ vector space. 
\end{proposition}
\begin{proof}
    The proof is  essentially the same as the proof in Corollary \ref{c:6.11}. An element in $TDef_L$ is a (actually an equivalence class of) closed element in $L^1$, hence  $TDef_L$ has a natural $\K$ vector space structure. After neglecting terms higher than two orders, two elements $\phi, \phi'$ in $TDef_L$ are gauge equivalent if and only if there are some $a\in L^0$ such that $\phi'=\phi-da$. 
\end{proof}

Now we see that the tangent space of a deformation functor has a natural $\K$ vector space structure as the first cohomology group of the DGLA . 

\begin{definition}
    A morphism of DGLAs  $f: L\longrightarrow L'$ from $L$ to $L'$ is called a quasi-isomorphism from $L$ to $L'$ if it induces an isomorphism from $L$ to $L'$. We say two DGLAs are  quasi-isomorphic if such a quasi-isomorphism between them exists.
    
    \end{definition}

\begin{corollary}
    Given two DGLAs $L,L'$, if they  are quasi-isomorphic, we have $H^k(L) \cong  H^k(L')$ for any $k$. In particular, we have an isomorphism between  tangent spaces of two deformation functors $Def_L$, $Def_{L'}$. 
\end{corollary}

\begin{proposition}
       If $Def_L$ is (pro-)representable by some  (pro-)Artin algebra $R$, write $T R$ as the Zariski tangent space of  germ of analytic varieties $R$, then  $TDef_L \cong TR$ as $\K$ vector spaces.
\end{proposition}
\begin{proof}
By definition, $TDef_L= Def_L\left( \frac{\K \{t\}}{(t^2)}  \right)$, and the Zariski tangent space of $R$ is $Hom(R, \frac{\K \{t\}}{(t^2)})$.  Consider $f: \frac{\K \{t\}}{(t^2)}\longrightarrow  \frac{\K \{t\}}{(t^2)}, f(x)=rx, k\in \K$. Since $Def_L \cong h_R$ with some natural isomorphism $\phi$  we have the commutative diagram  

$$\begin{tikzcd}
	{ Def_L\left( \frac{\mathbb{K} \{t\}}{(t^2)}  \right)} && {Hom(R, \frac{\mathbb{K} \{t\}}{(t^2)})} \\
	\\
	{ Def_L\left( \frac{\mathbb{K} \{t\}}{(t^2)}  \right)} && {Hom(R, \frac{\mathbb{K} \{t\}}{(t^2)})}
	\arrow["{\phi}", from=1-1, to=1-3]
	\arrow["{\times r}"', from=1-1, to=3-1]
	\arrow["{\times r}", from=1-3, to=3-3]
	\arrow["{\phi}"', from=3-1, to=3-3]
\end{tikzcd}
$$
where $\phi$ is bijection. Therefore $TDef_L$ is bijective to $T(0,R)$ as set and for any $a\in Def_L\left( \frac{\mathbb{K} \{t\}}{(t^2)}  \right)$,  we have $\phi(ra)=r\phi(a)$. 

When we consider first order deformation, $Def_L\left( \frac{\mathbb{K} \{t\}}{(t^2)} \oplus \frac{\mathbb{K} \{t\}}{(t^2)} \right) $ can be decomposed as \\ $ Def_L\left( \frac{\mathbb{K} \{t\}}{(t^2)}  \right) \oplus  Def_L\left( \frac{\mathbb{K} \{t\}}{(t^2)}  \right)$ since 

\begin{equation*}
    \begin{split}
    Def_L\left( \frac{\mathbb{K} \{t\}}{(t^2)} \oplus \frac{\mathbb{K} \{t\}}{(t^2)} \right) &= \frac{ \left\{  \phi=\phi_1 t_1+ \phi_2 t_2, d\phi =0 | \phi_1,\phi_2 \in L_1 \right\}}{\text{gauge action}} \\
    &=\frac{ \left\{  \phi=\phi_1 t_1, d\phi=0 | \phi_1 \in L_1 \right\}}{\text{gauge action}}\oplus \frac{ \left\{  \phi=\phi_2 t_2, d\phi=0 | \phi_2 \in L_1 \right\}}{\text{gauge action}}.
\end{split}
\end{equation*}

  Also since $Hom(R, \cdot)$ preserves direct sums, we have $$Hom(R,\frac{\mathbb{K} \{\}}{(t^2)} \oplus \frac{\mathbb{K} \{t\}}{(t^2)} )= Hom(R,\frac{\mathbb{K} \{t\}}{(t^2)}) \oplus Hom(R,\frac{\mathbb{K} \{t\}}{(t^2)}).$$ Because $\phi$ is a natural isomorphism, the diagram 

$$
\begin{tikzcd}
	{ Def_L\left( \frac{\mathbb{K} \{t\}}{(t^2)}  \right) \oplus  Def_L\left( \frac{\mathbb{K} \{t\}}{(t^2)}  \right)} && {Hom(R,\frac{\mathbb{K} \{t\}}{(t^2)}) \oplus Hom(R,\frac{\mathbb{K} \{t\}}{(t^2)})} \\
	{Def_L\left( \frac{\mathbb{K} \{t\}}{(t^2)}  \right)} && {Hom(R,\frac{\mathbb{K} \{t\}}{(t^2)})}
	\arrow["{\pi_i}"', from=1-1, to=2-1]
	\arrow["\phi", from=1-1, to=1-3]
	\arrow["\phi"', from=2-1, to=2-3]
	\arrow["{\pi_i}", from=1-3, to=2-3]
\end{tikzcd}
$$
commutes, where $\pi_i$ is the projection to the $i$-th component. Hence we have that for any $(a,b)\in  Def_L\left( \frac{\mathbb{K} \{t\}}{(t^2)} \oplus \frac{\mathbb{K} \{t\}}{(t^2)} \right)$,  $\phi(a,b)=(\phi(a),\phi(b))$. Now we define function 

\begin{equation*}
    \begin{split}
        +: Def_L\left( \frac{\mathbb{K} \{t\}}{(t^2)}  \right) \oplus  Def_L\left( \frac{\mathbb{K} \{t\}}{(t^2)}  \right) &\longrightarrow Def_L\left( \frac{\mathbb{K} \{t\}}{(t^2)}  \right) \\
        (a,b)&\mapsto a+b
    \end{split}
\end{equation*} 
Again by naturality of $\phi$ we have commutative diagram 
  $$
  \begin{tikzcd}
	{ Def_L\left( \frac{\mathbb{K} \{t\}}{(t^2)}  \right) \oplus  Def_L\left( \frac{\mathbb{K} \{t\}}{(t^2)}  \right)} && {Hom(R,\frac{\mathbb{K} \{t\}}{(t^2)}) \oplus Hom(R,\frac{\mathbb{K} \{t\}}{(t^2)})} \\
	{Def_L\left( \frac{\mathbb{K} \{t\}}{(t^2)}  \right)} && {Hom(R,\frac{\mathbb{K} \{t\}}{(t^2)})}
	\arrow["{+}"', from=1-1, to=2-1]
	\arrow["\phi", from=1-1, to=1-3]
	\arrow["\phi"', from=2-1, to=2-3]
	\arrow["{+}", from=1-3, to=2-3]
\end{tikzcd}
  $$
  which tell us for any $(a,b)\in  Def_L\left( \frac{\mathbb{K} \{t\}}{(t^2)}  \right) \oplus  Def_L\left( \frac{\mathbb{K} \{t\}}{(x^2)}  \right)$, we have $\phi(a+b)=\phi(a)+\phi(b)$. Hence $\phi$ is a morphism between vector spaces. 
\end{proof}

Therefore if a differential graded Lie algebra $L$ is pro-represenatable by $R$ (or it has formal moduli space $R$), we can identify the tangent space of $Def_L$ and the (Zariski) tangent space of its formal moduli space $R$. In particular,

\begin{corollary}
        If $L$ is (pro-)representable by $\K[[t_1,...,t_h]]$, then $h=dim H^1(L)$.
\end{corollary}
\begin{proof}
    By definition, the Zariski tangent space of $\K[[t_1,...,t_h]]$ is the space of local algebra morphisms from $\K[[t_1,...,t_h]]$ to $ \frac{\K\{t\}}{(t^2)}$. Define $\alpha_i\in Hom(\K[[t_1,...,t_h]],\frac{\K\{t\}}{(t^2)})$ for $i=1,...,h$, where  $\alpha_i(t_i)=t $ and $ \alpha(t_j)=0$ if $j\neq i$. It is easy to see that $\alpha_i$ form a basis of $Hom(\K[[t_1,...,t_h]],\frac{\K\{t\}}{(t^2)})$. Hence dim $TDef_L=h= $dim $H^1(L)$.
\end{proof}

By Lemma \ref{p:4.25} we see that the tangent space of a deformation functor has a natural $\K$-vector space structure as the first cohomology group of the DGLA .  When two DGLA $L,L'$ are \textbf{quasi-isomorphic}, i.e. the morphism betweem them induces an isomoprhism between their cohomology groups,  we have $H^k(L) \cong  H^k(L')$ for any $k$. In particular, we have an isomorphism between the tangent space of two deformation functors $Def_L,\, Def_{L'}$. Moreover we have

\begin{theorem} \label{T:7.28} 
    If DGLA $L,L'$ are quasi-isomorphic, then their deformation functors $Def_L, Def_{L'}$ are isomorphic.
\end{theorem}
\begin{proof}
 See \cite{Manetti}, chapter 6.
\end{proof}

\begin{corollary} \label{c:7.29}
    Given a DGLA $L$, if $L$ is quasi-isomorphic to an abelian DGLA $L'$, then $Def_L$ is pro-representable and smooth.
\end{corollary}
\begin{proof}
     Since $L'$ is abelian, for any Artin algebra $A$
$$
Def_{L'}(A)= \frac{\{x\in {L'}^1 \otimes \mathfrak{m}_A | dx=0 \}}{\{x\sim y  \text{ if } x-y=da \text{ for some } a\in {L'}^0\}} = H^1(L')\otimes \mathfrak{m}_A.
$$
 Recall that we assume every cohomology group is finite dimensional. Let $(t_1,...,t_h)$ be a basis of $H^1(L')$, then any element in  $Def_{L'}(A)$ can be written as $\sum\limits_i(t_i\otimes a_i)$ where $a_i\in \mathfrak{m}_A$. Let $\K[[t_1,...,t_h]]$ be the algebra of formal power series with variables $t_1,...,t_h$. Define a family of morphisms 
 \begin{align*}
 f_A:Def_{L'}(A)&\longrightarrow Hom(\K[[t_1,...,t_h]],A), \\
 f_A(\sum\limits_i(t_i\otimes a_i))&=\phi \in \K[[t_1,...,t_h]] 
  \end{align*}
 where $\phi(t_i)=a_i$ for any $i$ for any Artin algebra $A$. We show now that $f$ is an isomorphism (of functors) between $Def_{L'}$ and $Hom(\K[[t_1,...,t_h]])$.

 For any Artin algebra $A$ and any element $\psi\in Hom(\K[[t_1,...,t_h]],A)$, there exists an element $\sum\limits_i t_i \otimes \psi(t_i) \in Def_{L'}(A)$ such that $f_A(\sum\limits_i t_i \otimes \psi(t_i))=\psi$. Hence we proved that $f_A$ is surjective. 

 To prove injectivity of $f_A$, let $\sum\limits_i t_i\otimes a_i,\sum\limits_i t_i \otimes a'_i \in Def_{L'}(A)$. Denote  $f_A(\sum\limits_i t_i\otimes a_i) 
$ by $\alpha$,  $ f_A(\sum\limits_i t_i\otimes a'_i)$ by $\alpha'$.  Suppose $\alpha=\alpha'$, then $\alpha(t_i)=a_i=\alpha'(t_i)=a'_i$. We have that $$\sum\limits_i t_i\otimes a_i=\sum\limits_i t_i \otimes a'_i.$$ Hence $f_A$ is injective.

To prove $f$ is a natural isomorphism between $Def_{L'}$ and $Hom(\K[[t_1,...,t_h]], \cdot)$, suppose there are  Artin algebras $A,B$ and a  morphism $p\in Hom(A,B)$. For any $\sum\limits_i t_i\otimes a_i \in Def_{L'}(A)$, denote $$p\circ f_A (\sum\limits_i t_i\otimes a_i ) \in {Hom(\mathbb{K}[[t_1,...,t_h]],B)}$$ by $\eta$ and $$f_B \circ(\text{id}\otimes p)(\sum\limits_i t_i\otimes a_i )$$ by $\eta'$. Then $\eta(t_i)=p(t_i)=\eta'(t_i)$, that is, the following diagram commutes. 

$$\begin{tikzcd}
	{Def_{L'}(A)\cong H^1(L')\otimes\mathfrak{m}_A} & {Hom(\mathbb{K}[[t_1,...,t_h]],A)} \\
	{Def_{L'}(B)\cong H^1(L')\otimes\mathfrak{m}_B} & {Hom(\mathbb{K}[[t_1,...,t_h]],B)}
	\arrow["{f_A}", from=1-1, to=1-2]
	\arrow["{\text{id}\otimes p}"', from=1-1, to=2-1]
	\arrow["{ p \circ}", from=1-2, to=2-2]
	\arrow["{f_B}"', from=2-1, to=2-2]
\end{tikzcd}$$

Therefore $Def_{L'} \cong h_{\K[[t_1,...,t_h]]}$.  By last theorem, since $L$ is quasi-isomorphic to $L'$,  we know that $Def_L \cong Def_{L'}$. Hence $Def_L \cong h_{\K[[t_1,...,t_h]]}$, $Def_L$ is pro-representable. By Lemma \ref{l:7.23} we see that $Def_L$ is smooth.
\end{proof}
\medskip
\subsection{Formal Deformation of Complex Structure}

In this section we illustrate the construction from Sect. 4.1-4.5 on  the problem of deformation of complex structures.
Recall that our DGLA is given by  the Kodaira--Spencer algebra of $M$: 
$$KS_M=\bigoplus_i \Omega^{0,i}(T^{1,0}M).$$ From Definition \ref{d:3.13}, a formal deformation of complex structure is a formal solution of Maurer-Cartan equation in $KS_M\otimes \mathfrak{m}_P$ where $P=\K [[t]]$. A $k$ order deformation is a solution of Maurer-Cartan equation in $KS_M \otimes \mathfrak{m}_{A}$, where $A=\frac{\K\{t\}}{(t^{k+1})}$. Hence the set of $k$ order deformations of complex structures is $MC_{KS_M}(\frac{\K\{t\}}{(t^{k+1})})$. 

In Sect.4, we  define $Def_L(A_k)$ as the space of $k$ order solutions modulo gauge action. 

In Sect.3 a formal deformation of complex structures is defined as a formal solution to MC equation. 
In order to make both approaches in Sect.3 and Sect. 4 coincide we identify isomorphic deformations of complex structures. 
Recall that two deformations of complex structure are isomorphic if one is the push-forward of the other by the flow of a vector field. That is, deformations of complex structures $J_{1,t}$, $J_{2,t}$ are isomorphic if there is a holomorphic vector field $X$ which generates a one parameter group of diffeomorphisms $I_t$, we have 
$$ J_{1,t} \circ dI_t=dI_t \circ J_{2,t}.$$

\medskip

For simplicity we consider deformations only on one dimensional complex manifold $(M,J)$, the case for arbitrary dimension is similar. 
We will see that two deformations of complex structure are isomorphic if and if if they are gauge equivalent. Hence after we identifying isomorphic deformations, we have the same definition of deformation as in last section.

Recall that given a DGLA $L$ and $a\in L^0$, $b\in L^1$, the gauge action of $a$ on $b$ is defined by 
\begin{equation}\label{eq:4.2}
   exp(ad_a)*b =b+\sum_{n=0}^{\infty} \frac{{(L_a)}^n}{(n+1) !}([a, b]-\mathrm{d} a) .
\end{equation}

\begin{theorem}
    Two formal deformations of complex structures $J_{t,1},J_{t,2}$ on $M$ with encoding map $\phi_{1,t},\phi_{2,t}$ are isomorphic if and only if $\phi_{1,t}, \phi_{2,t}$ are gauge equivalent in $KS_M$.
\end{theorem}

\begin{proof}
We first show that if $J_t$ is a trivial deformation, then its encoding map $\phi_t$ is gauge equivalent to zero in $KS_M$.

    Recall that a deformation $J_t$ is trivial if it is induced by an one-parameter group of diffeomorphisms $I_t$, that is, $J_t \circ dI_t=dI_t \circ J$. Hence the anti-holomorphic part of $J_t$ is the push forward of the anti-holomorphic part of $J$, that is, $T^{0,1}_t M$=$dI_t (T^{0,1}M)$. If we define a vector field $-X \in \Omega^{0} (T^{1,0}M)$ as the generator of $I_t$, i.e. $-X_p=\frac{d (I_t(p))}{d t}$, then we have $$\frac{d(dI_t \left( \frac{\partial}{\partial \Bar{z}} \right)}{dt}|_{t=0}=L_X\left( \frac{\partial}{\partial \Bar{z}} \right),$$ also $$\frac{d^k(dI_t \left( \frac{\partial}{\partial \Bar{z}} \right)}{dt^k}|_{t=0}=(L_X)^k\left( \frac{\partial}{\partial \Bar{z}} \right)$$ for any $k>1$. Now formally we have 
\begin{equation} \label{eq:dI}
    dI_t\left( \frac{\partial}{\partial \Bar{z}} \right)=  \frac{\partial}{\partial \Bar{z}} + tL_X\left( \frac{\partial}{\partial \Bar{z}} \right)+ \frac{1}{2!}t^2(L_X)^2\left( \frac{\partial}{\partial \Bar{z}} \right)+ .... ,
\end{equation}

which spans the anti-holomorphic part of $J_t$.  Since $X$ is a holomorphic vector field, we write it as $X=x\frac{\partial}{\partial z}$ where $x$ is a smooth function on $M$, then 
\begin{equation} \label{4.3}
L_X\left( \frac{\partial}{\partial \Bar{z}} \right)=[x\frac{\partial}{\partial z},\frac{\partial}{\partial \Bar{z}}]=-\frac{\partial x}{\partial \bar{z}}\frac{\partial}{\partial z}  .  
\end{equation}

Hence $L_X (\frac{\partial}{\partial \Bar{z}})$ has no term with $\frac{\partial}{\partial \Bar{z}}$. By Theorem \ref{T: 2.14.}, $[X,L_X(\frac{\partial}{\partial \Bar{z}})]$ has also no term in $\frac{\partial}{\partial \Bar{z}}$, same for $(L_X)^n(\frac{\partial}{\partial \Bar{z}}) $ for any $n>1$. Therefore we have $\pi^{0,1}(dI_t (\frac{\partial}{\partial \Bar{z}}))= \frac{\partial}{\partial \Bar{z}}$. Now we get $$(\pi^{0,1}|_{T^{0,1}_t})^{-1}\left( \frac{\partial}{\partial \Bar{z}}  \right) = dI_t \left( \frac{\partial}{\partial \Bar{z}}\right).$$

Hence  we have 
\begin{equation} \label{4.2}
\begin{split}
        \phi_{t}\left( \frac{\partial}{\partial \Bar{z}} \right)&=\pi^{1,0} \circ (\pi^{0,1}|_{T^{0,1}_t})^{-1} \left(\frac{\partial}{\partial \Bar{z}}\right) = \pi^{1,0}\circ \, dI_t \left( \frac{\partial}{\partial \Bar{z}}\right) \\
        &=\sum^\infty_{n=0} \frac{{(L_{tX})}^n}{(n+1)!} [tX, \frac{\partial}{\partial \Bar{z}}] = \sum^\infty_{n=0} \frac{{(L_{tX})}^n}{(n+1)!} (-\Bar{\partial}{(tX)}\left( \frac{\partial}{\partial \Bar{z}} \right)).  
        \end{split}
\end{equation}

By definition of gauge action in Eq \ref{eq:4.2}, we see that $\phi_t$ is gauge equivalent to $0$ if we let $a=tX$. Therefore, if a deformation of  complex structure $\phi_t$ is trivial then it is gauge equivalent to $0$ where $a=tX$ for some holomorphic tangent field $X\in \Omega^0(T^{1,0}M)$  in DGLA $KS_M$. 

Conversely, if $\phi_t$ is gauge equivalent to $0$ by gauge action $a=tX$ for some holomorphic vector field $X$, that is, 
\begin{equation}
    \phi_{t}\left( \frac{\partial}{\partial \Bar{z}} \right)= \sum^\infty_{n=0} \frac{{(L_{tX})}^n}{(n+1)!} (-\Bar{\partial}{(tX)}\left( \frac{\partial}{\partial \Bar{z}} \right))=\sum^\infty_{n=0} \frac{{(L_{tX})}^n}{(n+1)!} [tX, \frac{\partial}{\partial \Bar{z}}]
\end{equation}

Then we have 
\begin{equation}
\begin{split}
    T_t^{0,1}M&=(\text{id}+\phi_t)\left(\frac{\partial}{\partial \Bar{z}} \right) \\
    &= \frac{\partial}{\partial \Bar{z}} + tL_X\left( \frac{\partial}{\partial \Bar{z}} \right)+ \frac{1}{2!}t^2(L_X)^2\left( \frac{\partial}{\partial \Bar{z}} \right)+ .... \\
    &=dI_t \left(\frac{\partial}{\partial \Bar{z}} \right)
\end{split}    
\end{equation}
where $I_t$ is the flow generated by $-X$. Therefore we see that 
$$
T_t^{0,1}M=dI_t(T^{0,1}M).
$$
Hence the deformation is generated by $I_t$, i.e. $J_t$ is a trivial deformation.

\smallskip

Now we turn to the general case. Given  two deformations $J_{1,t}, J_{2,t}$, denote their (anti-)holomorphic parts in local coordinates by $(\frac{\partial}{\partial \Bar{z_{i,t}}})  \frac{\partial}{\partial {z_{i,t}}}$ and their encoding map by $\phi_{i,t}$ for $i=1,2$ respectively.  By definition of $\phi_{i,t}$ (see Definition \ref{R:3.5}) we have 

$$
\frac{\partial}{\partial \Bar{z_{i,t}}}= \frac{\partial}{\partial \Bar{z}}+ \phi_{i,t} \left( \frac{\partial}{\partial \Bar{z}} \right).
$$

To say  $J_{t,1}$ are $J_{t,2}$ are isomorphic is equivalent to say that $\frac{\partial}{\partial \Bar{z_{2,t}}}=dI_t\left( \frac{\partial}{\partial \Bar{z_{1,t}}} \right)$ for some one parameter group of diffeomorphisms $I_t$ on $M$.  If $\frac{\partial}{\partial \Bar{z_{2,t}}}=dI_t\left( \frac{\partial}{\partial \Bar{z_{1,t}}} \right)$ for some one parameter group of diffeomorphisms on $M$. 
Denote the generator vector field of the one parameter group by $-X$. We have
 \begin{equation}
 \begin{split}
     \phi_{2,t}\left( \frac{\partial}{\partial \Bar{z}} \right)-\phi_{1,t}\left( \frac{\partial}{\partial \Bar{z}} \right)&= \frac{\partial}{\partial \Bar{z_{2,t}}}-\frac{\partial}{\partial \Bar{z}}-\phi_{1,t}\left(\frac{\partial}{\partial \Bar{z}} \right) \\
&= dI_t\left( \frac{\partial}{\partial \Bar{z_{1,t}}} \right)-\frac{\partial}{\partial \Bar{z}}-\phi_{1,t}\left( \frac{\partial}{\partial \Bar{z}} \right) \\
&=dI_t \left(\frac{\partial}{\partial \Bar{z}} + \phi_{1,t} \left( \frac{\partial}{\partial \Bar{z}}\right)\right)-\frac{\partial}{\partial \Bar{z}}-\phi_{1,t}\left(\frac{\partial}{\partial \Bar{z}} \right). \\
\end{split}
 \end{equation}

Then

\begin{equation} \label{eq:4.6}
     \phi_{2,t}-\phi_{1,t}=dI_t +dI_t \circ \phi_{1,t} -\text{id}- \phi_{1,t}
\end{equation}

We have seen that $dI_t- \text{id}=\sum_{n=1}^\infty \frac{(L_{tX})^n}{n!}(-\bp ({tX}))$.  Similarly  we have $$dI_t \circ \phi_{1,t}-\phi_{1,t}=\sum_{n=1}^\infty  \frac{(L_{tX})^n}{n!}(\phi_{1,t}).$$ Substitute them into Eq.\ref{eq:4.6}, we have 
  
 \begin{equation}
     \phi_{2,t}=\phi_{1,t}+\sum_{n=0}^{\infty} \frac{({L_a})^n}{(n+1) !}([c, \phi_{1,t}]-\bp a),
 \end{equation}
 where $a=tX$. Therefore if $\phi_{1,t}$ and $\phi_{2,t}$ are isomorphic as deformations of complex structure, they are equivalent under gauge action of $a=tX$ in DGLA $KS_M$. The converse direction is similar as in the trivial deformation case.
\end{proof}

Therefore if we identify isomorphic deformations of a complex structure, then a formal deformation of complex structure on complex manfiold $M$ is an element in $Def_{KS_M}(\K[[t]])$.

\begin{definition}
Two formal deformations of complex structures $\phi_1[[t]], \phi_2[[t]]\in  Def_{KS_M}(\K[[t]])$ are called isomorphic if  there is a formal holomorphic vector field $X\in \Omega^{0}(T^{1,0}M)[[t]]$ such that $\phi_1[[t]]$ is  equivalent to $\phi_2[[t]]$ under gauge action of $X$. 
\end{definition}

Now we can see that after identifying isomorphic formal deformations of complex structure,  $Def_{KS_M}(\K[[t]])$ is the set of formal deformations of complex structures. 

In this spirit we define $k$ order deformation of complex structure.
 
\begin{definition}
A $k$ order deformation of complex structure is an element in $Def_{KS_M}(A_k)$.  
\end{definition}

\begin{example} {\cite{Manetti}}
Given a one dimensional complex torus $\mathbb{T}$,  its Kodaira-Spencer algebra is $L=\Omega^{0,*}(T^{1,0}\mathbb{T})$. Clearly $\Omega^{0,q}(T^{1,0}\mathbb{T})=0$ if $q\neq 0,1$ since the (complex) dimension of the torus is $1$. Then $L$ can be written as   

\begin{tikzcd}
	{} & 0 & 0 & {\Omega^0(T^{1,0}\mathbb{T})} & {\Omega^{0,1}(T^{1,0}\mathbb{T})} & 0 & 0 &
	\\
	& {}
	\arrow[from=1-2, to=1-3]
	\arrow[from=1-1, to=1-2]
	\arrow[from=1-3, to=1-4]
	\arrow[from=1-4, to=1-5]
	\arrow[from=1-5, to=1-6]
	\arrow[from=1-6, to=1-7]
\end{tikzcd} \\
Consider DGLA $H$ defined as 
$$
\begin{tikzcd}
	{} & 0 & 0 & {H^0(T^{1,0}\mathbb{T})} & {H^{0,1}(T^{1,0}\mathbb{T})} & 0 & 0 &  \\
	& {}
	\arrow[from=1-2, to=1-3]
	\arrow[from=1-1, to=1-2]
	\arrow[from=1-3, to=1-4]
	\arrow[from=1-4, to=1-5]
	\arrow[from=1-5, to=1-6]
	\arrow[from=1-6, to=1-7]
\end{tikzcd}
$$
 with trivial differential and trivial Lie bracket. 

One can prove that $L$ and $H$ are quasi-isomorphic as DGLAs (See \cite{Manetti} for details. In this case we say that $L$ is formal),  then $Def_L \cong Def_H$. By Lemma \ref{c:7.29}, $Def_L$ is pro-representable and smooth and $\mathbb{T}$ is unobstructed.

\end{example}

\medskip

\subsection{Unobstructness of Compact Calabi-Yau Manifold}
Recall that in Remark \ref{compact},  we have assumed that all manifolds in thesis are compact.  Moreover, we assume in this thesis that all cohomology groups are finite dimensional, which is true in case of compact Kähler manifolds.

\begin{definition}
    An $n$ dimensional Kähler manifold $M$ is called Calabi-Yau if it has a nowhere-vanishing holomorphic top form  $\omega\in  \Omega^{n,0}(M)$.
\end{definition}

G. Tian (\cite{G.Tian}) and Andrey N. Todorov (\cite{A.Todorov}) proved that if $M$ is a Calabi-Yau manifold, the deformation functor of its Kodaira-Spencer algebra is smooth, in other words, $M$ is unobstructed. We will follow the strategy in \cite{BK} and  prove this fact via deformation theory using DGLA.  We will show that the Kodaira-Spencer algebra $KS_M$ of a Calabi-Yau manifold is quasi-isomorphic to an abelian DGLA, then by Corollary \ref{c:7.29} $KS_M$ is smooth and $M$ is unobstructed. To the end of this, we firstly extend Kodaira-Spencer algebra. The first step is to extend vector field $T^{1,0}M$ to multi-vector field $\wedge^pT^{1,0}M$ on $M$ .

\begin{definition}
     A degree $k$ multi-vector field on a smooth manifold is a smooth section of degree $(k+1)$ exterior power of tangent bundle. Explicitly, on a smooth manifold $M$, a degree $(k-1)$ multi-vector field  $V\in \Gamma (\wedge^k TM)$ in local coordinates $\{x_i\}_{i=1}^n$is 

$$
V=\sum_{i_1 ,\ldots , i_k=1}^n X^{i_1 \ldots i_k}(x) \frac{\partial}{\partial x_{i_1}} \wedge \ldots \wedge \frac{\partial}{\partial x_{i_k}},
$$   
In particular, we denote $C^\infty (M)$ as a degree $(-1)$ multi-vector field. We will denote the degree of a multi-vector field $X$ by $\Tilde{X}$.
\end{definition}

\begin{definition}[Schouten–Nijenhuis bracket]
    Consider a smooth manfiold $M$, for $v_1,...v_i ,w_1,...\\,w_j\in \Gamma (TM)$, $f \in C^\infty (M)$, we define Schouten–Nijenhuis bracket as a bilinear map  
$$
[\cdot, \cdot]_S: \Gamma (\wedge^i TM) \otimes \Gamma (\wedge^j TM) \rightarrow \Gamma (\wedge^{i+j-1} TM)
$$
for $v_n, w_m \in \Gamma(TM), n=1,...,i, m=1,...,j$ and 
$$ 
[v_1 \wedge \ldots \wedge v_i, w_1 \wedge \ldots \wedge w_j]_S= \sum_{p=1}^i \sum_{q=1}^j(-1)^{p+q}\left[v_p, w_q\right] \wedge v_1 \wedge \ldots \wedge \widehat{v}_p \wedge \ldots \wedge v_i \wedge w_1 \wedge \ldots \wedge \widehat{w}_q \wedge \ldots \wedge w_j.
$$
and for any $f \in C^\infty (M)$,
$$
[v_1 \wedge \ldots \wedge v_i, f]_S=-\sum_{p=1}^i(-1)^p {L}_{v_p}(f) v_1 \wedge \ldots \wedge \widehat{v}_p \wedge \ldots \wedge v_i,
$$

$$
[f,w_1 \wedge \ldots \wedge w_j]_S=-\sum_{q=1}^i(-1)^q {L}_{w_q}(f) w_1 \wedge \ldots \wedge \widehat{w}_p \wedge \ldots \wedge w_j
$$
\end{definition}
It is easy to see that when the degree of multi-vector fields are $0$ or $-1$, Schouten–Nijenhuis bracket is  normal Lie bracket, so without causing confusion we will omit the $S$ in $[\cdot, \cdot]_S$ and write Schouten-Nijenhuis brackt as normal Lie bracket $[\cdot, \cdot]$. 
\begin{proposition}
Let $X,Y$ be multi-vector fields on a smooth manifold $M$.  Denote degree of any multi-vector field $F$ by $\Tilde{F}$. Then we have
$$\begin{aligned}   
  &  {1.[X, Y]=-(-1)^{\Tilde{X} \Tilde{Y}}[Y, X] }\\
  &
 {2.[X, Y \wedge Z]=[X, Y] \wedge Z+(-1)^{(\Tilde{X}+1) \Tilde{Y}} Y \wedge[X, Z]} \\
 &
 {3.(-1)^{\Tilde{X}\Tilde{Z}}[X,[Y,Z]]+(-1)^{\Tilde{X}\Tilde{Y}}[Y,[Z,X]]+(-1)^{\Tilde{Y}\Tilde{Z}}[Z,[X,Y]]=0  \quad \text{(graded Jacobi identity)}}
\end{aligned}
$$
\end{proposition}
\begin{proof}
This can be directly computed by the definition of Schouten-Nijenhuis bracket. We omit it here.
    
\end{proof}

We now give the definition the extended  Kodaira-Spencer algebra.

\begin{definition}
    Let $M$ be an $n$ dimensional complex manifold. Let
$$
{KS_M}:=\bigoplus_k {KS_M}^k,\quad  {KS_M}^k=:\bigoplus_{q-p+1=k} \Omega^{0, q}(X, \wedge^p T^{0,1}M),
$$
be a $\mathbb{Z}$-graded algebra where the degree of a homogeneous element in $\Omega^{0, q}(\wedge^p T^{1,0}M)$ is $k=q-p+1$. Introduce the $\text { multi-index } I=(i_1, \ldots, i_k)$ and $J=(j_1, \ldots, j_k)$ where $1 \leq i_1<\ldots<i_k \leq n$,$1 \leq j_1<\ldots<j_k \leq n$. Also we use $d r_I$ to represent $  dr_{i_1} \wedge \ldots \wedge d r_{i_k}$,  $\frac{\partial}{\partial z^I}$ to represent  $\frac{\partial}{\partial z^{i_1}} \wedge \ldots \wedge \frac{\partial}{\partial z^{i_k}}$. For $\alpha=\sum\limits_{I,J} \alpha^{I J} d \bar{z}_I \frac{\partial}{\partial z_J}$ and $\beta =\sum\limits_{K,L} \beta^{K L} d \bar{z}_K \frac{\partial}{\partial z_L}$. The extended  Schouten–Nijenhuis bracket  $$[\cdot,\cdot]: \Omega^{0, q}(\wedge^p{T^{1,0}}M) \times \Omega^{0, q^\prime}(\wedge^{p^\prime} T^{1,0}M )\longrightarrow \Omega^{0, q+q^\prime}(\wedge^{p+p^\prime-1} T^{1,0}M)$$ is defined as 

$$
{[\alpha, \beta] }  = \sum_{I,J,K,L}[\alpha_{I J} d \bar{z}_I \otimes \frac{\partial}{\partial z^J}, \beta_{K L} d \bar{z}_K \otimes \frac{\partial}{\partial z^L}]=\sum_{I,J,K,L}(-1)^{(|K|)(|J-1|)}d \bar{z}_I \wedge d \bar{z}_K\wedge[\alpha^{I J} \frac{\partial}{\partial z^J}, \beta^{K L} \frac{\partial}{\partial z^L}],
$$ 
We extend differential $\bar{\partial} $ on $ {KS_M}$  as  $\bar{\partial}:\Omega^{0, q}(\wedge^pT^{1,0}M) \longrightarrow \Omega^{0, q+1}(\wedge ^{p}T^{1,0}M)$ by $$\bar{\partial} \alpha=\sum_{I,J,k} \frac{\partial \alpha^{IJ}}{\partial \bar{z}^k} d \bar{z}^k \wedge d \bar{z}^I \otimes \frac{\partial}{\partial z_J}.$$
\end{definition} 
By directly computation we have  \\

\begin{proposition}
For any $\alpha, \beta \in {KS_M}$, we have 
$$\bar{\partial}[\alpha, \beta]=[\bar{\partial} \alpha, \beta]+(-1)^{deg \, {\alpha}}[\alpha, \bar{\partial} \beta]$$
\end{proposition}
\begin{proof}{(sketch)} 
    We can compute this componentwisely,we denote $r:=(i_1,...i_m)$ be an index in $J$, and $s:=(l_1,...,l_m)$ be an index in $L$.  
    $$\bar{\partial}[\alpha, \beta]=\sum_{I,r,K,s}\bar{\partial}[\alpha^{Ir} dz_I \wedge \frac{\partial}{\partial z^{i_1}}\wedge... \wedge \frac{\partial}{\partial z^{i_m}}, \beta^{Ks} dz_K \wedge  \frac{\partial}{\partial z^{l_1}}\wedge... \wedge \frac{\partial}{\partial z^{l_n}}]$$
    
    $$[\bar{\partial} \alpha, \beta]=\sum_{I,r,K,s}[\bar{\partial}\alpha^{Ir} dz_I \wedge \frac{\partial}{\partial z^{i_1}}\wedge... \wedge \frac{\partial}{\partial z^{i_m}}, \beta^{Ks} dz_K \wedge \frac{\partial}{\partial z^{l_1}}\wedge... \wedge \frac{\partial}{\partial z^{l_n}}]$$

    $$[\alpha, \bar{\partial} \beta]=\sum_{I,r,K,s}[\alpha^{Ir} dz_K \wedge \frac{\partial}{\partial z^{i_1}}\wedge... \wedge \frac{\partial}{\partial z^{i_m}}, \bar{\partial}\beta^{Ks} dz_K \wedge \frac{\partial}{\partial z^{l_1}}\wedge... \wedge \frac{\partial}{\partial z^{l_n}}]$$
    expand them  , the point is that the partial differential in $C$ produces $d\bar{z_k}$ and if we move it to the left then we get a coefficient of $(-1)^{deg \, \alpha}$.
    
\end{proof}

Hence we have
\begin{corollary}
    ${KS_M}$ is a DGLA, with extended Schouten-Nijenhui bracket as Lie bracket and $\bp$ as differential. We call ${KS_M}$ the extended Kodaira-Spencer algebra of $M$. We omit ``extended'' when there is no confusion.
\end{corollary}

\begin{remark} \label{r:7.40}
${KS_M}$ carries also a natural exterior algebra  structure. If we  rewrite $\alpha, \beta \in \Omega^{0,q}(\wedge^p T^{1,0}M)$ as  $\alpha= \sum\limits_{I,J} f_{IJ} d\Bar{z}_I \wedge \frac{\partial}{\partial z_J}$,  $\beta=\sum\limits_{I',J'} g_{I'J'} d\Bar{z}_{I'} \wedge \frac{\partial}{\partial z_{J'}}$, we can define wedge product as $$\alpha \wedge \beta = \sum\limits_{I,J,I',J'} (-1)^{I'J} f_{IJ}g_{I'J'}  d\Bar{z}_I\wedge d\Bar{z}_{I'} \wedge \frac{\partial}{\partial z_J}\wedge \frac{\partial}{\partial z_{J'}}.$$ 
\end{remark}

\medskip

Given an $n$ dimensional Calabi--Yau manifold $M$ with nowhere-vanishing holomorphic $n$ form $\omega$, we have an isomorphism between the space of degree $(p-1)$ multi-vector fields on $M$ and the space of $(n-p)$-holomorphic forms on $M$,   $$\tilde{\eta}: \Gamma(\wedge^p T^{1,0}M) \longrightarrow \Omega^{n-p,0}(M),$$
$$v_1\wedge...\wedge v_p\mapsto \iota{(v_1\wedge...\wedge v_p)}(\omega),$$ where  $\iota$ is contraction map, i.e. $(\iota{(v_1\wedge...\wedge v_p)}(\omega))(v_{p+1},...,v_{n})=\omega(v_1,...,v_p,v_{p+1},...,v_{n})$. 

By this, $\tilde{\eta}$ induces an isomorphism 
\begin{align*}
 \eta: \Omega^{0,q}(\wedge^pT^{1,0}M)&\longrightarrow \Omega^{n-p,q}(M) \\
           \alpha \otimes v_1\wedge...\wedge v_p &\mapsto \alpha \otimes \iota{(v_1\wedge...\wedge v_p)}(\omega)
\end{align*}

\begin{lemma} \label{l:9.8}
   The morphism $\eta$ induces an isomorphism of complexes $$(\Omega^{0,*}(\wedge^p T^{1,0}M),\Bar{\partial}) \cong (\Omega^{*,*}(M), \bar{\partial}),$$ i.e. $\bar{\partial}$ commutes with $\eta$.
\end{lemma}
\begin{proof}
    Consider the diagram 
    $$
    \begin{tikzcd}
	\ldots & {\Omega^{0,q-1}(\wedge^{p} T^{1,0}M)} & {\Omega^{0,q}(\wedge^p T^{1,0}M)} & {\Omega^{0,q+1}(\wedge^{p} T^{1,0}M)} & \ldots \\
	\ldots & {\Omega^{n-p,q-1}(M)} & {\Omega^{n-p,q}(M)} & {\Omega^{n-p,q+1}(M)} & \ldots
	\arrow["{\bar{\partial}}", from=1-2, to=1-3]
	\arrow["{\bar{\partial}}", from=1-1, to=1-2]
	\arrow["{\bar{\partial}}", from=1-3, to=1-4]
	\arrow[from=1-4, to=1-5]
	\arrow["{\bar{\partial}}"', from=2-1, to=2-2]
	\arrow["{\bar{\partial}}"', from=2-2, to=2-3]
	\arrow["{\bar{\partial}}"', from=2-3, to=2-4]
	\arrow["\eta"', from=1-2, to=2-2]
	\arrow["\eta"', from=1-3, to=2-3]
	\arrow["\eta"', from=1-4, to=2-4]
	\arrow[from=2-4, to=2-5]
\end{tikzcd}
    $$
    For any $\mu=\sum\limits_{IJ}f_{IJ}d\Bar{z_{IJ}}d\Bar{z_I}\wedge \frac{\partial}{\partial z_J} \in \Omega^{0,q}(\wedge^p T^{1,0}M)$, $$\eta \cdot \Bar{\partial} \mu=\eta \cdot \sum\limits_{IJ}\frac{\partial f_{IJ}}{\partial \Bar{z_k}}d\Bar{z_k} \wedge d\Bar{z_I}\wedge \frac{\partial}{\partial z_J}=\sum\limits_{IJ}\frac{\partial f_{IJ}}{\partial \Bar{z_k}}d\Bar{z_k} \wedge d\Bar{z_I}\wedge \iota(\frac{\partial}{\partial z_J})(\omega).$$ 
    On the other hand, $$\Bar{\partial}\cdot \eta \mu =\Bar{\partial} \cdot \sum\limits_{IJ} f_{IJ}d\Bar{z_I}\wedge \iota(\frac{\partial}{\partial z_J})(\omega)=\sum\limits_{IJ}\frac{\partial f_{IJ}}{\partial \Bar{z_k}}d\Bar{z_k} \wedge d\Bar{z_I}\wedge \iota(\frac{\partial}{\partial z_J})(\omega).$$ Hence the diagram above commutes and $(\Omega^{0,*}(\wedge^* T^{1,0}M),\Bar{\partial}) \cong (\Omega^{*,*}(M), \bar{\partial})$.
\end{proof}

Consider  operations on $\Omega^{0,*}(\wedge^* T^{1,0}M)$ below 
$$
\begin{tikzcd}
	{\Omega^{0,q}(\wedge^p T^{1,0}M)} & {\Omega^{n-p,q}(M)} & {\Omega^{n-p+1,q}(M)} & {\Omega^{0,q}(\wedge^{p-1} T^{1,0}M)}
	\arrow["\eta", from=1-1, to=1-2]
	\arrow["\partial", from=1-2, to=1-3]
	\arrow["{\eta^{-1}}", from=1-3, to=1-4]
\end{tikzcd}
$$
and  define $\Delta:=\eta^{-1}\cdot \partial \cdot \eta$. Notice that $\Delta^2=\eta^{-1} \cdot \partial \cdot \eta \cdot \eta^{-1} \cdot \partial \cdot \eta =\eta^{-1} \cdot \partial^2 \cdot \eta=0$,  we define a new complex $(\Omega^{0,q}(\wedge^{*} T^{1,0}M), \Delta)$ for any $q$ as
$$\begin{tikzcd}
	\ldots & {\Omega^{0,q}(\wedge^{p+1} T^{1,0}M)} & {\Omega^{0,q}(\wedge^p T^{1,0}M)} & {\Omega^{0,q}(\wedge^{p-1} T^{1,0}M)} & \ldots
	\arrow["\Delta", from=1-2, to=1-3]
	\arrow["\Delta", from=1-1, to=1-2]
	\arrow["\Delta", from=1-3, to=1-4]
	\arrow[from=1-4, to=1-5]
\end{tikzcd}$$
and by definition of $\Delta$, we see that the diagram 
$$\begin{tikzcd}
	\ldots & {\Omega^{0,q}(\wedge^{p+1} T^{1,0}M)} & {\Omega^{0,q}(\wedge^p T^{1,0}M)} & {\Omega^{0,q}(\wedge^{p-1} T^{1,0}M)} & \ldots \\
	\ldots & {\Omega^{n-p-1,q}(M)} & {\Omega^{n-p,q}(M)} & {\Omega^{n-p+1,q}(M)} & \ldots
	\arrow["\Delta", from=1-2, to=1-3]
	\arrow["\Delta", from=1-1, to=1-2]
	\arrow["\Delta", from=1-3, to=1-4]
	\arrow[from=1-4, to=1-5]
	\arrow["\partial"', from=2-1, to=2-2]
	\arrow["\partial"', from=2-2, to=2-3]
	\arrow["\partial"', from=2-3, to=2-4]
	\arrow["\eta"', from=1-2, to=2-2]
	\arrow["\eta"', from=1-3, to=2-3]
	\arrow["\eta"', from=1-4, to=2-4]
	\arrow[from=2-4, to=2-5]
\end{tikzcd}
$$
commutes, hence $(\Omega^{0,q}(\wedge^* T^{1,0}M),\Delta) \cong (\Omega^{*,q}(M), \partial)$ as complexes for any $q$.  

\begin{lemma} {(Tian-Todorov)} \label{l:tt}
    Let $\alpha\in \Omega^{0,q}(\wedge^p T^{1,0}M)$, $\beta \in \Omega^{0,q'}(\wedge^{p'} T^{1,0}M)$, then
    $$
    (-1)^p[\alpha,\beta]=(-1)^{deg(\alpha)+1}(\Delta(\alpha \wedge \beta)- \Delta(\alpha) \wedge \beta) -  \alpha \wedge \Delta (\beta).
    $$

    Recall that $deg(\alpha)=q-p+1, deg(\beta)=q'+p'+1$. 
\end{lemma}
\begin{proof}
    See \cite{D.Huybrechts}, chapter 6.
\end{proof}

\begin{lemma} \label{l:9.8}
    For any $p$, the complexes $K^p=(\Omega^{0,*}(\wedge^p T^{1,0}M) \cap ker\Delta, \Bar{\partial})$, $I^p=(\Omega^{0,*}(\wedge^p T^{1,0}M) \cap Im\Delta, \Bar{\partial})$ are well-defined. In particular, $\bp \Delta=-\Delta \bp$  and  $H^i(I^p)=0$ for any $i$.
\end{lemma}
\begin{proof}
    To see that  $K^p$ is a complex, we  need to check that $\Bar{\partial}$ maps $\Omega^{0,q}(\wedge^p T^{1,0}M) \cap ker\Delta$ to $\Omega^{0,p+1}(\wedge^q T^{1,0}M) \cap ker\Delta$. Recall that by Corollary \ref{c:2.25}, $\partial \Bar{\partial}=-\bp \partial$, also $\Bar{\partial}$ commutes with $\eta$ and $\eta^{-1}$. If $\Delta x=0$, then $$\Delta \Bar{\partial} x=\eta^{-1} \partial  \bar{\partial} \eta x=-\eta^{-1} \bar{\partial}  {\partial} \eta x= - \bar{\partial}  \eta^{-1}  {\partial} \eta x= -\Bar{\partial} \Delta x =0,$$ we get $\bp(x)\in \Omega^{0,p+1}(\wedge^q T^{1,0}M) \cap ker\Delta$. 
    
    The proof for the case of $I^p$ is similar. For any $x\in \Omega^{0,*}(\wedge^* T^{1,0}M)$, $\bar{\partial}\Delta x=\bar{\partial}\eta^{-1} \partial \eta x$. By $\partial \bar{\partial}$-lemma \ref{l:pp}, we have $\partial \eta x= \Bar{\partial} y$ for some $y\in \Omega^{*,*}(M)$, hence $$\bar{\partial}\Delta x=\bar{\partial}\eta^{-1}\Bar{\partial} y= \eta^{-1}\bar{\partial} \Bar{\partial} y=0= \Delta \bar{\partial} x.$$

    To see $I^p$ is acyclic, given arbitrary $\Bar{\partial}-$closed element $\Delta x\in \Omega^{0,q}(\wedge^p T^{1,0}M) \cap Im\Delta$,  by $\partial \bar{\partial}$-lemma \ref{l:pp} we have $$\Delta x = \eta^{-1} \partial \eta x = \eta^{-1}\bar{\partial} \eta z= \bp z$$ for some $z\in \Omega^{*,*}(M)$. Therefore it is exact.
\end{proof}

\begin{corollary} \label{c:9.10}
    For any  $p$, $K^p$ is a DGLA, hence a sub-DGLA of ${KS_M}$, that is, the inclusion map $K\longrightarrow {KS_M}$ is an injective morphism between DGLAs.
\end{corollary}
\begin{proof}
    We just need to check that $K^p$ is closed under Lie bracket. By Tian-Todorov lemma \ref{l:tt}, for any $\alpha, \beta \in K^p$, which implies that 
 $\Delta \alpha= \Delta \beta =0$,  we have $$\Delta[\alpha, \beta]=(-1)^{deg(\alpha)+1+p}({\Delta^2}(\alpha \wedge \beta)- \Delta (\Delta(\alpha) \wedge \beta)) -  \Delta(\alpha \wedge \Delta (\beta))=0.$$ Hence $[\alpha, \beta]\in K^p $.
\end{proof}

\begin{corollary} {(Tian–Todorov)} \label{c:BTT}
    Let $M$ be a Calabi-Yau manifold, then $Def_{KS_M}$, the deformation functor of its Kodaira-Spencer algebra is smooth. In particular, let $h=dim(H^{0,1}(T^{1,0}M))$, $KS_M$ has formal moduli space $\C[[t_1,...t_h]]$.
\end{corollary}
\begin{proof}
    For arbitrary $p$, we have a natural inclusion of complexes $$i: K^p=(\Omega^{0,*}(\wedge^p T^{1,0}M) \cap ker\Delta, \Bar{\partial}) \hookrightarrow  (\Omega^{0,*}(\wedge^p T^{1,0}M) , \Bar{\partial}). $$
    Clearly $i_*$, the induced map on cohomology,  maps $H^q(K^p)$ to $H^{0,q}(\wedge^p T^{1,0}M)$ injectively.  $i_*$ is also surjective since for any closed $x\in \Omega^{0,q}(\wedge^p T^{1,0}M)$, $\Bar{\partial}x=0$ implies $$\Delta \Bar{\partial} x=-\Bar{\partial}\Delta x=0,$$ i.e. $\Delta x$ is a closed element in $I^p$, hence by Lemma \ref{l:9.8}, $\Delta x$ is also exact in $I^p$, i.e., we have $$\Delta x= \Bar{\partial}\Delta z= -\Delta \Bar{\partial} z $$ for some $z\in \Omega^{0,q-1}(\wedge^{p}T^{1,0}M)$. Let $y=x+\Bar{\partial}z \in \Omega^{0,q}(\wedge^p T^{1,0}M)$. We have $\Bar{\partial} y= \Delta y=0$ therefore $y \in \Omega^{0,q}(\wedge^p T^{1,0}M) \cap ker\Delta=K^p$ is closed and $x-y=\Bar{\partial z}$, that is, $$i_*(y)=x$$ in $H^{0,q}(\wedge^p T^{1,0}M)$. By Corollary \ref{c:9.10} $K^p$ is a DGLA. Therefore $i_*$ is a  quasi-isomorphism of DGLAs from $K^p$ to $(\Omega^{0,*}(\wedge^p T^{1,0}M), \Bar{\partial})$. 

 Define $$H=\bigoplus\limits_{k=p-q+1}H^{0,q}(\wedge^p T^{1,0}M)$$ 
 as the total cohomology space with respect to differential $\Delta$, equip it with differential $d=0$ and Lie bracket $[\cdot, \cdot]=0$. From this we construct a DGLA $(H^p,d=0, [\cdot,\cdot]=0):=(H^{0,*}(\wedge^pT^{1,0}M),d=0,[\cdot,\cdot]=0)$ for arbitrary $p$. 
 
 For any $x\in \Omega^{0,q}(\wedge^p T^{1,0}M) \cap ker\Delta$, $\eta \Bar{\partial}x=\Bar{\partial} \eta x$ is $\Bar{\partial}$-exact in $\Omega^{*,*}(M)$. By $\partial \Bar{\partial}-$lemma \ref{l:pp}, $\Bar{\partial}\eta x=\partial y$ for some $y\in \Omega^{*,*}(M)$. Let $z=\eta^{-1}y$,  we have $\Bar{\partial}x=\Bar{\partial}\eta \eta^{-1}x=\eta^{-1} \eta  \Bar{\partial}x=\eta^{-1} \partial y = \eta^{-1} \partial \eta z= \Delta z$. Hence any $\Bar{\partial}$-closed element in $K^p$ is $\Delta$-exact and the diagram 
 $$
 \begin{tikzcd}
	{K^p:} & \ldots & {\Omega^{0,q}(\wedge^p T^{1,0}M) \cap ker\Delta} & {\Omega^{0,q+1}(\wedge^p T^{1,0}M) \cap ker\Delta} & \ldots \\
	{H^p:} & \ldots & {H^{0,q}(\wedge^p T^{1,0}M)} & {H^{0,q+1}(\wedge^p T^{1,0}M)} & \ldots
	\arrow["{\bar{\partial}}", from=1-2, to=1-3]
	\arrow["{\bar{\partial}}", from=1-3, to=1-4]
	\arrow["{\bar{\partial}}", from=1-4, to=1-5]
	\arrow["d", from=2-2, to=2-3]
	\arrow["d", from=2-3, to=2-4]
	\arrow["d", from=2-4, to=2-5]
	\arrow["p", from=1-3, to=2-3]
	\arrow["p", from=1-4, to=2-4]
\end{tikzcd}
 $$
where $p$ is the canonical projection from $ker\Delta$ to its cohomology class of $\Delta$, commutes. Hence $p$ is a morphism of complexes from $K^p$ to $H^p$, and it induces an isomorphism on their cohomology groups (or we say that $p$ is a quasi-isomorphism of complexes). By Tian-Todorov lemma, for any $\alpha,\beta \in K^p$, $$[\alpha,\beta]=(-1)^{deg(\alpha)+1+p}(\Delta(\alpha \wedge \beta)- \Delta(\alpha) \wedge \beta) -  \alpha \wedge \Delta (\beta)=(-1)^{deg(\alpha)+1+p}\Delta(\alpha \wedge \beta)$$ which vanishes in its projection in $H$. Therefore, $[p(\alpha),p(\beta)]_K=0=p[\alpha,\beta]$, i.e., $p$ is a morphism of DGLA. Hence $(K^p,\Bar{\partial},[\cdot,\cdot])$ is  caonically quasi-isomorphic to $(H^p,d=0,0)$ as a DGLA.  

Now we consider the DGLAs $$({KS_M}^p, \Bar{\partial}, [\cdot,\cdot]):=(\Omega^{0, *}( \wedge^p T^{1,0}M),\Bar{\partial},[\cdot,\cdot]),$$ $$(K^p,\Bar{\partial},[\cdot,\cdot])= (\bigoplus\limits_{q-p+1=*} \Omega^{0, q}( \wedge^p T^{1,0}M) \cap ker\Delta,\Bar{\partial}, [\cdot,\cdot]),$$  $$(H^p,d=0)=(\bigoplus\limits_{q-p+1=*} H^{0, q}( \wedge^p T^{1,0}M),0,0).$$

By the arguments above,   $(K^p,\Bar{\partial},[\cdot,\cdot])$ is quasi-isomorphic to $({KS_M}^p, \Bar{\partial} , [\cdot,\cdot])$ and  $(K^p,\Bar{\partial},[\cdot,\cdot])$ is quasi-isomorphic to $(H^p,0,0)$. Hence their deformation functors are isomorphic. Since $H^p$ is abelian, by Corollary \ref{c:7.29}, we have that $Def_{\Omega^p_M}$ is pro-represestable by $\C[[t_1,...,t_h]]$ and smooth, where $$h=dim H^{1}(H^p)= \sum\limits_{q=p}dim H^{0,q}(\wedge^pT^{1,0}M)$$. In particular, when $p=1$, $({KS_M}^p, \Bar{\partial},[\cdot,\cdot] )$ is the Kodaira-Spencer algebra of $M$, $Def_{KS_M}$ is pro-represestable by $\C[[t_1,...,t_h]]$, where $h=dim(H^{0,1}(T^{1,0}M))$, therefore it is smooth.

\end{proof}

\newpage

\section{Generalized Moduli Space of Complex Structures on Calabi-Yau Manifold } 

We have seen that a formal deformation of a  complex structure on $(M,J)$ is a formal solution of Maurer-Cartan equation in DGLA $\Omega^{0, *}(T^{1,0}M)$ modulo gauge equivalence and the space of first order deformations is  $H^{0, 1}(T^{1,0}M)$. But we actually only used the DGLA only in degree 0,1,2.  One can enhance the tool we used by introducing a graded version of deformation functor. With such a graded deformation functor, we make full use of the DGLA and can see ``classical deformation problem" as a special case of ``graded deformation problem". 
(We are mainly the following \cite{BK} in this chapter)

\subsection{Graded Deformation Functor, Graded Moduli Space}
The definitions in this section are almost the same as in last chapter but replace algebra with graded algebra. The Artin local ring we used $A_k=\frac{\K \{t\}}{(t^{k+1})}$ can be considered as a graded Artin local ring concentrated in degree $0$.

\begin{definition}
   Let $L$ be a DGLA. Denote  the category of local $\mathbb{Z}$-graded Artin $\K$ algebra by $g\mathbf{Art}_{\mathbb{K}}$.   We call  the functor 
\begin{align*}
        Def^{\mathbb{Z}}_L: \quad \quad \quad &g\mathbf{Art}_{\mathbb{K}} \longrightarrow   \mathbf{Set} \\
                           &A \mapsto \frac{\left\{x \in (L \otimes \mathfrak{m}_A)^1 \mid d x+\frac{1}{2}[x, x]=0\right\}}{exp((L \otimes \mathfrak{m}_A )^0)}        
\end{align*}
the \textbf{graded deformation functor} of $L$. Recall that DGLA $L \otimes \mathfrak{m}_A$ is constructed as in Lemma \ref{e:7.3}, i.e. 
\begin{equation*}
    \begin{split}
        (L \otimes \mathfrak{m}_A)^k &=\bigoplus\limits_{i-j=k}(L^i\otimes \mathfrak{m}^j_A) \\
        d(l\otimes a)&=dl\otimes a \\
        [l_1\otimes a_1,l_2\otimes a_2]&=(-1)^{deg(l_2)deg(a_1)}[l_1,l_2]\otimes a_1a_2
    \end{split}
\end{equation*}
 for any $l\in L, a\in \mathfrak{m}_A$.
\end{definition}

\begin{definition}
   Given DGLA $L$,  we call an element in $Def^{\mathbb{Z}}_L\left(\frac{\K \{t\}}{(t^{k+1})}\right)$  a $k-$order graded deformation. If  every $x_k\in Def^{\mathbb{Z}}_L(\frac{\K \{t\}} {(t^{k+1})})$ can be extended to an element $x_{k+1}\in Def^{\mathbb{Z}}_L(\frac{\K \{t\}} {(t^{k+2})})$, we say $Def^{\mathbb{Z}}_L$ is smooth. 
\end{definition}

Given a $\mathbb{Z}$-graded $\K$ vector space $V=\bigoplus_{i\in \mathbb{Z}} V^i$, the  algebra of formal power series $\K[[t_1,...,t_n]]$  has a natural graded structure if we let the variables $t_1,...,t_n$ be elements in $V$. We call such $\K[[t_1,...,t_n]]$ an algebra of formal power series on graded variables.

\begin{definition}
    We call a graded  algebra pro-Artin if it has form $\frac{\K[[t_1,...,t_n]]}{I}$, where $\K[[t_1,...,t_n]]$ is an algebra of formal power series on graded variables and $I$ is an ideal of it. 
\end{definition}

For any graded (pro-)Artin algebra $A$, we can construct a new (pro-)Artin algebra $A[-1]$ by shifting the graded by one, i.e., $A[-1]^i=A^{i-1}$. 

Given a DGLA $L$,  we define
\begin{definition}
    The graded deformation functor $Def^{\mathbb{Z}}_L$ is called pro-representable if there is a pro-Artin graded algebra $R$ such that $Def^{\mathbb{Z}}_L \cong h_R$.
\end{definition}

\begin{definition}
    If  $Def^{\mathbb{Z}}_L$ is  pro-representable by some pro-Artin graded algebra $R$, then we call $R$ the graded formal moduli space of $Def^{\mathbb{Z}}_L$. 
\end{definition}

From now on, we will call a graded formal moduli space just moduli space unless specified.

\begin{theorem} \label{t:9.6}
    If two DGLAs  $L,L'$  are quasi-isomorphic, then $Def^{\mathbb{Z}}_L \cong Def^{\mathbb{Z}}_{L'}$
\end{theorem}
\begin{proof}
    See \cite{Lurie1}
\end{proof}

\begin{proposition}
    Given a DGLA $L$, if  $Def^{\mathbb{Z}}_L$ is pro-representable, then it is smooth.
\end{proposition}
\begin{proof}
    The proof is essentially same as the proof in Lemma \ref{l:7.23}, just replace Artin algera by graded Artin algebra.
\end{proof}

\begin{proposition} \label{p:9.8}
    If the DGLA $L$ is abelian, then $L$ is pro-representable. 
\end{proposition}
\begin{proof}
    $L$ is abelian implies that the Lie bracket vanishes. Similar as the case in Corollary  \ref{c:7.29}, given a graded Artin algebra $A$, we have
    \begin{align*}
            Def^{\mathbb{Z}}_L(A)&=\frac{ker \, d \cap (L\otimes \mathfrak{m}_A)^1 }{d(L\otimes \mathfrak{m}_A)^0)}=\frac{\bigoplus\limits_{i-j=1}(Z^i(L)\otimes \mathfrak{m}^j_A)}{\bigoplus\limits_{i-j=1}(d(L^{i+1})\otimes \mathfrak{m}^j_A)}\\
            &=\bigoplus\limits_{i-j=1} (H^i(L)\otimes \mathfrak{m}^j_A)
        \end{align*}

Suppose $H^i(L)$ has graded basis $\{ t^i_1,...,t^i_{h_i} \}$ and $\mathfrak{m}^j_A$ is generated by $\{ a^j_1,...,a^j_{g_j} \}$. 

For any element $\phi=\sum_{i-j=1} \sum_{k,l} t^i_k \otimes a^j_l \in Def^{\mathbb{Z}}_L(A)$, we can define a  graded Artin algebra morphism  $\psi=f_A(\phi)\in Hom^{\mathbb{Z}}(\K[[...,t^i_1,...,t^i_{h_i},...,t^{i+1}_1,...,t^{i+1}_{h_{i+1}},...]][1],A)$ by setting $$\psi(t^i_k)=a^j_l.$$ Notice that since $j=i-1$, to get a morphism of graded algebra we have to shift the algebra of formal power series on graded variables by degree $1$. Conversely, given  any graded Artin algebra morphism $\psi\in Hom(\K[[...,t^i_1,...,t^i_{h_i},...,t^{i+1}_1,...,t^{i+1}_{h_{i+1}},...]][1],A)$ with $\psi(t^i_k)=a^j_l$,  the element  $$\phi=\sum\limits_{i-j=1} \sum\limits_{k,l} t^i_k \otimes a^j_l \in Def^{\mathbb{Z}}_L(A)$$ satisfies $f_A(\phi)=\psi$. Hence we have a bijection between $Def^{\mathbb{Z}}_L(A)$ and $Hom^{\mathbb{Z}}(\K[[...,t^i_1,...,t^i_{h_i},...,\\ t^{i+1}_1,...,t^{i+1}_{h_{i+1}},...]][1],A)$. Moreover we have $$Def^{\mathbb{Z}}_L \cong Hom^\mathbb{Z}(\K[[...,t^i_1,...,t^i_{h_i},...,t^{i+1}_1,...,t^{i+1}_{h_{i+1}},...]][1], \cdot ).$$
\end{proof}

In last proposition, for simplicity, let $H=\bigoplus\limits_i H^i(L)$ be the total  cohomology space of $L$ and denote $t_H$ as the graded variables in $H$, we can write the isomorphism above as 
$$ Def^{\mathbb{Z}}_L \cong h_{\K[[t_H]][1]}.$$

In particular,  $ Def^{\mathbb{Z}}_L$ has graded moduli space $\K[[t_H]][+1])$. The Zariski tangent space of $\K[[t_H]][+1]$ is the graded vector space $H[1]=\bigoplus\limits_i H^i(L)[1]$.

\begin{remark}
    Now we see that by extending classical deformation to a graded version, we can use more information about DGLA and the total cohomology space can be explained as the tangent space of the graded moduli space. More about  this ``enhanced" deformation can be found in  \cite{Lurie1} or \cite{Toen}. 
\end{remark}

\medskip

\subsection{Extended Moduli Space of Calabi-Yau Manifold} 

We now go back to the case of Calabi-Yau manifold and define  extended moduli spaces following \cite{BK}. Let $M$ be a $n$ dimensional Calabi-Yau manifold with nowhere-vanishing holomorphic form $\omega$. Consider the graded deformation functor of  ${KS_M}$, which is the extended Kodaira-Spencer algebra of $M$, 
$$
{KS_M}:=\bigoplus\limits_k {KS_M}^k,\quad  {KS_M}^k=:\bigoplus\limits_{q-p+1=k} \Omega^{0, q}(X, \wedge^p T^{0,1}M).
$$

\begin{theorem}{(\cite{BK})} \label{t:9.10}
 $Def^{\mathbb{Z}}_{{KS_M}}$ is pro-representable. 
\end{theorem}
\begin{proof}
    In the proof of Corollary \ref{c:BTT} we introduced the natural inclusion and projection morphism between three DGLAs and showed they are quasi-isomorphisms of DGLAs. In diagram,
$$   \begin{tikzcd}
	{(\Omega^{0,*}(\wedge^p T^{1,0}M) , \bar{\partial},[\cdot,\cdot])} & {(\Omega^{0,*}(\wedge^p T^{1,0}M)\cap ker\Delta , \bar{\partial},[\cdot,\cdot])} & {(H^{0,*}(\wedge^pT^{1,0}M),0,0)}
	\arrow["i"', from=1-2, to=1-1]
	\arrow["p", from=1-2, to=1-3]
\end{tikzcd}$$
  Since $i$ and $p$ are canonical maps induced by inclusion and projection, respectively, after direct sum of all degrees we have again a quasi-isomorphism of the following three DGLAs:
   \begin{equation}\label{quasi}
           \begin{tikzcd}
	{(\bigoplus\limits_{*=q-p+1}\Omega^{0,q}(\wedge^p T^{1,0}M) , \bar{\partial},[\cdot,\cdot])} \\
	{(\bigoplus\limits_{*=q-p+1}\Omega^{0,*}(\wedge^p T^{1,0}M)\cap ker\Delta , \bar{\partial},[\cdot,\cdot])} \\
	{(\bigoplus\limits_{*=q-p+1}H^{0,*}(\wedge^pT^{1,0}M),0,0)}
	\arrow["p"', from=2-1, to=3-1]
	\arrow["i", from=2-1, to=1-1]
\end{tikzcd}
       \end{equation} 

Recall that  $H=\bigoplus\limits_{q,p}H^{0,q}\wedge^p(T^{1,0}M)$ is the total cohomology space of ${KS_M}$ with respect to $\Delta$, which is a graded vector space and an element in  $H^{0,q}(\wedge^pT^{1,0}M)$ has degree $q-p+1$. Denote $t_H$ as the graded variables in $H$. By Lemma \ref{t:9.6} and Proposition \ref{p:9.8},  we have 
$$
Def^{\mathbb{Z}}_{{KS_M}}\cong h_{\C[[t_H]][1]}
$$
\end{proof} 

\begin{remark}  
The ``classical" moduli space of complex structures on $M$ in Lemma  \ref{c:BTT} is $\C[[t_1,...t_h]]$, where $h=dim(H^{0,1}(T^{1,0}M))$, is a subspace of extended moduli space  $\C[[t_H]]$, where $H$ is the whole cohomology space.
\end{remark}

Now we have 
\begin{definition} {(Kontsevich)}  \label{d:dgbv}
    Let $M$ be a Calabi-Yau manifold, and $H=\bigoplus\limits_{q,p}H^{0,q}(\wedge^pT^{1,0}M)$ be the total cohomology space, then  the extended moduli space of complex structure of $M$ is the graded Artin algebra $\C[[t_H]][1]$. In particular, the tangent space of the extended moduli space is the graded vector space $H[1]$
\end{definition}

\begin{example}
We compute an example of extended moduli space.    Let $M$ be a submanifold of $\mathbb{P}^4$ defined by  the equation 
    $$
    \{(z_1,z_2,z_3,z_4,z_5)\in \mathbb{P}^4, z_1^5+z_2^5+z_3^5+z_4^5+z_5^5=0\}.
    $$ 
 Denote  the dimension of $H^{p,q}(M)$ by $h^{p,q}$. The dimensions of cohomology groups of $M$ are (see \cite{M.Gross} for the detail of  computation of the cohomology groups) 
$$
h^{0,0}=h^{3,3}=h^{3,0}=h^{0,3}=h^{1,1}=h^{2,2}=1, \quad h^{1,2}=h^{2,1}=101 
$$

and any other cohomology group is trivial.  By Lemma \ref{l:9.8}, we have $H^{0,q}(\wedge^p T^{1,0}M)\cong H^{3-p,q}(M)$,  then
$$
dim H^{0,0}(\wedge^3 T^{1,0}M)=dimH^0(M)=dimH^{0,3}(M)=dim H^{0,3}(\wedge^3 T^{1,0}M)=1 
$$
and 
$$
dim H ^{0,1}(\wedge^2 T^{1,0}M)=dim H^{0,2}(T^{1,0}M)=101
$$
Then by the definition above, the extended moduli space of $M$ is $\C[[\mathbf{t}_{03},\mathbf{t}_{00},\mathbf{t}_{30},\mathbf{t}_{33},\mathbf{t}_{12},\mathbf{t}_{21}]]$, where $\mathbf{t}_{ij}$ is the graded basis of $H^{0,i}(\wedge^jT^{1,0}M)$, with grading $j-i+2$ (recall the degree $1$ shift in the definition). In particular  the extended moduli space of $M$ has dimension $1+1+1+1+101+101=206$.

\end{example}
\newpage

\section{A Constuction of (Formal) Frobenius Manifold from Extended Moduli Space}

In \cite{BK}, Baranikov an Kontsevich proved that the extended moduli space of a Calabi-Yau manifold has a special structure called formal Frobenius manifold structure. By \cite{Dubrovin}, a formal Frobenius manifold can be seen as a deformation of Frobenius algebra. Since two dimensional topological quantum field theories are characterized by Frobenius algebras (see an introduction to two dimensional topological quantum field theory in Appendix \ref{A}), a formal Frobenius manifold can be seen as a deformation of 2D topological quantum field theory. In this chapter we will follow \cite{BK}, \cite{Manin} to prove that the extended moduli space of a Calabi-Yau manifold is a formal Frobenius manifold.

\subsection{Formal Frobenius Manifold}
A Frobenius manifold is a geometrization of the WDVV equation (see, for instance \cite{Manin} for details). Roughly speaking, it is a manifold equipped with a smooth varying Frobenius algebra structure on its tangent space and  locally everywhere there exists a potential function.

We introduce the following. 
\begin{itemize}
    \item Let $H$ be a  finite-dimensional $\Z_2$-graded vector space over $\C$, generated by a graded basis $\{t_1,...,t_n \}$. 

\item     Let  $<\cdot,\cdot>$ be a non-degenerate symmetric inner product on $H$. 
We  define a $(2,0)-$tensor $g$ such that  $g_{ij}=<t_i,t_j>$. 
\item Denote $H^*$ as the canonical dual space of $H$ under $g$. Precisely, for any $v\in H$, its dual vector is defined as $v^*=<v,\cdot>$.

\item  Let  $K=\C[[t_H]]$ be the algebra of  formal power series with variables in $H^*$.

\smallskip

\item Let $*$ be an associative and  commutative  multiplication on $K \otimes H$. By this multiplication we define a $(1,2)-$tensor $(A_{ij}^k)$ with value in $K$ on $H$ such that $t_i * t_j = \sum_k A_{ij}^k t_k$. Hence in term of $(A_{ij}^k)$, where $A_{ij}^k \in K$, commutativity and associativity are expressed as follows (the degree of any  $t_i$ is denoted by $\Tilde{i}, \text{where} \,  i=1,...,n$)

\begin{itemize}
    \item For any $a,b,c$, $A_{ba}^c=(-1)^{\Tilde{a}\Tilde{b}}A_{ab}^c$ \hfill(commutativity)
    \item For any $a,b,c,d$, $\sum\limits_e A_{ab}^e A_{ec}^d=(-1)^{\Tilde{a}(\Tilde{b}+\Tilde{c})}\sum\limits_e A_{bc}^e A_{ea}^d$ \hfill(associativity)

    \end{itemize}

\item Define $(3,0)-$tensor $(A_{ijk})$ by $A_{ijk}=\sum_l A_{ij}^l g_{lk}$

\end{itemize}

\begin{definition}(Manin \cite{Manin})
  We call $(H,g, *)$ a formal Frobenius manifold if the objects $H,g,\wedge$ are defined above and  there is a potential function $\Phi \in K$ such that $A_{ijk}=\frac{\partial^3 \Phi}{\partial t_i \partial t_j \partial t_k}$. In this case, we say that the tensor $(A_{ij}^k)$ is potential.
\end{definition}

\begin{definition}
   We call a finite dimensional commutative, associative $\C$ algebra $A$ with a non-degenerate  bilinear map $<\cdot,\cdot>: A \times A \longrightarrow A$ such that $<ab,c>=<a,bc>$ a \textbf{Frobenius algebra}. 
\end{definition}

    Given a formal Frobenius manifold $(H, *, g)$. For arbitrary point $x=\left(x_1, \ldots, x_n\right) \in$ $H$, that is, $x=\sum_i x_i t_i$. If we see $H$ as the tangent space of $H$ at $x$, we have a multiplication on tangent space at $x$, denoted also by $*$. Explicitly, the multiplication on $H$ viewed as the tangent space at $x$ is given by $t_i * t_j=\sum_k A_{i j}^k(x) t_k$. In this way, we construct a multiplication on the tangent space at any point on $H$. Since $\left(A_{i j}^k\right)$ is smooth, we have a smooth family of multiplications on tangent spaces of $H$ (or on $H$ since we identify $H$ with its tangent space at any point), which varies smoothly with respect to the base point of the tangent space.

It is easy to check that $T_x H$ is always a Frobenius algebra with algebra multiplication $*$ and inner product from the inner product on for the formal Frobenius manifold. Since $x$ is an arbitrary point on $H$, in this way we have constructed a smooth family of Frobenius algebras $(T_x H $ over $H$, where $x \in H$.

\begin{remark}
In \cite{Manin}, a formal Frobenius manifold is defined as the formal completion of $H$ at the origin. In this definition, a formal Frobenius manifold can be seen as a formal manifold which contains only one point or a germ of Frobenius manifold at a point. 

\end{remark}

\medskip

\subsection{Construction of Formal Frobenius Manifold from Extended Moduli Space }{
Now we give the sketch of the construction of formal Frobenius manifold in \cite{BK},  with a more general method introduced by Y.Manin in  \cite{Manin}. \\

 Let $M$ be a Calabi-Yau manifold with holomorphic top form $\omega$, $H=\bigoplus\limits_{q,p}H^{0,q}(\wedge^pT^{1,0}M)$ with basis $\{ t_i\}$ be the total cohomology space of $M$.  

Recall that in Remark \ref{r:7.40}, ${KS_M}$ carries  a natural exterior algebra  structure. Let $\alpha= \sum\limits_{I,J} f_{IJ} d\Bar{z}_I  \wedge \frac{\partial}{\partial z_J}$,  $\beta=\sum\limits_{I',J'} g_{I'J'} d\Bar{z}_{I'} \wedge \frac{\partial}{\partial z_{J'}} \in {KS_M}$, we can define wedge product as $$\alpha \wedge \beta = \sum\limits_{I,J,I',J'} (-1)^{I'J} f_{IJ}g_{I'J'}  d\Bar{z}_I\wedge d\Bar{z}_{I'} \wedge \frac{\partial}{\partial z_J}\wedge \frac{\partial}{\partial z_{J'}}$$

 \begin{lemma} \label{l:10.11}
  Define a multiplication on $H$ as a wedge product induced from the wedge product on ${KS_M}$ and an inner product $<\cdot,\cdot>$ as $<t_i,t_j>=\int_M \iota{(t_i\wedge t_j)}(\omega) \wedge \omega$, then $(H,\wedge, <\cdot,\cdot>)$ is a Frobenius algebra.
 \end{lemma}
\begin{proof}
    Associativity and $<t_i\wedge t_j, t_k>=<t_i,t_j\wedge t_k>$ are obvious, we prove  now  that  $<\cdot,\cdot>$ is non-degenerate. 

    If $<x,y>=0$ for any $y$, we denote a representative of cohomology class of $x$ by $\hat{x}\in \Omega^{0,q}(\wedge^p T^{1,0}M)$, similar for $y,z,...$ \, .  
 Suppose that locally  $\hat{x}=fd\Bar{z}^Q\otimes \frac{\partial}{\partial z^P}$.   Since $M$ is oriented, we can define $\frac{\partial}{\partial z^{N-P}}$ as the completment of $\frac{\partial}{\partial z^{P}}$ in $\frac{\partial}{\partial z^1}\wedge...\wedge \frac{\partial}{\partial z^n}$, i.e. $$ 
\frac{\partial}{\partial z^{N}}\wedge \frac{\partial}{\partial z^{N-P}}=\frac{\partial}{\partial z^1}\wedge...\wedge \frac{\partial}{\partial z^n}.$$ 

Similarly we can define   $d\Bar{z}^{N-Q}$ such that $$d\Bar{z}^{N}\wedge d\Bar{z}^{N-Q}=\Bar{\omega}.$$ 

Now let $y=\bar{f} d\Bar{z}^{N-Q}\otimes \frac{\partial}{\partial z^{N-P}}$, then  $<x,y>=\int_M |f|^2 \Bar{\omega}\wedge \omega$, since $\omega$ is nowhere-vanishing on $M$, $\Bar{\omega}\wedge \omega$ is a volume form on $M$. Therefore $<x,y>=0$ only if $f=0$, i.e. $x=0$. Hence the inner product is non-degenerate. 
\end{proof}

 To see the formal Frobenius manifold structure on $H$, we  construct  a family of Frobenius algebras over  $H$.

\begin{lemma} \label{l:10.11}
    Let $(L,d,[\cdot,\cdot])$ be a DGLA,   $a\in L^1$ be a solution of Maurer-Cartan equation, then the deformed differential $d_a$, defined by $d_a(b)=db+[a,b]$ is indeed a differential, i.e. $d_a^2=0$. Moreover, $(L,d_a,[\cdot,\cdot])$ is a DGLA.
\end{lemma}
\begin{proof}
    By definition, 

    \begin{equation}
          \begin{split}
        d_a^2(b)&=d(db+[a,b])+[a,db+[a,b]] \\
        &=d[a,b]+[a,db]+[a,[a,b]] \\
        &=[da,b]-[a,db]+[a,db]+[a,[a,b]] \\
        &=[da,b]+[a,[a,b]].
    \end{split}
    \end{equation}

    By graded Jacobi identity $[a,[a,b]]=[[a,a],b]-[a,[a,b]]$, we have $$[a,[a,b]]=\frac{1}{2}[[a,a],b]$$. Since $a\in MC(L)$,  we have $$d_a^2(b)= [da+\frac{1}{2}[a,a],b]=0.$$  

    To prove  $(L,d_a,[\cdot,\cdot])$ is DGLA, what is left is to check that $d_a[X,Y]=[d_aX,Y]+(-1)^{deg(X)}[X,d_aY]$. By definition, 
    \begin{equation}
    \begin{split}
      d_a[X,Y]&=d[X,Y]+[a,[X,X]] \\
      &=[dX,Y]+(-1)^{deg(X)}[X,dY]+[[a,X],Y]+(-1)^{deg(X)}[X,[a,Y]]\\
      &=[dX+[a,X],Y]+[X,(-1)^X(dY+ [a,Y])]\\
      &=[d_aX,Y]+(-1)^{deg(X)}[X,d_aY]. 
      \end{split} 
  \end{equation}

\end{proof}

\smallskip

Given a DGLA $L$, suppose $Def_L^\Z$ is pro-representable. Let $K:=\C[[t]]$ be the (graded) algebra of formal power series in $H$. We have a formal solution $\gamma(t)\in {KS_M} \otimes \C[[t]]$ of Maurer-Cartan equation  
$$
\bar{\partial}\gamma(t)+\frac{1}{2}[\gamma(t),\gamma(t)]=0 
$$
and $\gamma(t)=\sum\limits_a \gamma_{a}t_a+\sum\limits_{a_1,a_2} \gamma_{a_1 a_2}t_{a_1}t_{a_2}+...$ . We construct a DGLA $({KS_M} \otimes K, \Bar{\partial},[,])$ as in Lemma \ref{e:7.3}.  Now by  Lemma \ref{l:10.11} we can construct a family of new DGLAs $({KS_M} \otimes K, \Bar{\partial}_{\gamma(t)},[,])$ with deformed differential $\Bar{\partial}_{\gamma(t)}= \Bar{\partial} + [\gamma(t), \cdot]$. The point is, the deformed differential  would produce isomorphic cohomology group as $H$. Hence it will give us a family of multiplications on $H$(hence a family of Frobenius algebras). 

\begin{theorem}(Manin \cite{Manin})
   There is a  quasi-isomoprhism of DGLAs from $({KS_M} \cap ker\Delta \otimes K , \bar{\partial}_{\gamma_t},[\cdot,\cdot])$ to $(H\otimes K,d=0,[\cdot,\cdot]_H=0) $. In particular, the cohomology of complex $({KS_M} \otimes K , \bar{\partial}_{\gamma(t)})$ is canonically isomorphic to $(H,d=0)$.
\end{theorem}
\begin{proof}
Notice that $K=\C[[t]]$ is flat, the cohomology of  $({KS_M} \otimes K,\bp)$ is just $H\otimes K$. Recall that in Theorem \ref{t:9.10},  cohomology of $(\Omega\otimes K)$ is canonically isomorphic to the cohomology of $({KS_M} \cap ker\Delta \otimes K , \bar{\partial})$. Hence to prove the isomorphism of cohomology we can replace $({KS_M} \otimes K,\bp)$ with $({KS_M} \cap ker\Delta \otimes K , \bar{\partial})$. In Theorem \ref{t:9.10}, we have a canonical quasi-isomorphism $$p: ({KS_M} \cap ker\Delta  , \bar{\partial},[\cdot,\cdot]) \longrightarrow (H,0,0).$$ Again since $K$ is flat, we have a   canonical quasi-isomorphism $$p\otimes id: ({KS_M}\cap ker\Delta\otimes K,\bp) \longrightarrow (H\otimes K,0,0).$$ Without causing confusion we will write also $p$ for $p\otimes id$.  

We want to prove that $p$ is a quasi-isomorphism from $ ({KS_M}\cap ker\Delta\otimes K,\bp_{\gamma(t)},[\cdot,\cdot]) $ to $H\otimes K,d=0,0)$. If it is  true then we have a canonical isomorphism of the  cohomomolgy of $({KS_M} \otimes K , \bar{\partial}_{\gamma_t})$ and $(H,d=0)$. 

For any $\phi \in {KS_M}\cap ker\Delta$,  $$p \bp_{\gamma(t)}\phi=p(\bp \phi+[\gamma(t),\phi])=[p\phi,p\phi]=0=d\bp_{\gamma(t)}.$$ Hence $p$ is a morphism of DGLA (compability with Lie bracket can be checked in the same way as in undeformed differential case).    Since $p$ is projection from ${KS_M}\cap ker\Delta$ to $H$, $p:\Omega\cap ker\Delta\otimes K \longrightarrow H\otimes K$ is surjective. Hence $p$ induces surjective map between cohomology groups of DGLAs. What is left to prove is the injectivity of $p$ on cohomology, that is, we want to prove $ker \bp_{\gamma(t)} \cap ker \, p\subset Im \bp_{\gamma(t)}$. 

To prove this, let $c=\sum\limits_{i\geq 0}c_i \in ker \bp_{\gamma(t)} \cap ker \, p \subset {KS_M}\otimes K$ where $c_i$ represents degree $i$ (in $K$ part) homogeneous part of $c$. Explicitly, for example, we can write $c_i=\gamma_{a_1...a_i} t_{a_1}...t_{a_i}$. Similarly we write $\gamma= \sum \limits_{i\geq 0} \gamma_i, \gamma_0=0$  also in this way and omit $t$ for simplicity. $\bp_{\gamma }(c)=0$ implies  $c_0=0$ and 
 \begin{equation} \label{eq:14}
       \bp c_n=-\sum_{i+j=n}[\gamma_i,c_j].  
\end{equation}

Since $p$ is also a quasi-isomorphism from $({KS_M}\cap ker\Delta,\bp)$ to $(H,0)$, $dc_0=0$ implies that $p\circ dc_0=0$ is exact in $(H,0)$. Hence $c_0$ is exact in $({KS_M}\cap ker\Delta \otimes K, \bp)$.  To prove $c\in Im \bp_{\gamma (t)}$ is equivalent to find $a=\sum \limits_{i\geq 0}a_i \in ker \bp_{\gamma (t)}$ such that 

\begin{equation} \label{eq:15}
      c_{n}=\bp a_{n}+ \sum_{i+j=n}[\gamma_i,a_j]
 \end{equation}
  which is satisfied when $n=0$. We find  $a$ by induction. Assume we have already $a_0,...a_n$ satisfying Eq. \ref{eq:15}. We want to find $a_{n+1}$ such that $\bp a_{n+1}=c_{n}- \sum \limits_{i+j=n}[\gamma_i,a_j]$. Hence it is  enough to check if $c_{n+1}- \sum \limits_{i+j=n+1}[\gamma_i,a_j]$ is exact in $(\Omega\cap ker\Delta\otimes K,\bp)$. 
  \smallskip

We have $$\bp (c_{n}- \sum \limits_{i+j=n}[\gamma_i,a_j])=\bp c_{n+1}-\sum \limits_{i+j=n+1}[\bp \gamma_i,a_j]+\sum \limits_{i+j=n+1}[\gamma_i,\bp a_j]$$, substitute Eq. \ref{eq:14} into it to replace $\bp \gamma_i$,  substitute equation \ref{eq:15} to replace $\bp a_j$ and using \\ $\bp \gamma_i= -\frac{1}{2}\sum \limits_{i+j=n}[\gamma_i,\gamma_j] $ to replace $\bp \gamma_i$,  we have 
  \begin{equation}
        \bp (c_{n}- \sum_{i+j=n}[\gamma_i,a_j])=  -\sum_{i+j+k=n+1}[\gamma_i[\gamma_j,a_k]]+\frac{1}{2}\sum_{i+j+k=n+1}[[\gamma_i,\gamma_j],a_k]
  \end{equation}
which is zero due to graded Jacobi identity. Hence $(c_{n}- \sum\limits_{i+j=n}[\gamma_i,a_j])$ is closed in $(\Omega\cap ker\Delta\otimes K,\bp)$.  Since $c\in ker \, p $ and $p$ is a quasi-isomorphism from $({KS_M}\cap ker\Delta\otimes K,\bp,[\cdot,\cdot])$ to $(H\otimes K,0,0)$, we see that $c$ is exact. 
\end{proof}

Denote $H_{\gamma(t)}\otimes K$ as the cohomology of  ${KS_M} \cap ker \Delta \otimes K$ with respect to $\bp_{\gamma(t)}$.  By last theorem we have $H_{\gamma(t)}\otimes K\cong H\otimes K$.  Since we have a natural wedge product on ${KS_M} \cap ker \Delta \otimes K$, the wedge product induces a product on $H_{\gamma(t)}\otimes K$, which is the family of multiplication of Frobenius algebra we are looking for.

\begin{definition}
    Multiplication $\wedge_t$ on $H\otimes K$ is defined as the following. Let $x,y\in H\otimes K$ and $\hat{x},\hat{y}$ be their representatives in ${KS_M} \otimes K$,
    $$
    x\wedge_t y= \hat{x}\wedge \hat{y}  \, mod \, \bp_{\gamma(t)}
    $$

Given such a multiplication $\wedge_t$, we define a $(2,1)-$tensor with value in $H^*$, written  as $(A_{ij}^k)$ where $A_{ij}^k \in H^*$ by $t_i \wedge_t t_j =\sum \limits_k A_{ij}^k t_k $.
\end{definition}

    With the definition above,  we have constructed a family of  Frobenius algebras over $H$ and when $t=0$,  $\bp_{\gamma(t)}$ is just $\bp$, i.e. we have the ordinary wedge product at the origin of $H$. Hence  we call $\wedge_t$ \textbf{deformed wedge product}. 

\medskip

    To prove that $(H,\wedge_t, <\cdot,\cdot>)$ is a formal Frobenius manifold, what is left is only to prove that $(A_{ijk})$ is potential. The proof is quite technical so we omit it here, interested reader can found the proof in \cite{Manin}.

\subsection{Conclusion}
In last section we have seen that the extended moduli space of complex structures of a Calabi-Yau manifold admits a structure of formal Frobenius manifold.  In Gromov-Witten theory, a deformation of the cohomology ring of a Calabi-Yau manifold (called quantum cohomology) also admits a formal Frobenius manifold structure on it, where the deformed product is defined by Gromov-Witten invariant, see detail in \cite{Salamon}.  It is conjectured in \cite{BK} that this formal Frobenius manifold from quantum cohomology is related to the formal Frobenius manifold structure on the extended moduli space.

\newpage

\begin{appendices}

\section{$2D$ Topological Quantum Field Theory and Frobenius Algebra } \label{A}
For simplicity we only give an unrigorous definition of topological quantum field theory following \cite{Lurie2}.
\begin{definition}{(unprecise)}
 Define $\mathbf{Cob}(d)$ as a category whose objects are closed $(d-1)-$dimensional closed manifolds    and  morphisms between two objects $M,N$ are defined as $Hom(M,N):=\{B \text{ is a } d- \text{dimensional manifold },  \partial B =\Bar{M}\sqcup N \} /\sim $, where the equivalence relation is diffeomorphism relative to boundary. $\text{id}_M:=M\times [0,1]$. The composition of two morphisms is given by gluing two bordisms. We call $\mathbf{Cob}(d)$ the bordism category of dimension $d$.
 
 A \textbf{topological quantum field theory} (\textbf{TQFT}) is a functor $Z: \mathbf{Cob}(d) \longrightarrow \mathbf{\C Vek}$, from the category of bordism to the category of (finite-dimensional) $\C$ vector spaces, which satisfies,  
$$
Z(M\sqcup N) \cong Z(M) \otimes Z(N), \quad Z(\emptyset)\cong \C , \quad Z(\Bar{M})=Z(M)^*
$$
 
\end{definition}

\begin{example}
    When $d=2$, a $(d-1)$ dimensional closed manifold can only be disjoint union of $S^1$. Let $A=S^1\sqcup S^1, B=S^1, C=S^1\sqcup S^1$, naturally we have a  morphism $\phi$ from $A$ to $B$, a morphism  $\psi$ from $B$ to $C$ and their composition which are described by the diagram below, see figure \ref{fig:gluing},

Given a TQFT, i.e. a functor $Z: \mathbf{Cob}(d) \longrightarrow \mathbf{\C Vek}$, we denote $V=Z(S^1)\in \mathbf{\C Vek} $, then $Z(A)=V\otimes V$, $Z(B)=V$, $Z(C)=A\otimes A$. $\hat{\phi}:=Z(\phi)$ is a linear map from $V\otimes V$ to $V$. Hence the bilinear map $\hat{\phi}:V \times V \longrightarrow V$ gives $V$ an $\C$-algebra structure. 
\end{example}

\begin{lemma}
With the notation in the last example, $V$ as a $\C$ algebra is commutative and associative.
\end{lemma}
\begin{proof} (sketch)
    Let $a,b\in V$,  we can flip the pant without changing $\hat{\phi}$, then for any $a,b \in V$, $\hat{\phi}(b,a)=\hat{\phi}(a,b)$, see figure \ref{fig:commutativity}.
  
\end{proof}

    Similar for associativity, changing of the order of multiplication is achieved again by flipping  vertically: 
   $\hat{\phi}(\hat{\phi}(a,b),c)=\hat{\phi}(a,\hat{\phi}(b,c))$, see figure \ref{fig:Associativity}.

If we interpret unit disk $D=\{ x\in \C, |x|\le 1 \}$ as a bordism from empty set  to $S^1$, then $Z(D)$ defines a morphism from $\C$ to $V$ and the image of $1\in \C$ under $Z(D)$ gives a unit in $V$ as $\C$ algebra. Conversely if we interpret $D$ as a bordism from  $S^1$ to $\emptyset$, we have a linear map $tr: V\longrightarrow \C$, see figure \ref{fig:disk}.  

With the operations above, it is not hard to prove that
\begin{lemma}
    $tr\circ\hat{\phi}$, written as $\textlangle \cdot , \cdot \textrangle$, is a non-degenerate map from $V\times V$ to $\C$. Moreover for any $a,b,c \in V, <\hat{\phi}(a,b),c>=<a,\hat{\phi}(b,c)>$.  
\end{lemma}
\begin{proof}
    See \cite{Lurie2}
\end{proof}

Now we see that a  $2D$ TQFT gives us an algebra which is commutative, associative with a non-degenerate bilinear map $\textlangle \cdot , \cdot \textrangle$ which satisfies $\textlangle ab , c \textrangle=  \textlangle a , bc \textrangle$. 

\begin{definition} \label{A:5}
    We call a commutative, associative algebra $F$ with a non-degenerate bilinear map satisfying $\textlangle ab , c \textrangle=  \textlangle a , bc \textrangle$ for any $a,b,c$ in the algebra \textbf{Frobenius algebra}, denoted by $(F,*, <\cdot,\cdot>)$ where $*$ is the algebra multiplication.
\end{definition}

Hence a $2D$ TQFT determines a Frobenius algebra. The converse is also true. 

\begin{theorem}
     Given a Frobenius algebra,  one can uniquely determine a 2D TQFT.
\end{theorem}
\begin{proof}
    See \cite{Lurie2}.
\end{proof}
Hence two dimensional TQFT can be  completely characterized by Frobenius algebra.

\begin{remark}
    One can replace the target $\mathbf{\C Vek}$ of TQFT with the category of finite dimensional $\Z_2$ graded vector space (or super-vector space) $\mathbf{\C Vek}^{\mathbb{Z}_2}$ and require  $\hat{\phi}(a,b)=(-1)^{deg(a)deg(b)}\hat{\phi}(b,a)$. This super version of TQFT  is denoted as  super TQFT. Define the super-Frobenius algebra as a $\Z_2$ graded Frobenius algebra such that it satisfies all requirement of a Frobenius algebra. The commutativity in $\Z_2$ graded-case means super-commutativity, i.e., $ab=(-1)^{deg(a)deg(b)}ba$.  Similarly, we have a one to one correspondence between  super TQFT  and  super-Frobenius algebra. We omit the  ``super'' if there is no confusion.
\end{remark}

  \begin{figure}[h]
    \centering
    \includegraphics[width=13cm]{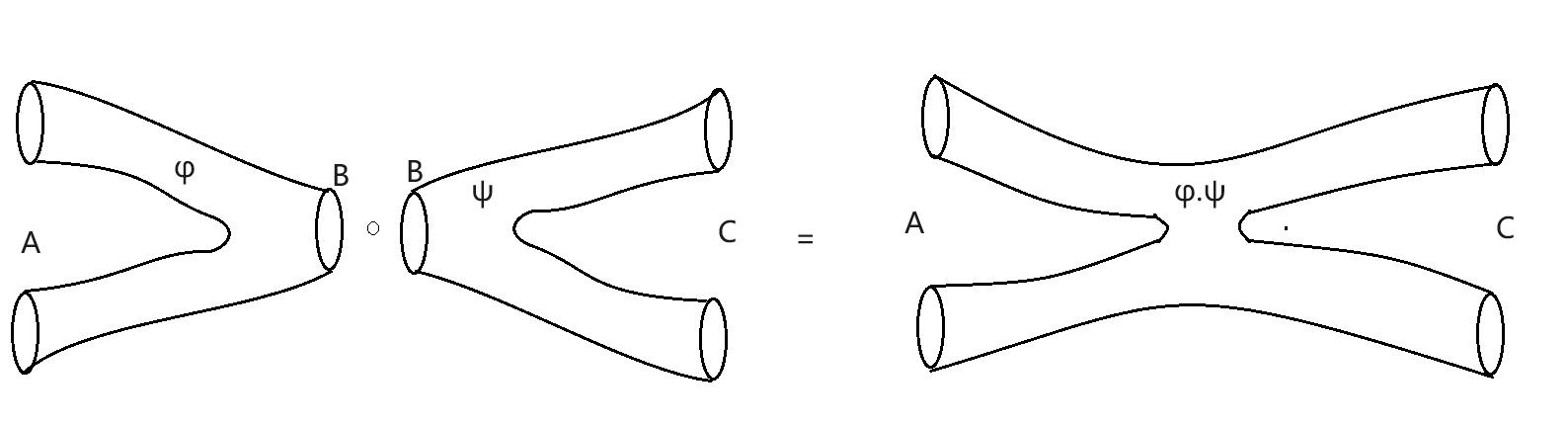}
    \caption{Gluing}
    \label{fig:gluing}
\end{figure}

\tikzset{every picture/.style={line width=0.75pt}} 
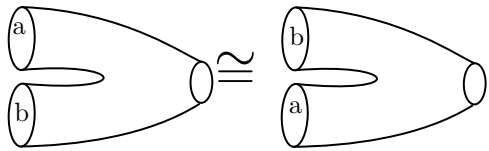
\begin{figure}[h]
    \centering
   \begin{tikzpicture}  
[x=0.75pt,y=0.75pt,yscale=-0.75,xscale=0.7]
\centering

\draw   (16.13,56.6) .. controls (21.29,56.69) and (25.3,66.22) .. (25.09,77.91) .. controls (24.87,89.59) and (20.51,98.99) .. (15.35,98.91) .. controls (10.18,98.83) and (6.17,89.29) .. (6.39,77.61) .. controls (6.61,65.93) and (10.97,56.52) .. (16.13,56.6) -- cycle ;
\draw   (16.42,108.12) .. controls (21.58,108.26) and (25.46,117.84) .. (25.08,129.52) .. controls (24.71,141.2) and (20.22,150.55) .. (15.06,150.41) .. controls (9.9,150.27) and (6.02,140.69) .. (6.39,129.01) .. controls (6.77,117.33) and (11.26,107.98) .. (16.42,108.12) -- cycle ;
\draw    (15.35,98.91) .. controls (89.67,92.54) and (98.96,115.22) .. (16.42,108.12) ;
\draw    (16.13,56.6) .. controls (91,59.05) and (136.09,90.94) .. (145.37,93.39) ;
\draw    (15.06,150.41) .. controls (54.53,149.81) and (106.91,143.68) .. (144.05,121.6) ;
\draw   (145.37,93.39) .. controls (149.71,93.43) and (153.16,99.65) .. (153.09,107.28) .. controls (153.02,114.91) and (149.44,121.07) .. (145.11,121.03) .. controls (140.77,120.99) and (137.31,114.78) .. (137.39,107.15) .. controls (137.46,99.52) and (141.04,93.36) .. (145.37,93.39) -- cycle ;
\draw   (213.43,56.8) .. controls (218.59,56.91) and (222.53,66.48) .. (222.22,78.16) .. controls (221.92,89.84) and (217.49,99.21) .. (212.32,99.1) .. controls (207.16,98.98) and (203.22,89.42) .. (203.53,77.74) .. controls (203.83,66.06) and (208.27,56.68) .. (213.43,56.8) -- cycle ;
\draw   (213.32,108.31) .. controls (218.48,108.49) and (222.29,118.1) .. (221.82,129.77) .. controls (221.36,141.45) and (216.8,150.77) .. (211.64,150.6) .. controls (206.48,150.42) and (202.67,140.81) .. (203.13,129.14) .. controls (203.6,117.46) and (208.16,108.14) .. (213.32,108.31) -- cycle ;
\draw    (212.32,99.1) .. controls (286.69,93.22) and (295.8,115.96) .. (213.32,108.31) ;
\draw    (213.43,56.8) .. controls (288.28,59.74) and (333.12,91.93) .. (342.39,94.44) ;
\draw    (211.64,150.6) .. controls (251.11,150.25) and (303.54,144.47) .. (340.84,122.64) ;
\draw   (342.38,94.44) .. controls (346.72,94.51) and (350.13,100.74) .. (350,108.38) .. controls (349.87,116.01) and (346.24,122.14) .. (341.91,122.08) .. controls (337.57,122.01) and (334.16,115.77) .. (334.29,108.14) .. controls (334.43,100.51) and (338.05,94.38) .. (342.38,94.44) -- cycle ;

\draw (152.2,84.55) node [anchor=north west][inner sep=0.75pt]  [font=\huge] [align=left] {$\displaystyle \cong $};
\draw (7.08,65.28) node [anchor=north west][inner sep=0.75pt]   [align=left] {a};
\draw (8.41,118.63) node [anchor=north west][inner sep=0.75pt]   [align=left] {b};
\draw (206.67,67.73) node [anchor=north west][inner sep=0.75pt]   [align=left] {b};
\draw (206.01,118.01) node [anchor=north west][inner sep=0.75pt]   [align=left] {a};

\end{tikzpicture}
    \caption{Commutativity}
    \label{fig:commutativity}
\end{figure}

\begin{figure}[b]
    \centering

\tikzset{every picture/.style={line width=0.75pt}} 

\begin{tikzpicture}[x=0.75pt,y=0.75pt,yscale=-0.7,xscale=0.7]

\draw   (21.13,18.6) .. controls (26.29,18.69) and (30.3,28.22) .. (30.09,39.91) .. controls (29.87,51.59) and (25.51,60.99) .. (20.35,60.91) .. controls (15.18,60.83) and (11.17,51.29) .. (11.39,39.61) .. controls (11.61,27.93) and (15.97,18.52) .. (21.13,18.6) -- cycle ;
\draw   (21.42,70.12) .. controls (26.58,70.26) and (30.46,79.84) .. (30.08,91.52) .. controls (29.71,103.2) and (25.22,112.55) .. (20.06,112.41) .. controls (14.9,112.27) and (11.02,102.69) .. (11.39,91.01) .. controls (11.77,79.33) and (16.26,69.98) .. (21.42,70.12) -- cycle ;
\draw    (20.35,60.91) .. controls (94.67,54.54) and (103.96,77.22) .. (21.42,70.12) ;
\draw    (21.13,18.6) .. controls (96,21.05) and (141.09,52.94) .. (150.37,55.39) ;
\draw    (21.12,111.84) .. controls (60.59,111.24) and (112.97,105.11) .. (150.11,83.03) ;
\draw   (150.37,55.39) .. controls (154.71,55.43) and (158.16,61.65) .. (158.09,69.28) .. controls (158.02,76.91) and (154.44,83.07) .. (150.11,83.03) .. controls (145.77,82.99) and (142.31,76.78) .. (142.39,69.15) .. controls (142.46,61.52) and (146.04,55.36) .. (150.37,55.39) -- cycle ;
\draw   (19.96,133.49) .. controls (24.57,133.31) and (28.68,142.66) .. (29.15,154.38) .. controls (29.62,166.1) and (26.26,175.75) .. (21.66,175.93) .. controls (17.05,176.12) and (12.93,166.77) .. (12.47,155.05) .. controls (12,143.33) and (15.35,133.68) .. (19.96,133.49) -- cycle ;
\draw    (19.96,133.49) .. controls (82,112.6) and (249,146.6) .. (240,107.6) ;
\draw    (150.11,83.03) .. controls (202,80.6) and (226,90.6) .. (240,107.6) ;
\draw    (150.37,55.39) .. controls (215,46.6) and (286,94.6) .. (299,102.6) ;
\draw    (21.66,175.93) .. controls (95,163.6) and (264,152.6) .. (299,131.6) ;
\draw   (299.32,102.6) .. controls (301.48,102.63) and (303.17,109.14) .. (303.08,117.14) .. controls (302.99,125.15) and (301.16,131.62) .. (299,131.6) .. controls (296.84,131.58) and (295.15,125.07) .. (295.24,117.06) .. controls (295.33,109.05) and (297.15,102.58) .. (299.32,102.6) -- cycle ;
\draw   (380.09,180.59) .. controls (385.24,180.18) and (388.64,170.41) .. (387.68,158.76) .. controls (386.73,147.12) and (381.78,138.01) .. (376.64,138.42) .. controls (371.49,138.82) and (368.09,148.59) .. (369.04,160.24) .. controls (370,171.89) and (374.94,181) .. (380.09,180.59) -- cycle ;
\draw   (377.12,129.16) .. controls (382.26,128.69) and (385.53,118.88) .. (384.42,107.25) .. controls (383.31,95.62) and (378.24,86.57) .. (373.1,87.03) .. controls (367.96,87.5) and (364.69,97.31) .. (365.8,108.94) .. controls (366.91,120.57) and (371.98,129.62) .. (377.12,129.16) -- cycle ;
\draw    (376.64,138.42) .. controls (451.22,140.08) and (459.05,116.85) .. (377.12,129.16) ;
\draw    (380.09,180.59) .. controls (454.65,173.42) and (497.64,138.75) .. (506.75,135.71) ;
\draw    (374.19,87.54) .. controls (413.62,85.65) and (466.29,88.46) .. (504.74,108.15) ;
\draw   (506.75,135.71) .. controls (511.07,135.4) and (514.13,128.98) .. (513.58,121.37) .. controls (513.02,113.76) and (509.07,107.84) .. (504.74,108.15) .. controls (500.41,108.46) and (497.36,114.88) .. (497.91,122.49) .. controls (498.47,130.1) and (502.42,136.02) .. (506.75,135.71) -- cycle ;
\draw   (371.67,66) .. controls (376.28,65.89) and (379.79,56.3) .. (379.52,44.58) .. controls (379.25,32.85) and (375.29,23.43) .. (370.68,23.54) .. controls (366.07,23.65) and (362.55,33.24) .. (362.83,44.96) .. controls (363.1,56.69) and (367.06,66.11) .. (371.67,66) -- cycle ;
\draw    (371.67,66) .. controls (434.9,82.94) and (599.42,38.46) .. (592.9,77.95) ;
\draw    (504.74,108.15) .. controls (556.68,107.3) and (580,95.8) .. (592.9,77.95) ;
\draw    (506.75,135.71) .. controls (571.8,140.41) and (639.63,88.02) .. (652.1,79.22) ;
\draw    (370.68,23.54) .. controls (444.66,31.22) and (614.01,31.52) .. (650.27,50.27) ;
\draw   (652.42,79.19) .. controls (654.58,79.03) and (655.85,72.43) .. (655.25,64.44) .. controls (654.66,56.46) and (652.43,50.11) .. (650.27,50.27) .. controls (648.11,50.43) and (646.84,57.04) .. (647.44,65.02) .. controls (648.03,73.01) and (650.26,79.35) .. (652.42,79.19) -- cycle ;

\draw (312.08,89.26) node [anchor=north west][inner sep=0.75pt]  [font=\huge,rotate=-359.43] [align=left] {$\displaystyle \cong $};
\draw (12.08,27.28) node [anchor=north west][inner sep=0.75pt]   [align=left] {a};
\draw (13.41,80.63) node [anchor=north west][inner sep=0.75pt]   [align=left] {b};
\draw (14,144) node [anchor=north west][inner sep=0.75pt]   [align=left] {c};
\draw (365,36) node [anchor=north west][inner sep=0.75pt]   [align=left] {a};
\draw (369,97) node [anchor=north west][inner sep=0.75pt]   [align=left] {b};
\draw (373,148) node [anchor=north west][inner sep=0.75pt]   [align=left] {c};
\draw (137,24) node [anchor=north west][inner sep=0.75pt]   [align=left] {$\displaystyle \phi ( a,b)$};
\draw (250,63) node [anchor=north west][inner sep=0.75pt]   [align=left] {$\displaystyle \phi ( \phi ( a,b) ,c)$};
\draw (498,142) node [anchor=north west][inner sep=0.75pt]   [align=left] {$\displaystyle \phi ( b,c)$};
\draw (644,84) node [anchor=north west][inner sep=0.75pt]   [align=left] {$\displaystyle \phi ( a,\phi ( b,c))$};
\end{tikzpicture}
    \caption{Associativity}
    \label{fig:Associativity}
\end{figure}
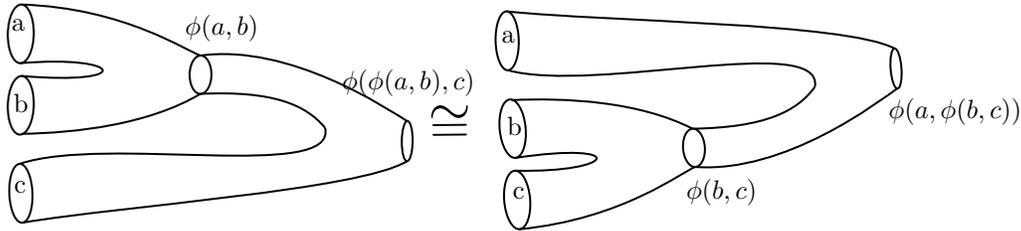

\begin{figure}[]
    \centering
        \label{fig:disk}

\begin{tikzpicture}[x=0.75pt,y=0.75pt,yscale=-0.3,xscale=0.4]

\draw   (135.74,39.26) .. controls (153.95,38.35) and (170.16,84.24) .. (171.95,141.77) .. controls (173.74,199.3) and (160.43,246.68) .. (142.22,247.59) .. controls (124.01,248.5) and (107.8,202.61) .. (106.01,145.08) .. controls (104.23,87.55) and (117.54,40.17) .. (135.74,39.26) -- cycle ;
\draw    (135.75,39.26) .. controls (639.01,52.03) and (639.01,245.93) .. (142.21,247.59) ;

\draw (62.98,137.22) node [anchor=north west][inner sep=0.75pt]   [align=left] {$\displaystyle S^{1}$};
\draw (546.62,129.41) node [anchor=north west][inner sep=0.75pt]   [align=left] {$\displaystyle \emptyset $};

\end{tikzpicture}

    \caption{Disk as a morphism from $S^1$ to $\emptyset$}
    \vspace{128in}

\end{figure}
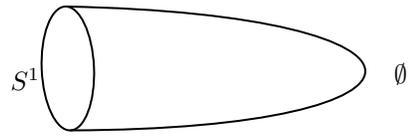
\newpage

\end{appendices}
\printbibliography
\end{document}